\documentclass[12pt]{article}
\usepackage{amsfonts}

\usepackage{latexsym}
\usepackage{amssymb}
\usepackage{epsfig}
\usepackage{amsmath}
\usepackage{amsthm}
\usepackage[mathscr]{eucal}
\usepackage{graphicx}

\usepackage{caption}
\usepackage{subcaption}

\newtheorem{Lemma}{Lemma}

\newtheorem{Theorem}[Lemma]{Theorem}

\newtheorem{Definition}{Definition}

\renewcommand{\qed}{\hfill{\ \ \rule{2mm}{2mm}} \vspace{0.2in}}

\begin{document}

\title{Duality in percolation via outermost boundaries III: Plus connected components}
\author{ \textbf{Ghurumuruhan Ganesan}
\thanks{E-Mail: \texttt{gganesan82@gmail.com} } \\
%EndAName
\ \\
New York University, Abu Dhabi}
\date{}
\maketitle

\begin{abstract}

%In this paper, we study the structure of left-right crossings
%of the random geometric graph \(G = G(n,r_n)\) of \(n\) nodes
%uniformly distributed in \(S = [0,1]^2\) with \(r_n = \epsilon\sqrt{\frac{\log{n}}{n}}\)
%for some \(\epsilon > 0.\) Tiling \(S\) horizontally and
%vertically into rectangles of length \(1\) and width \(Mr_n,\) we
%show that each rectangle has a left-right crossing of edges with
%high probability if \(M\) is sufficiently large.
%We call the resulting subgraph to be a ``backbone" of \(G.\)

%The techniques we use to construct the backbone has quite a few applications.
%As a first, we show that the diameter of second largest component in \(G\)
%is \(O(1)\) with high probability. Secondly,

Tile~\(\mathbb{R}^2\) into disjoint unit squares \(\{S_k\}_{k \geq 0}\) with the origin being the centre of \(S_0\)
and say that \(S_i\) and \(S_j\) are star adjacent if they share a corner and plus adjacent if they share an edge.
Every square is either vacant or occupied. In this paper, we use the structure of the outermost boundaries derived in
Ganesan~(2015) to alternately obtain duality between star and plus connected components in the following sense:
There is a star connected cycle of vacant squares attached to and surrounding
the finite plus connected component containing the origin. %We also obtain the mutual exclusivity of

\vspace{0.1in} \noindent \textbf{Key words:} Star and plus connected components, duality.

\vspace{0.1in} \noindent \textbf{AMS 2000 Subject Classification:} Primary:
60J10, 60K35; Secondary: 60C05, 62E10, 90B15, 91D30.
\end{abstract}

\bigskip

\renewcommand{\theequation}{\thesection.\arabic{equation}}
\setcounter{equation}{0}
\section{Introduction} \label{intro}

Tile \(\mathbb{R}^2\) into disjoint unit squares \(\{S_k\}_{k \geq 0}\) with origin being the centre of~\(S_0.\) Every square in \(\{S_k\}\) is assigned one of the two states, occupied or vacant and the square \(S_0\) containing the origin is always occupied. For \(i \neq j,\) we say that~\(S_i\) and~\(S_j\) are \emph{adjacent} or \emph{star adjacent} if they share a corner between them. We say that \(S_i\) and \(S_j\) are \emph{plus adjacent}, if they share an edge between them. Here we follow the notation of Penrose (2003).

\subsection*{Model Description}
We first discuss star connected components. We say that the square~\(S_i\) is connected to the square~\(S_j\) by a \emph{star connected \(S-\)path} if there is a sequence of distinct squares \((Y_1,Y_2,...,Y_t), Y_l \subset \{S_k\}, 1 \leq l \leq t\) such that~\(Y_l\) is star adjacent to~\(Y_{l+1}\) for all \(1 \leq l \leq t-1\) and \(Y_1 = S_i\) and \(Y_t = S_j.\) If all the squares in \(\{Y_l\}_{1 \leq l \leq t}\) are occupied, we say that~\(S_i\) is connected to~\(S_j\) by an \emph{occupied} star connected \(S-\)path.

Let \(C(0)\) be the collection of all occupied squares in \(\{S_k\}\) each of which is connected to the square~\(S_0\) by an occupied star connected \(S-\)path. We say that \(C(0)\) is the  star connected occupied component containing the origin. Throughout we assume that \(C(0)\) is finite and we study the outermost boundary for finite star connected components containing the origin. By translation, the results hold for arbitrary finite star connected components.

%//COLLECT ALL DEFNS OF STAR PATHS ETC TOGETHER!!!

%We use a piecewise cycle merging algorithm to the derive the
%structure of outermost boundaries for finite occupied star and plus connected components:
%In the plus connected case, the outermost boundary is a single cycle;
%and in the star connected case, we prove that the outermost boundary is a connected union
%of cycles with mutually disjoint
%interiors. Using the above, we then obtain duality
%between star and plus connected components.

Let \(G_0\) be the graph with vertex set being the set of all corners of the squares of \(\{S_k\}\) in the component~\(C(0)\) and edge set consisting of the edges of the squares of \(\{S_k\}\) in~\(C(0).\) Two vertices in the graph~\(G_0\) are said to be adjacent if they share an edge between them. Two edges in \(G_0\) are said to be adjacent if they share an endvertex between them.

Let \(P = (e_1,e_2,\ldots,e_t)\) be a sequence of distinct edges in \(G_0.\) We say that~\(P\) is a \emph{path} if~\(e_i\) and~\(e_{i+1}\) are adjacent for every \(1 \leq i \leq t-1.\)  Let~\(a\) be the endvertex of~\(e_1\) not common to~\(e_2\) and let~\(b\) be the endvertex of~\(e_t\) not common to~\(e_{t-1}.\) The vertices \(a\) and \(b\) are the \emph{endvertices} of the path~\(P.\)

We say that~\(P\) is a \emph{self avoiding path} if the following three statements hold: The edge \(e_1\) is adjacent only to \(e_2\) and no other~\(e_j, j \neq 2.\) The edge~\(e_t\) is adjacent only to~\(e_{t-1}\) and no other~\(e_j, j \neq t-1.\) For each~\(1 \leq i \leq t-1,\) the edge~\(e_i\) shares one endvertex with~\(e_{i-1}\) and another endvertex with~\(e_{i+1}\) and is not adjacent to any other edge~\(e_j, j\neq i-1,i+1.\)

We say that~\(P\) is a \emph{circuit} if \((e_1,e_2,\ldots,e_{t-1})\) forms a path and the edge~\(e_t\) shares one endvertex with \(e_1\) and another endvertex with~\(e_{t-1}.\) We say that~\(P\) is a \emph{cycle} if \((e_1,e_2,...,e_{t-1})\) is a self avoiding path and the edge~\(e_t\) shares one endvertex with~\(e_1\) and another endvertex with~\(e_{t-1}\) and does not share an endvertex with any other edge~\(e_j, 2 \leq j \leq t-2.\) Any cycle~\(C\) contains at least four edges and divides the plane~\(\mathbb{R}^2\) into two disjoint connected regions. As in Bollob\'as and Riordan~(2006), we denote the bounded region to be the \emph{interior} of~\(C\) and the unbounded region to be the \emph{exterior} of~\(C.\)

We use cycles to define the outermost boundary of star connected components. Let \(e\) be an edge in the graph~\(G_0\) defined above. We say that \(e\) is adjacent to a square \(S_k\) if it is one of the edges of~\(S_k.\) We say that \(e\) is a  \emph{boundary edge} if it is adjacent to a vacant square and is also adjacent to an occupied square of the component~\(C(0).\) Let \(C\) be any cycle of edges in \(G_0.\) We say that the edge~\(e\) lies in the interior (exterior) of the cycle~\(C\) if both the squares in~\(\{S_j\}\) containing~\(e\) as an edge lie in the interior (exterior) of~\(C.\)

%We say that~\(e\) lies in the exterior of~\(C\) if both the squares in~\(\{S_j\}\) containing~\(

%// CAN PUSH THE BELOW TO NEXT SEC?? PRKVMM +eTC...

We have the following definition.
\begin{Definition} \label{out_def} We say that the edge \(e\) in the graph~\(G_0\) is an \emph{outermost boundary} edge of the component \(C(0)\) if the following holds true for every cycle \(C\) in~\(G_0:\) either \(e\) is an edge in \(C\) or \(e\) belongs to the exterior of \(C.\)

We define the outermost boundary \(\partial _0\) of \(C(0)\) to be the set of all outermost boundary edges of~\(G_0.\)
\end{Definition}

The following result combines Theorems~\(1\) and~\(2\) from Ganesan~(2017).
\begin{Theorem}\label{thm_out} Suppose \(C(0)\) is finite.
The outermost boundary \(\partial_0\) of \(C(0)\) is the union of a unique set of cycles \(C_1,C_2,\ldots,C_n\) in \(G_0\) with the following properties:\\
\((i)\) Every edge in \(\cup_{1 \leq i \leq n} C_i\) is an outermost boundary edge.\\
\((ii)\) The graph \(\cup_{1 \leq i \leq n}C_i\) is a connected subgraph of \(G_0.\)\\
\((iii)\) If \(i \neq j,\) the cycles \(C_i\) and \(C_j\) have disjoint interiors and have at most one vertex in common.\\
\((iv)\) Every occupied square \(J_k \in C(0)\) is contained in the interior of some cycle \(C_{j}.\)\\
\((v)\) If \(e \in C_{j}\) for some \(j,\) then \(e\) is a boundary edge adjacent to an occupied square of~\(C(0)\) contained in the
interior of~\(C_j\) and also adjacent to a vacant square lying in the exterior of all the cycles in~\(\partial_0.\)\\
If the component~\(C(0)\) is plus connected, then the outermost boundary is a single cycle satisfying all the above conditions.
\end{Theorem}
%Similarly in Theorem~2 of Ganesan (2015), we obtained that the outermost boundary for the finite plus connected component is a single cycle.

The following result is Lemma~\(6\) of Ganesan~(2017) and collects the necessary properties regarding outermost boundary cycles.
\begin{Lemma}\label{outer} Let~\(C(0) = \cup_{1 \leq k \leq M} J_k \subset \{S_j\}\) be finite. For every occupied square~\(J_k \in C(0), 1 \leq k \leq M,\) there exists a unique cycle~\(D_k\) in the graph~\(G_0\) satisfying the following properties.\\
\((a)\) The square \(J_k\) is contained in the interior of \(D_k.\)\\
\((b)\) Every edge in the cycle \(D_k\) is a boundary edge adjacent to one occupied square of \(C(0)\) in the interior and one vacant square in the exterior.\\
\((c)\) If \(C\) is any cycle in \(G_0\) that contains \(J_k\) in the interior, then every edge in~\(C\) either belongs to~\(D_k\) or is contained in the interior.\\
\((d)\) Every edge in \(D_k\) is an outermost boundary edge in the graph~\(G_0.\)\\
We denote \(D_k\) to be the \emph{outermost boundary cycle} containing the square \(J_k \in C(0).\)
\end{Lemma}

\subsection*{Duality}
To study duality between star and plus connected components, we first define plus connected~\(S-\)cycles. Let \(\{J_i\}_{1 \leq i \leq m}\) be a set of squares in the set~\(\{S_k\}.\) We say that the sequence \(L_J = (J_1,...,J_m)\) is a \emph{plus connected \(S-\)path} if for each \(1 \leq i \leq m-1,\) we have that the square~\(J_i\) is plus adjacent i.e., shares an edge, with the square~\(J_{i+1}.\) We say that \(L^+\) is a \emph{plus connected \(S-\)cycle} if \((J_1,\ldots,J_{m-1})\) is a plus connected \(S-\)path and in addition, the square~\(J_m\) is plus adjacent to both~\(J_{m-1}\) and~\(J_1.\)

\begin{figure}[tbp]
\centering
%\fbox{
\includegraphics[width=2.5in, trim= 110 290 200 210, clip=true]{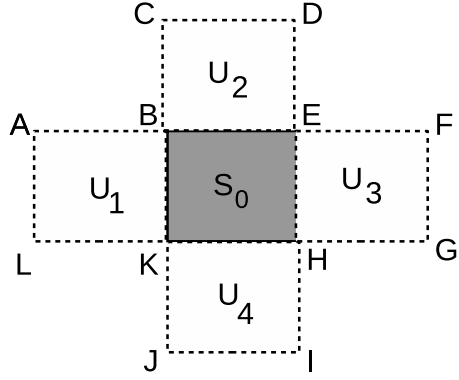}
%}
\caption{The sequence of vacant squares \((U_1,U_2,U_3,U_4)\) form the vacant star connected~\(S-\)cycle~\(H_{out}\) (see Theorem~\ref{thm5}) surrounding the plus connected component~\(C^+(0) = S_0.\)}
\label{skel_fig_ill}
\end{figure}

%We illustrate the result of Theorem~\ref{thm4} in Figure~\ref{skel_fig_ill} for the simplest case when \(C(0) = S_0,\) the square with origin as centre. The square~\(S_0\) is denoted as the dark grey square and the sequence of surrounding dotted squares form the \(S-\)cycle~\(G_{out}\) of vacant squares. The cycle~\(abcdefgha\) is the dual skeleton~\(SK(G_{out}).\) The dotted square formed by the corners~\(w,x,y\) and~\(z\) forms the outermost boundary~\(\partial_G\) of the \(S-\)cycle~\(G_{out}.\)

%From Theorem~\(1\) of Ganesan (2015), we have that the outermost boundary~\(\partial_0\) of the finite star connected component \(C(0)\) is a connected union of cycles with disjoint interiors. Moreover, any two cycles in~\(\partial_0\) share at most one vertex in common and every occupied square of \(C(0)\) is contained in the interior of some cycle of~\(\partial_0.\) In Figure~\ref{skel_fig_ill}\((b),\) we represent the cycles of \(\partial_0\) by the grey cycles marked~\(A\) and~\(B.\) All the occupied squares of the component \(C(0)\) are present inside one of the grey circles \(A\) or \(B.\) The dotted curve surrounding \(A\cup B\) represents the dual skeleton~\(SK(G_{out}).\) Finally, the continuous curve surrounding~\(SK(G_{out})\) is the outermost boundary cycle~\(\partial_G.\)

%\begin{figure}[tbp]
%\centering
%\fbox{
%\includegraphics[width=2in, trim= 100 280 180 50, clip=true]{skel_c0.eps}
%}
%\caption{Dual skeleton of an~\(S-\)cycle of squares represented by the dotted cycle~\(abcdefgha.\)}
%\label{skel_fig_ill}
%\end{figure}

We have an analogous definition for star connected \(S-\)cycle. As before, let \(\{Q_i\}_{1 \leq i \leq n}\) be a set of distinct squares in the set \(\{S_k\}.\) We say that the sequence \(L_Q = (Q_1,...,Q_n)\) is a \emph{star connected \(S-\)path} if for each \(1 \leq i \leq n-1,\) we have that~\(Q_i\) is star adjacent, i.e., shares a corner with~\(Q_{i+1}.\) We say that~\(L_Q\) is a \emph{star connected \(S-\)cycle} if \(n \geq 3,\) the sequence~\((Q_1,\ldots,Q_{n-1})\) forms a star connected \(S-\)path and in addition~\(Q_n\) is also star adjacent to~\(Q_1.\) We consider only star \(S-\)cycles containing at least three squares.

Let \(\Lambda^+_0\) denote the set of all vacant squares that is plus adjacent, i.e. shares an edge, with some occupied square in the plus connected component~\(C^+(0).\) We have the following result.
\begin{Theorem}\label{thm5} Suppose \(C^+(0)\) is finite. There exists a unique star connected \(S-\)cycle \(H_{out} = (U_1,...,U_q)\) with the following properties:\\
\((i)\) For every \(i, 1 \leq i \leq q,\) the square~\(U_i\) is in~\(\Lambda^+_0.\)\\
\((ii)\) The outermost boundary of~\(H_{out}\) is a single cycle~\(\partial_H.\) Every square in~\(H_{out}\) has an edge in~\(\partial_H.\)\\
\((iii)\) Every occupied square in the component~\(C^+(0)\) is contained in the interior of~\(\partial_H.\) Every vacant square in~\(\Lambda^+_0\) is contained in the interior of~\(\partial_H.\)
\end{Theorem}
Thus the corresponding sequence of squares~\((U_{1},...,U_{q})\) form a star connected cycle of vacant squares containing all squares of~\(C(0)\) in the interior. In Figure~\ref{skel_fig_ill}, the sequence of dotted squares~\(U_1,U_2,U_3\) and~\(U_4\) form the star connected sequence~\(H_{out}\) of vacant squares surrounding the plus connected component~\(C^+(0) = S_0,\) denoted by the dark grey square. The sequence of points~\(ABCDEFGHIJKLA\) forms the outermost boundary~\(\partial_H.\)

\subsection*{Auxiliary Results}
In this subsection, we state three results that are of independent interest and are used for the proof of Theorem~\ref{thm5}.
The results are intuitive and for completeness we provide the proofs in Section~\ref{misc} and~\ref{pf_thm8}.

Let \(G\) be the graph with vertex set as the vertices of corners of the squares in \(\{S_k\}\) and edge set as the edges of \(\{S_k\}.\) In the first result, we determine the outermost boundary of a star connected \(S-\)cycle (we refer to previous subsection for definition). We recall that any star connected \(S-\)cycle is a star connected component containing at least three distinct squares. From Theorem~\(1\) of Ganesan~(2017) we have that the outermost boundary of any star connected component is a connected union of cycles.  For the particular case of star \(S-\)cycles, we have the following slightly stronger result.
\begin{Theorem}\label{thm6} Let \(L_H = (H_1,\ldots,H_n), \{H_i\} \subseteq \{S_k\}, n \geq 3\) be a star connected \(S-\)cycle containing at least three squares. The outermost boundary~\(\partial_0(L_H)\) of~\(L_H\) is a single cycle. %in the graph \(G.\)
\end{Theorem}
For a plus connected \(S-\)cycle, the above result is true directly by Theorem~\(2\) of Ganesan~(2017).

The second result is regarding the union of the squares contained in a cycle. Consider a cycle \(C\) in the graph~\(G\) and suppose that \(\{Y_i\}_{1 \leq i \leq k} \subseteq \{S_j\}\) are the squares contained in the interior of \(C.\)
\begin{Theorem}\label{thm7}
The union of the squares \(\cup_i Y_i\) contained in the interior of a cycle~\(C \subset G\) is a plus connected component whose outermost boundary is~\(C.\)
\end{Theorem}
%//done...TO WRIT ABT OUTERMOST BNDRY THERE ETC...

The final result concerns the merging of two cycles with mutually disjoint interiors. Suppose that \(C\) and \(D\) are two cycles with mutually disjoint interiors that share more than one vertex in common. From Theorem~3 of Ganesan (2017), there exists a unique cycle \(E,\) consisting only of edges of \(C\) and~\(D,\) that contains both \(C\) and \(D\) in  its interior. We recall that the merging algorithm in Theorem~3 of Ganesan (2017) proceeds iteratively using paths of \(C\) that lie in the exterior of \(D.\)

\begin{figure}[tbp]
\centering
%\fbox{
\includegraphics[width=2.5in, trim= 60 300 140 120, clip=true]{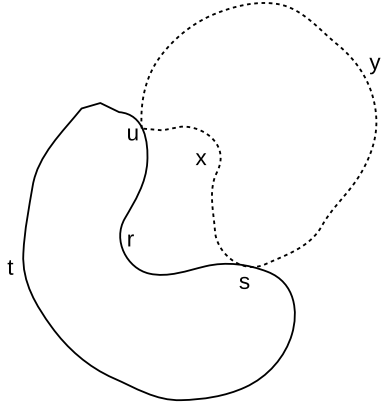}
%}
\caption{The segments~\(uxs\) and~\(uys\) of the dotted cycle~\(C = uxsyu\) are bridges for the solid cycle~\(D = urstu.\) The base for both the bridges is the segment~\(urs.\)}
\label{cyc_fig}
\end{figure}

The above merging algorithm is general and holds even for cycles with intersecting interiors. In case the cycles \(C\) and \(D\) have mutually disjoint interiors, we have a slightly stronger result. Suppose that \(P\) is a path of edges contained in the cycle \(C.\) We say that \(P \subset C\) is a \emph{bridge} for \(D\) if its endvertices lie in \(D\) and every other vertex in \(P\) lies in the exterior of \(D.\) Similarly we define bridges for the cycle \(C.\) In Figure~\ref{cyc_fig}, the segments~\(uxs\) and~\(uys\) of the dotted cycle~\(C = uxsyu\) are bridges for the solid cycle~\(D = urstu.\)

We have the following result.
\begin{Theorem}\label{thm8} Suppose that cycles \(C\) and \(D\) are two cycles with mutually disjoint interiors that contain more than one vertex in common. There is a unique cycle~\(E\) satisfying the following properties.\\
\((i)\) The cycle \(E = P_1 \cup P_2\) where \(P_1 \subset D\) is a bridge for cycle~\(C\) and \(P_2 \subset C\) is the bridge for \(D.\)\\
\((ii)\) The bridges~\(P_1\) and \(P_2\) are unique and the interior of cycles~\(C\) and~\(D\) lies in the interior of~\(E.\)
\end{Theorem}
We define \(E = P_1 \cup P_2\) to be the \emph{bridge decomposition} of the cycle~\(E.\) In Figure~\ref{cyc_fig}, the union of the bridges~\(P_1 = uts\) and~\(P_2 = uys\) contains both the cycles~\(C = uysxu\) and~\(D = urstu\) in its interior and is the desired bridge decomposition.

The paper is organized as follows: We prove the auxiliary results of Theorem~\ref{thm6} in Section~\ref{misc}, Theorem~\ref{thm7} in Section~\ref{misc2} and~Theorem~\ref{thm8} in Section~\ref{pf_thm8}. Finally, we prove the Theorem~\ref{thm5} in~\ref{pf5}. %We first prove the main Theorem~\ref{thm4} in Section~\ref{pf4}. Using Theorem~\ref{thm4}, we prove the Theorem~\ref{thm_lr} regarding left right crossings in Section~\ref{pf_lr}. We then prove the auxiliary results of Theorems~\ref{thm6} and~\ref{thm7} in Section~\ref{misc} and we prove the auxiliary~Theorem~\ref{thm8} in Section~\ref{pf_thm8}. Finally, we prove the Theorem~\ref{thm5} in~\ref{pf5}.

%be the square in \(C(0)\) whose centre has the largest \(x\) and

%To see for circuit...

%now to see that every edge is outer most boundary....seems ok...and to write crflly...fr our bnft etc... also to obtain circuit...and write carefully that cycle graph is a tree etc...for our bnft etc..and to write uniqueness by construction the outermost boundary is unique?? by construction shud be unique...suppose there exists a different set of cycles ?? adn argue briefly?? to see dtls crflly..for our bnt etc........

%\begin{figure}[tbp]
%\centering
%\includegraphics[width=3.5in, trim= 0 250 0 275, clip=true]{ckt_mg.eps}
%\caption{Merging cycle \(ABCDA\) with the segment \(AEC\)}
%\label{merg_cyc}
%\end{figure}

%//DO WE STATE THE THEOREM 1 OF GAN 2015?? PRKVMM...since we seem to use it often??prkvmm...

\section{Proof of Theorem~\ref{thm6}}\label{misc}
We label the squares in the sequence~\(L_Y\) as with label~\(1\) and every other square with label~\(0.\) We then apply Theorem~1 of Ganesan~(2017) with ``occupied" replaced by label \(1\) and ``vacant" replaced by label \(0.\)

%//DONE...NEED TO USE THE FACT THAT THERE ARE THREE OR MOE SQUARES IN S-CYCLE...

We first see that there cannot be three or more cycles in the outermost boundary~\(\partial_0(L_Y).\) We eliminate the two cycle case later. Suppose~\(\partial_0(L_Y) = \cup_{i=1}^{q} C_i\) where each~\(C_i, 1 \leq i \leq q\) is a cycle and~\(q \geq 3.\) The cycles~\(\{C_i\}\) satisfy the conditions of Theorem~1 of Ganesan~(2017) and so they have mutually disjoint interiors and~\(C_i\) and~\(C_j\) intersect at at most one point for~\(i \neq j.\)

We recall the construction of the cycle graph~\(H_{cyc}\) of~\(\partial_0(L_Y).\)  The vertex set of~\(H_{cyc}\) is~\(\{1,2,\ldots,q\}\) and we draw an edge between~\(i\) and~\(j\) if the corresponding cycles~\(C_i\) and~\(C_j\) intersect. From~\((3.1)\) of Ganesan~(2017) we have that~\(H_{cyc}\) is acyclic.

We use the cycle graph \(H_{cyc}\) to arrive at a contradiction. Without loss of generality, we assume that the first square \(Y_{1}\) lies in the interior of cycle~\(C_1.\) Let \(j_1\) be the smallest index \(j \geq 1\) such that \(Y_j\) does not belong to the interior of \(C_1.\) Again without loss of generality, we assume that \(Y_{j_1}\) belongs to the interior of the cycle \(C_2.\)  Thus the cycles \(C_1\) and \(C_2\) meet at some point \(v \in \mathbb{R}^{2}\) and there is an edge between vertices \(1\) and \(2\) in \(H_{cyc}.\) Also, the squares \(Y_{j_1-1}\) and \(Y_{j_1}\) meet at the point~\(v\) and since no other cycle \(C_i, 3 \leq i \leq q\) meets~\(C_1\) or~\(C_2\) at~\(v,\) we have that no other square of~\(L_Y\) contains~\(v.\)

For the next step suppose that~\(j_2\) is the smallest index~\(j > j_1\) such that the square~\(Y_{j}\) does not belong to the interior of~\(C_2.\) There exists such an index~\(j_2\) since any \(S-\)cycle contains at least three distinct squares. The square \(Y_{j_2}\) cannot belong to the interior of \(C_1.\) Because if \(Y_{j_2}\) belongs to the interior of \(C_1,\) then \(Y_{j_2-1}\) and \(Y_{j_2}\) necessarily meet at~\(v.\) Since \(j_1 \leq j_2-1 < j_2\) there are at least three distinct vacant squares in the set \(\{Y_{j_1-1},Y_{j_1},Y_{j_2-1},Y_{j_2}\}\) each containing the point~\(v.\) This is a contradiction to the last line of the previous paragraph.

%But this is a contradiction since from the above paragraph, the point \(v\) is only common to the vacant squares \(Y_{j_1-1}\) and \(Y_{j_1}\) and no other square of \(L_Y.\) //WRT MORE...and this is a contradiction (see Figure...) since three distinct vacant squares cannot meet at a single point. In Figure.....

From the above paragraph, we therefore have that the square~\(Y_{j_2}\) does not belong to the interior of the cycle~\(C_1\) and does not belong to the interior of the cycle~\( C_2.\) As before, without loss of generality, we assume that~\(Y_{j_2}\) is contained in the interior of the cycle~\(C_3.\) Correspondingly, the cycle graph~\(H_{cyc}\) contains an edge between the vertices~\(2\) and~\(3.\) In the next step of the iteration, suppose that~\(j_3\) is the smallest index~\(j > j_2\) such that~\(Y_{j}\) does not belong to the interior of~\(C_3.\)

Suppose that \(Y_{j_3}\) belongs to the interior of the cycle \(C_{k_4}.\) Arguing as above, we again have that three vacant squares cannot meet at a point and so \(k_4 \neq 2.\) If \(k_4 =1\) we have a contradiction since this means that vertices \(1\) and \(3\) are connected by an edge and so \(H_{cyc}\) contains a cycle. Else without loss of generality we let \(k_4 = 4.\)

Continuing this way, we obtain cycles \(C_{k_j}, j \geq 5\) and at each step of the iteration we check if we obtain a cycle in the graph \(H_{cyc}.\) Suppose the iteration continues until the final square \(Y_{n}\) of the sequence \(L_Y\) and suppose without loss of generality that \(Y_{n}\) belongs to the final cycle~\(C_{q}.\) By our iterative process, we have that the path \(P = (1,2,\ldots,q)\) is contained in the cycle graph \(H_{cyc}.\) Since \(Y_n\) and \(Y_1\) are star adjacent, the cycles \(C_1\) and \(C_q\) share a corner. Therefore there is an edge between vertices \(1\) and \(q\) in \(H_{cyc}\) and so again we obtain a cycle in \(H_{cyc},\) a contradiction.

%//done...THEN SEE then see that there cannot be two cycles...\\

From the above discussion we have that the outermost boundary \(\partial_0(L_Y)\) cannot have three or more cycles. Suppose now that \(\partial_0(L_Y) = C_1 \cup C_2\) contains exactly two cycles \(C_1\) and \(C_2.\) Suppose \(C_1\) and \(C_2\) meet at \(v \in \mathbb{R}^2.\) From Theorem~1 of Ganesan (2017), the cycles \(C_1\) and \(C_2\) do not have any other vertex in common. As before, we assume that the first square \(Y_{1}\) lies in the interior of cycle~\(C_1.\) Let \(j_1\) be the smallest index \(j > 1\) such that \(Y_j\) does not belong to the interior of \(C_1.\) Since there are two cycles, we have that \(j_1 < n\) strictly and the square \(Y_{j_1}\) belongs to the interior of the cycle~\(C_2.\)

Let \(j_2\) be the smallest index \(j > j_1\) such that \(Y_j\) does not belong to the interior of \(C_2.\) If no such \(j\) exists, then all the remaining squares \(\{Y_{j}\}_{j_1 \leq j \leq n}\) of the sequence \(L_Y\) belong to the interior of the cycle \(C_2.\) In particular, the square \(Y_{n}\) is star adjacent to \(Y_{1}.\) But this leads to a contradiction since \(Y_{1}\) shares a corner \(v \in \mathbb{R}^2\) with \(Y_{j_1}\) and is not star adjacent to any other square contained in the interior of \(C_2.\)

If on the other hand suppose that \(j_2 < n\) strictly. We then have that the square~\(Y_{j_2-1}\) (contained in the interior of the cycle~\(C_2\)) is star adjacent to~\(Y_{j_2}\) (contained in the interior of the cycle~\(C_1\)). This means that \(Y_{j_2-1}\) meets \(Y_{j_2}\) at the point~\(v.\) Arguing as above, we have a contradiction since there are at least three distinct vacant squares in the set \(\{Y_{j_1-1},Y_{j_1},Y_{j_2-1},Y_{j_2}\}.\) Therefore the outermost boundary~\(\partial_0(L_Y)\) is a single cycle. \(\qed\)

\section{Proof of Theorem~\ref{thm7}}\label{misc2}
We prove by induction on the number \(n\) of squares contained in the interior of a cycle \(C.\)  Strictly speaking, we apply induction assumption only on cycles containing the origin in the interior. The statement is true for arbitrary cycles by translation.

The statement is true for cycles containing a single square in the interior. Suppose it is true for all cycles containing \(n\) squares in the interior. Let~\(C\) be any cycle containing \(n+1\) cycles in its interior. Label all squares in the interior of \(C\) with label \(1\) and all other squares with label \(0.\) Let \(\{J_i\}_{0 \leq i \leq n}, J_i \subset \{S_k\}\) be the squares contained in the interior of \(C.\) Let \(J_0\) be the topmost rightmost square contained in the interior of~\(C\) obtained as follows. If \(z_k = (x_k,y_k)\) denotes the centre of square \(J_k,\) then \( y_0 =  \max_{0 \leq k \leq n} y_k.\) Further if \({\cal S} \supset J_0\) is the set of squares in \(\{J_i\}\) whose centre has \(x-\)coordinate \(y_0,\) then \(x_0 = \max_{J_i \in {\cal S}} x_i.\)

Let \(h_0(t), h_0(r), h_0(l)\) and \(h_0(b)\) be the top, right, left and bottom edges of the square \(J_0.\) We have the following properties.\\
\((t1)\) An edge \(e\) belongs to the cycle~\(C\) if and only if it is adjacent to a square~\(J_k\) in the interior of \(C\) and a square \(Q_k \notin \{J_i\}\) in the exterior of \(C.\) \\
\((t2)\) If \((x_i,y_i)\) is the centre of the square~\(J_i\) in the interior of the cycle~\(C,\) then \(y_i \leq y_0.\) Also if \(y_i = y_0,\) then \(x_i \leq x_0.\)\\
\((t3)\) The edges \(h_0(t)\) and \(h_0(r)\) belong to the cycle~\(C\) and \((h_0(t),h_0(r))\) is a subpath in~\(C.\)\\

%The squares \(W_i, 5 \leq i \leq 8\) lie in the exterior of~\(C.\)\\
%\((t4)\) Let \(x_{tl}\) be the vertex common to the edges \(h_0(t)\) and \(h_0(l).\) If~\(h_{prev}\) is the edge in the cycle~\(C\) sharing the vertex~\(x_{tl}\) with~\(h_0(t),\) then \(h_{prev} \in \{h_0(l), h_1(t)\}.\)\\
%//PRKVMM WHT PROPTEIS NEEDED ETC...
%//WRT HERE...

\emph{Proof of \((t1)-(t3)\)}: To see property \((t1)\) is true, let \(e \in C\) be an edge common to two squares in \(\{S_k\}.\) Exactly one of the squares is contained in the interior of the cycle~\(C\) and the other is contained in the exterior of~\(C.\) This proves property~\((t1)\) and property~\((t2)\) is also true by definition. %The third statement follows from the first two statements by considering the centre of the squares \(\{W_i\}_{5 \leq i \leq 8}.\) For example, the square \(W_5\) with centre \((u_5,v_5)\) satisfies \(v_5 = y_0\) and \(u_5 = x_0+1.\) An analogous argument holds for the other squares.

We prove \((t3)\) by contradiction. Suppose that \(h_0(t)\) does not belong to the cycle~\(C.\) The centre of the square~\(J_0\) is \((x_0,y_0)\) and so the centre of the top edge~\(h_0(t)\) is \(\left(x_0,y_0 + \frac{1}{2}\right).\) Since \(J_0\) is contained in the interior of~\(C,\) the edge~\(h_0(t)\) is contained in the interior of \(C.\) Therefore some edge \(e\) in~\(C\) cuts the line \(x = x_0\) at some point \((x_0,y)\) where \( y \geq y_0 + \frac{3}{2}.\) The final statement is true since edges of squares in \(\{S_k\}\) intersect the line \(x = x_{0}\) at \(y_{0} + \frac{k}{2}\) for integer \(k \neq 0.\)

The edge \(e\) is adjacent to some square \(J_k\) in the interior of the cycle~\(C\) and the centre~\((x_k,y_k)\) of~\(J_k\) satisfies \(y_k \geq y_0+1.\) This contradicts the first statement of property~\((t2)\) above. An analogous analysis using the second statement of property~\((t2)\) proves that the edge~\(h_0(r)\) also belongs to~\(C.\) Since the edges~\(h_0(t)\) and~\(h_0(r)\) have a common endvertex, they are consecutive edges in \(C.\)~\(\qed\)

%To prove \((t4),\) we only need to see that the ``top" edge \(h_8(r)\) sharing an endvertex with \(h_0(t)\) cannot belong to \(C.\) If it does, one of the two squares \(W_7\) or \(W_8\) containing \(h_8(r)\) necessarily belongs to the interior of \(C,\) a contradiction to property \((t2).\)~\(\qed\)

To apply the induction assumption, we modify the cycle \(C\) and obtain a cycle \(C_1\) containing exactly \(n\) squares in the interior. To do this, we use property \((t3)\) above to determine that at least two edges (\(h_0(t)\) and \(h_0(r)\)) of the square~\(J_0\) belong to~\(C.\) All four edges of~\(J_0\) cannot belong to~\(C\) since this would mean that~\(C\) contains only~\(J_0\) in its interior.

We therefore consider two possible cases: \((a)\) Exactly three edges of the square~\(J_0\) belong to the cycle~\(C\) and \((b)\) Exactly two edges of~\(J_0\) belong to~\(C.\)\\\\
\underline{Case \((a)\)}: We suppose that the edges~\(h_0(l), h_0(t)\) and~\(h_0(r)\) belong to the cycle~\(C.\) An analogous analysis holds if~\(h_0(t),h_0(r)\) and~\(h_0(b)\) belong to~\(C.\)

We modify the cycle \(C\) and define
\[C_1 = (C \setminus \{h_0(l),h_0(t),h_0(r)\}) \cup h_0(b).\] We have the following properties.\\
\((t4)\) The graph \(C_1\) is a cycle and the top most right most square \(J_0\) lies in the exterior of \(C_1.\) \\
\((t5)\) The cycle \(C_1\) contains exactly exactly \(n\) squares \(\{J_i\}_{1 \leq i \leq n}\) in its interior.

\emph{Proof of \((t4)-(t5)\)}: We prove that \(C_1\) is a cycle as follows. From property~\((t3)\) above, we have that \((h_0(t), h_0(r))\) is a subpath of the cycle~\(C.\) Since~\(h_0(l)\) and~\(h_0(t)\) share a common endvertex, we have that the sequence of edges~\((h_0(l),h_0(t),h_0(r))\) form a subpath in~\(C.\) Therefore the bottom edge~\(h_0(b)\) of the square \(J_0\) cannot belong to the cycle \(C\) since this would mean that \(C\) consists only one square~\(J_0\) in its interior, a contradiction.

Let \(x\) and \(y\) be the endvertices of the bottom edge \(h_0(b).\) The subpath of edges\\\((h_0(l),h_0(t),h_0(r))\) also has endvertices~\(x\) and~\(y\) and therefore the set\\\(C \setminus \{h_0(l),h_0(t),h_0(r)\}\) is a path of edges with endvertices~\(x\) and~\(y.\) Since \(h_0(b) \notin C\) (see previous paragraph), we have that~\(C_1\) as defined above is a cycle.

To see that the square~\(J_0\) lies in the exterior of the cycle~\(C_1,\) let \(z_{0} = (x_{0}, y_{0})\) be the centre of the square~\(J_0.\) If \(J_0\) lies in the interior of \(C_1,\) then some edge \(e\) in the cycle~\(C_1\) intersects the line \(x = x_{0}\) at \((x_{0},y)\) for some \(y \geq y_{0} + \frac{1}{2}.\) Since the edge \(h_0(t) \notin C_1,\) we have that \(e \neq h_0(t)\) and so \(y \geq y_0 + \frac{3}{2}.\)

The edge \(e\) also necessarily  belongs to the original cycle~\(C\) and is therefore adjacent to some square \(J_k\) contained in the interior of \(C.\) The centre \((x_k,y_k)\) of \(J_k\) satisfies \(y_k \geq y_0 + 1.\) This contradicts property~\((t2)\) above. Thus \(J_0\) lies in the exterior of \(C_1\) and this proves \((t4).\)

We prove \((t5)\) as follows. The square \(J_0\) shares only the edge \(h_0(b)\) with the cycle~\(C_1\) and no other endvertex with \(C_1.\) Thus every point in the interior of \(C\) either belongs to \(C_1\) or the square \(J_0.\) This means that \(C_1\) has exactly~\(n\) squares \(\{J_i\}_{1 \leq i \leq n}\) in its interior.\(\qed\)

We now apply the induction assumption on the cycle \(C_1.\) The set of squares contained in \(C_1\) form a plus connected component and so the squares \(\{J_i\}_{1 \leq i \leq n}\) form a plus connected component. The square~\(J_0\) in the exterior of~\(C_1\) shares the edge~\(h_0(b)\) with some square~\(J_k, k \geq 1\) lying in the interior of~\(C_1.\) This means that \(J_0 \cup \{J_i\}_{1 \leq i \leq n}\) form a plus connected component. This proves the theorem for the case where exactly three edges of the square~\(J_0\) belong to the cycle~\(C.\)\\\\
\underline{Case \((b)\)}: We now consider the other possibility where there are exactly two edges of the square~\(J_0\) in the cycle~\(C.\)
Using property \((t3)\) above, we have that the only edges of \(J_0\) belonging to \(C\) are \(h_0(t)\) and \(h_0(r).\)

We modify the cycle \(C\) and define
\[C_2 = (C \setminus \{h_0(t),h_0(r)\}) \cup \{h_0(l),h_0(b)\}.\] Let \(z_{lb}\) be the endvertex common to \(h_0(l)\) and \(h_0(b).\)
We consider two subcases separately depending on whether \(z_{lb} \in C\) are not.

\((b1)\) If the vertex~\(z_{lb}\) does not belong to the cycle~\(C,\) then the argument is analogous to Case~\((a)\) above. We have that properties \((t4)-(t5)\) above are true for the graph~\(C_2\) also. Therefore~\(C_2\) is a cycle and the squares \(\{J_i\}_{1 \leq i \leq n}\) in the interior of~\(C_2\) form a plus connected component. The square~\(J_0\) shares an edge~\(h_0(b)\) with~\(C_2\) and so the squares \(J_0 \cup \{J_i\}_{1 \leq i \leq n}\) form a plus connected component.

%//WRT MORE DETAILS +ETC..??prkvmm...

\((b2)\) If the vertex~\(z_{lb}\) belongs to the cycle~\(C,\) then the modified graph~\(C_2\) is a union of two cycles~\(D_1\) and~\(D_2\) sharing a common vertex~\(z_{lb}.\) We apply induction assumption on each of~\(D_1\) and~\(D_2.\) Without loss of generality assume that squares \(\{J_i\}_{1 \leq i \leq p}\) belong to the interior of~\(D_1\) and the rest \(\{J_i\}_{p+1 \leq i \leq n}\) belong to the interior of~\(D_2.\) The square~\(J_0\) shares an edge with some square \(J_{i_1}, 1 \leq i_1 \leq p\) lying in the interior of \(D_1\) and another square \(J_{i_2}, p+1 \leq i_2 \leq n\) lying in the interior of~\(D_2.\)

By induction assumption the squares \(\{J_i\}_{1 \leq i \leq p}\) form a plus connected component. Thus \(J_0 \cup \{J_i\}_{1 \leq i \leq p}\) forms a plus connected component. Similarly, the squares \(J_0 \cup \{J_i\}_{p+1 \leq i \leq n}\) also form a plus connected component. Thus the squares~\(\{J_i\}_{0 \leq i \leq n}\) form a plus connected component. This proves case \((b).\)

To prove the statement regarding the outermost boundary, we argue as follows. From Theorem~\(2\) of Ganesan (2017), we have that the outermost boundary of the plus connected component~\(\{J_i\}_{0 \leq i \leq n}\) is a single cycle~\(C_J.\) Also, every square \(J_i, 0 \leq i \leq n\) belongs to the interior of \(C_J.\) Since the cycle~\(C\) contains only of edges of~\(\{J_i\}_{0 \leq i \leq n},\) we have that every edge in the cycle~\(C\) either belongs to or lies in the interior of~\(C_J.\) If \(C_J \neq C,\) then there exists an edge~\(e \in C_J,\) belonging to the exterior of~\(C.\) The edge~\(e\) is the edge of some square~\(J_k, 0 \leq k \leq n\) and so \(J_k\) lies in the exterior of~\(C,\) a contradiction.\(\qed\)

\section{Proof of Theorem~\ref{thm8}}\label{pf_thm8}
Let \(\{P_i\}_{1 \leq i \leq t} \subset C\) be the set of all bridges for the cycle~\(D,\) contained in the cycle~\(C.\) Fix \(1 \leq i \leq t\) and suppose that~\(a_i\) and~\(b_i\) are the endvertices of the path~\(P_i\) and let~\(Q_i\) and~\(R_i\) be the two subpaths in the cycle~\(D\) with endvertices~\(a_i\) and~\(b_i.\) Among the two cycles~\(Q_i \cup P_i\) and~\(R_i \cup P_i,\)  exactly one of them contains no squares of~\(D\) in its interior. Suppose
\begin{equation}\label{gap_def}
g(P_i, D) = Q_i \cup P_i
\end{equation}
contains no squares of~\(D\) in its interior. We define~\(g(P_i,D)\) to be the \emph{gap} between the bridge~\(P_i\) and the cycle~\(D.\) We also define~\(Q_i = Ba(P_i,D)\) to be the \emph{base} for the bridge~\(P_i\) contained in the cycle~\(D.\) By definition, the union of the paths \(P_i \cup R_i\) contains both the cycles~\(g(P_i,D)\) and~\(D\) in its interior.

In Figure~\ref{cyc_fig}, the cycles~\(C\) and~\(D\) are denoted by the dotted and solid curves, \(uxsyu\) and \(urstu,\) respectively. The segments~\(P_1 = uxs\) and \(P_2 = uys\) are bridges with the same base \(Q_1 = Q_2 = urs.\) The gap \(g(P_1,D)\) is the union of the segments~\(uxs\) and~\(urs.\) The paths~\(R_1 = R_2 = uts\) and union of the paths~\(P_1 \cup R_1 = \{uxs\} \cup \{uts\}\) contains both the cycles~\(g(P_1,D) = uxsru\) and~\(D = urstu\) in its interior.

Let \(\{J_i\}_{1 \leq i \leq n} \subset \{S_k\}\) be the squares in the interior of~\(C.\) We first prove the result of the Theorem when~\(t = 1,\) i.e., there is exactly one bridge,~\(P_1.\) In this case, the desired bridge decomposition is \(P_1 \cup R_1\) and we argue as follows. The interior of the cycle~\(P_1 \cup R_1\) contains the interior of  the cycle~\(D\) in its interior. If there exists a square~\(J_k\) that lies in the interior of \(C\) but in the exterior of \(P_1 \cup R_1,\) then there exists an edge \(e \notin P_1\) of the cycle~\(C\) lying in the exterior if \(P_1 \cup R_1.\) Since the interior of the cycle~\(D\) is contained in the interior of~\(P_1 \cup R_1,\) the edge~\(e\) also lies in the exterior of~\(D.\) But this is a contradiction, since the path~\(P_1\) contains all edges of the cycle~\(C,\) lying in the exterior of the cycle~\(D.\)

%We have the following properties. \\
We henceforth assume that \(t \geq 2\) and therefore there are at least two distinct bridges and two (not necessarily distinct) bases.
We have the following properties.\\
\((w1)\) If edge \(e \in C\) lies in the exterior of the cycle \(D,\) then \(e \in \cup_{i=1}^{t} P_i.\) If \(i \neq j,\) the path~\(P_i\) has no edges in common with the cycle~\(g(P_j,D).\)\\
\((w2)\) For every \(1 \leq i \leq t,\) the gap \(g(P_i,D)\) and the cycle \(D\) have mutually disjoint interiors. Either all the squares \(\{J_k\}_{1 \leq k \leq n}\) contained in the interior of the cycle \(C\) also lie in the interior of the cycle \(g(P_i,D)\) or all the squares \(\{J_k\}_{1 \leq k \leq n}\) lie in the exterior of \(g(P_i,D).\)\\
\emph{Proof of \((w1)-(w2)\)}: The first part property~\((w1)\) and the first part of property~\((w2)\) is true by construction. We prove the second part of property~\((w1)\) as follows. We fix \(i \neq j\) and show separately that \(P_j\) is edge disjoint with~\(P_i\) and with~\(Q_i.\) Let \(P_i = (e_{1},\ldots,e_{k})\) and suppose \(e_{1} \in P_i \cap P_j.\) Traversing the path~\(P_i\) let~\(f\) be the smallest index such that the edges \(e_{l}, 1 \leq l \leq f\) belong to \(P_i \cap P_j\) and the edge~\(e_{f+1} \in P_i \setminus P_j.\)

The edges~\(e_f\) and~\(e_{f+1}\) share an endvertex \(v_f.\) Moreover, there is also an edge \(g \in P_j\) with endvertex \(v_f\) and \(g \notin \{e_f,e_{f+1}\}.\) This is a contradiction since the edges \(e_{f},e_{f+1}\) and \(g\) all belong to the cycle~\(C\) and every vertex in~\(C\) is adjacent to exactly two edges of~\(C.\) This proves that~\(P_i\) and~\(P_j\) are edge disjoint. The path \(Q_i\) is a subpath of the cycle~\(D\) and every edge in the path~\(P_j,\) by definition, lies in the exterior of~\(D.\) So the paths~\(P_j\) and~\(Q_i\) are also disjoint. This proves the second part of \((w1).\)

We prove the second part of \((w2)\) for \(i = 1\) and an analogous proof holds for all~\(i.\) Suppose that the gap \(g(P_1,D) = P_1 \cup Q_1\) where~\(Q_1\) is the base for the bridge \(P_1\) as defined in (\ref{gap_def}). From Theorem~\ref{thm7}, we have that the squares \(\{J_j\}_{1 \leq j \leq n}\) contained in the interior of the cycle~\(C\) form a plus connected component.

If there exists a square \(J_{i_1}\) of \(C\) in the interior of \(g(P_1,D)\) and a square~\(J_{i_2}\) in the exterior, then there is a plus connected path from~\(J_{i_1}\) to~\(J_{i_2}\) consisting only of squares in~\(\{J_i\}.\) In particular some edge \(e\) of the cycle~\(g(P_1,D)\) is common to a square~\(J_{k_1}\) contained in the interior of \(g(P_1,D)\) and a square~\(J_{k_2}\) in the exterior of \(g(P_1,D).\) If \(e \in P_1,\) this is a contradiction since \(P_1 \subset C\) and every edge in the cycle \(C\) is adjacent to one square in the interior of~\(C\) and one square in the exterior of~\(C.\)

If \(e \in Q_1 \subseteq D,\) then \(e\) is adjacent to one square in the interior of \(D\) and one square in the exterior of \(D.\) Thus one of the squares \(J_{k_1}\) or \(J_{k_2}\) lies in the interior of \(D.\) But this is a contradiction since \(C\) and \(D\) have mutually disjoint interiors and so no square \(J_m, 1 \leq m \leq n\) is contained in the interior of the cycle~\(g(P_1,D).\) This proves the second part of~\((w2).\) \(\qed\)

We have the following Lemma regarding the interior of the gaps.
\begin{Lemma}\label{gap_n_empty}
There exists \(1 \leq i_0 \leq t\) such that the gap~\(g(P_{i_0},D)\) contains the cycle~\(C\) in its interior.
\end{Lemma}
Our argument is by contradiction. Suppose every gap is empty, i.e., if \(1 \leq i \leq t,\) then the gap~\(g(P_i,D)\) contains no square of~\(\{J_k\}_{1 \leq k \leq n}\) in its interior. We recall that \(\{J_k\}_{1 \leq k \leq n}\) are the squares contained in the interior of the cycle~\(C.\) Consider two bridges~\(P_{i}\) and~\(P_{j}, i \neq j\) and the corresponding bases~\(Q_i\) and~\(Q_j.\) We have the following property.
\begin{eqnarray}
&&\text{Every edge of~\(P_j\) is contained in the exterior of the cycle~\(g(P_i,D).\)} \nonumber\\
&&\;\;\;\;\;\;\text{Every edge of~\(Q_j\) is contained in the exterior of~\(g(P_i,D).\)}\label{w3_prop}
\end{eqnarray}
In other words, the cycles \(g(P_i,D)\) and \(g(P_j,D)\) are edge disjoint and have mutually disjoint interiors.\\\\
\emph{Proof of (\ref{w3_prop})}: To see that the first statement is true, we assume that some edge~\(e\) in \(P_j\) is contained in the interior of the cycle~\(g(P_i,D) = P_i \cup Q_i.\) The edge~\(e\) is the edge of a square~\(J_k\) contained in the interior of the cycle~\(C.\) Thus the square \(J_k\) lies in the interior of the cycle~\(g(P_i,D)\) and from property~\((w2),\) we have that all squares \(\{J_k\}_{1 \leq k \leq n}\) are contained in the interior of the cycle~\(g(P_i,D),\) a contradiction.

We prove the second statement of~(\ref{w3_prop}) as follows. Suppose some edge~\(f\) in~\(Q_j\) lies in the interior of the cycle~\(g(P_i,D).\) The edge~\(f\) is the edge of a square~\(S_k\) lying in the interior of the cycle~\(D.\) The square \(S_k\) lies in the interior of the cycle \(g(P_i,D)\) which contradicts the construction of the gap defined above.

Thus every edge in \(Q_j\) either belongs to~\(g(P_i,D)\) or lies in the exterior of~\(g(P_i,D).\) In particular, the cycles \(g(P_i,D)\) and \(g(P_j,D)\) have mutually disjoint interiors. Suppose some edge \(e \in Q_j\) belongs to~\(g(P_i,D).\) The edge \(e \in Q_j\) is an edge of the cycle \(D\) and by definition, every edge in the path \(P_i\) lies in the exterior of \(D.\) Thus \(e \in Q_i\) and let \(A_1\) and \(A_2\) be the two squares containing \(e\) as an edge. Suppose \(A_1\) lies in the interior of \(g(P_i,D)\) and \(A_2\) lies in the interior of \(g(P_j,D).\) One of the squares in \(\{A_1,A_2\}\) lies in the interior of the cycle \(D\) and this contradicts the definition of the gaps in (\ref{gap_def}).

From the above discussion, we have that every edge in \(Q_j\) lies in the exterior of \(g(P_i,D).\) This proves (\ref{w3_prop}).\(\qed\)

%Indeed~\(R_i\) contains only edges of~\(D\) that do not belong to the cycle~\(g(P_i,D).\) From property \((w2),\) we have that \(g(P_i,D)\) and the cycle \(D\) have mutually disjoint interiors and so every edge of~\(R_i\) lies in the exterior of the cycle~\(C.\) If a vertex \(v \in C\) belongs to \(R_i,\) then \(v\) also belongs to \(Q_i.\) This means that \(v \in \{a_i,b_i\}\) is the endvertex of the path~\(Q_i\) common to~\(R_i.\)

\subsection*{Merging bridges with \(D\)}
To prove Lemma~\ref{gap_n_empty}, we first merge the bridges~\(\{P_k\}_{1 \leq k \leq t}\) with the cycle~\(D\) one by one to obtain a final cycle~\(D_{fin}.\) Since we have assumed that the gaps~\(\{g(P_i,D)\}\) are all empty, we obtain that the final cycle~\(D_{fin}\) does not contain the cycle~\(C\) in its interior. This contradicts Theorem~\(3\) of Ganesan~(2017), where we obtain that the cycles obtained by merging~\(C\) and~\(D\) is unique and contains both~\(C\) and~\(D\) in its interior.

At the beginning of the iteration, we set \(E_0 = D\) and in the first iteration we consider the bridge~\(P_1.\) We write \(P_1 = B(P_1,E_0)\) to emphasize the fact that~\(P_1\) is a bridge for the cycle~\(E_0.\) Every edge in the bridge~\(B(P_1,E_0)\) lies in the exterior of~\(E_0\) and the endvertices of~\(B(P_1,E_0)\) belong to the cycle~\(E_0.\) We set \(E_1 = P_1 \cup R_1\) to be the cycle obtained in the first iteration after merging~\(B(P_1,E_0)\) with the cycle~\(E_0.\) Here \(R_1 = E_0 \setminus Q_1\) is a subpath of \(E_0\) and the path \(Q_1 = Ba(P_1,E_0)\) is base for the bridge~\(P_1\) contained in the cycle~\(E_0.\)

For \(i = 1,\) the cycle \(E_i\) has the following properties. \\
\((f1)\) The cycle \(E_i = \left(\cup_{k=1}^{i} P_k\right) \cup T_i,\) where \(T_i = D\setminus \left(\cup_{k=1}^{i} Q_k\right).\)\\
\((f2)\) The interior of the cycle~\(E_i\) is the union of the interior of the cycle~\(E_{i-1}\) and the interior of the cycle~\(g(P_i,E_{i-1}) = g(P_i,E_0).\) Thus every square in the interior of \(E_i\) either belongs to the interior of \(D\) or belongs to the interior of one of the cycles in~\(\cup_{k=1}^{i}g(P_k,E_0).\)\\
\((f3)\) All squares in the interior of the cycle~\(C\) are contained in the exterior of~\(E_i.\)\\
\((f4)\) For \(1 \leq i \leq t-1,\) we have that \(P_{i+1} = B(P_{i+1},E_i)\) and \(Q_{i+1} = Ba(P_{i+1},E_i).\) Also every square in the interior of the cycle~\(g(P_{i+1},E_i) = g(P_{i+1},E_0)\) lies in the exterior of~\(E_i.\)\\\\
From the final property we obtain that the path~\(P_{i+1}\) is also a bridge for the new cycle~\(E_i,\) the path~\(Q_{i+1}\) is still the base for~\(P_{i+1}\) contained in the cycle~\(E_i\) and the gap between~\(P_{i+1}\) and~\(E_i\) remains unchanged and lies in the exterior of~\(E_i.\) This allows us to proceed to the next step of the iteration.

\emph{Proof of \((f1)-(f4)\) for \(i = 1\)}: The properties \((f1)\) and~\((f2)\) are true by construction. The property \((f3)\) is true by the definition of~\(E_i.\) %//WRT MORE HERE??? PRKVMM+eTC...

To prove \((f4),\) we first verify the following properties.\\
\((a)\) The subpath \(Q_{i+1}\) is contained in the cycle \(E_i.\)\\
\((b)\) Every edge in the path~\(P_{i+1}\) lies in the exterior of the cycle~\(E_i.\)\\
\((c)\) The cycle \(g(P_{i+1},E_0) = P_{i+1} \cup Q_{i+1}\) and the cycle~\(E_i\) have mutually disjoint interiors.

For \((a)\) we use the fact that~\(Q_j\) and~\(Q_k\) are edge disjoint (see~(\ref{gap_def}) and~(\ref{w3_prop})) for~\(j \neq k\) and so \(T_i = D \setminus \cup_{l=1}^{i} Q_l \subset E_i\) contains the subpath~\(Q_{i+1}.\) %To see that~\(Q_{i+1}\) is still the base for~\(P_{i+1}\) contained in the cycle \(E_i,\) we verify the following two statements:

We prove \((b)\) as follows. Using~\((f2),\) we obtain that the interior of the cycle~\(E_i\) is the union of the interiors of the cycle~\(D\) and the interiors of the cycles~\(g(P_k,D), 1 \leq k \leq i.\) By definition, every edge in the path~\(P_{i+1}\) lies in the exterior of the cycle~\(D.\) From property~(\ref{w3_prop}) of the gaps, every edge in~\(P_{i+1}\) lies in the exterior of the cycle~\(g(P_k,D)\) for~\(1 \leq k \leq i.\) Thus every edge in~\(P_{i+1}\) lies in the exterior of the cycle~\(E_i.\) This proves \((b).\)

For \((c),\) we first use property~\((f3)\) to obtain that every square in the interior of the cycle~\(E_i\) either belongs to the interior of one of the cycles in \(\cup_{k=1}^{i} g(P_k,E_0)\) or belongs to the interior of the cycle~\(E_0 = D.\) From~(\ref{w3_prop}), we have that the cycle~\(g(P_{i+1},E_0)\) has mutually disjoint interior with every cycle~\(g(P_k,E_0), 1 \leq k \leq i\) and from property~\((w2),\) we have that the cycle~\(g(P_{i+1},E_0)\) and~\(D\) have mutually disjoint interiors. Thus \((c)\) is also true.

We prove \((f4)\) now. The endvertices~\(a_{i+1}\) and~\(b_{i+1}\) of the path~\(Q_{i+1}\) are also the endvertices of~\(P_{i+1}.\) Using~\((a)\) and~\((b),\) we have that the path~\(P_{i+1}\) is also a bridge for the new cycle~\(E_i.\) Also using~\((c),\) we have that the path \(Q_{i+1} = Ba(P_{i+1},E_i)\) is the base for~\(P_{i+1}\) contained in the cycle~\(E_i.\) Statements~\((a)-(c)\) also obtain that the gap~\(g(P_{i+1},E_i) = g(P_{i+1},E_0).\)~\(\qed\)

%//WRT MORE + ETC???

%//WRT MORE +ETC ABV AND BELOW??PRKVMM..

Using property \((f4),\) we proceed to the next step of the iteration and merge \(P_2\) with \(E_1\) to get a new cycle \(E_2\) which again satisfies properties \((f1)-(f4).\) This procedure continues and after~\(n\) steps, we obtain the final cycle~\(E_n.\) We use the properties of the cycle \(E_n\) to prove Lemma~\ref{gap_n_empty}.

\emph{Proof of Lemma~\ref{gap_n_empty}}: From property \((w1),\) we have that every edge of the cycle \(C\) lying in the exterior of the cycle \(D\) belongs to one of the bridges in \(\{P_i\}_{1 \leq i \leq t}.\) Therefore the cycle~\(E_n\) is precisely the cycle obtained by the merging algorithm of Theorem~\(3\) of Ganesan~(2017). From Theorem~\(3\) of Ganesan (2017), we also have that the cycle~\(E_n\) contains the cycles~\(C\) and~\(D\) in its interior. But this contradicts property~\((f2),\) which states that every square in the interior of cycle~\(C\) lies in the exterior of the cycle~\(E_n.\) This proves Lemma~\ref{gap_n_empty}. \(\qed\)

\emph{Proof of Theorem~\ref{thm8}}: From Lemma~\ref{gap_n_empty}, we obtain that there exists \(1 \leq i_0 \leq t\) such that the gap \(g(P_{i_0},D) = P_{i_0} \cup Q_{i_0}\) (see (\ref{gap_def})) contains the cycle \(C\) in its interior. We recall that the path \(Q_{i_0}\) is a subset of the cycle \(D\) and the path \(R_{i_0} = D \setminus Q_{i_0}\) is such that the union of the paths~\(P_{i_0} \cup R_{i_0}\) is a cycle that contains both the cycles~\(g(P_{i_0},D)\) and~\(D\) in its interior. In particular, the union~\(P_{i_0} \cup R_{i_0}\) is a cycle containing both the cycles~\(C\) and~\(D\) in its interior. The path~\(P_{i_0} \subset C\) is a bridge for cycle~\(D.\) Reversing the roles of \(C\) and \(D\) in the above proof, we obtain that~\(R_{i_0}\subset D\) is also a bridge for cycle~\(C.\)~\(\qed\)

 %Either \(f_k\) lies in the left edge of \(R\) or at least one endvertex of \(g_j\) lies in the interior of \(R.\) The edge \(g_j\) is also the edge of a vacant square~\(A_j\) lying in the interior of the rectangle \(R.\)

\section{Proof of Theorem~\ref{thm5}}\label{pf5}
%(show uniqueness later...after obtaining the outermost boundary...)

%(?? write later..

%(see fig...) (to write more..?? see crflly..)

\emph{Proof of Theorem~\ref{thm5}}:
Let \(G^+_C\) denote the graph with vertex set as corners of squares of \(C^+(0)\) and edge set as edges of such squares. Let
\begin{equation}\label{do_plus}
\partial^+_0 = (e_1,...,e_t)
\end{equation}
be the outermost boundary cycle in \(G^+_C\) for the component \(C^+(0)\) obtained from Theorem~2 of Ganesan (2017).
Every edge \(e_i\) is adjacent to a vacant square~\(Y_i\) contained in the exterior of \(\partial_0^+\) and an occupied square~\(W_i\) contained in the interior of~\(\partial_0^+.\)

\begin{figure}[tbp]
\centering
%\fbox{
\includegraphics[width=3in, trim= 80 240 220 330, clip=true]{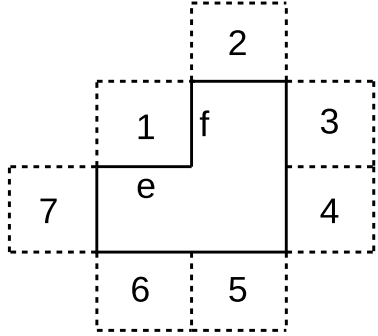}
%}
\caption{The solid curve denotes the outermost boundary cycle~\(\partial^+_0\) and the vacant squares labelled~\(i, 1 \leq i \leq 7,\) share an edge with~\(\partial_0^+\) and lie in the exterior of~\(\partial^+_0.\)}
\label{yj_nt_sm_fig}
\end{figure}

For \(i \neq j,\) the squares \(Y_i\) and \(Y_j\) need not be distinct. This is illustrated in Figure~\ref{yj_nt_sm_fig}, where the solid curve represents the outermost boundary cycle~\(\partial^+_0.\) The vacant squares labelled \(i, 1 \leq i \leq 7\) share an edge with~\(\partial_0^+\) and lie in the exterior of \(\partial^+_0.\) The edges~\(e\) and~\(f\) of the outermost boundary cycle~\(\partial_0^+\) are both edges of the same vacant square labelled~\(1.\)

Let \(\Lambda^+ = \bigcup_{i=1}^{t}\{Y_i\}\) be the set of all vacant squares lying in the exterior of \(\partial_0^+\) and sharing an edge with~\(\partial_0^+.\)

\subsubsection*{Outline of the Proof}
The proof consists of seven steps. In the \textbf{first} step, we merge vacant squares in \(\Lambda^+\) one by one with the outermost boundary \(\partial_0^+\) in an iterative manner to obtain a final cycle~\(D_{fin}.\) This cycle \(D_{fin}\) is in fact the desired outermost boundary \(\partial_H\) in the statement of the Theorem. To see this, we extract a star connected component \(L\) of vacant squares by exploring~\(D_{fin}\) edge by edge. We show that the outermost boundary of \(L\) is \(D_{fin}\) and \(D_{fin}\) contains at least one edge from \emph{every} square in~\(L.\) This forms the \textbf{second} step in the proof.

If we establish that \(L\) is a star connected \(S-\)cycle, then we are done since~\(L\) would be the desired~\(H_{out}.\) Indeed, in the \textbf{third} step of our proof, we see that the component \(L\) \emph{contains} a star connected \(S-\)cycle. In the remaining three steps, we assume that~\(L\) is not a star \(S-\)cycle and arrive at a contradiction. In the \textbf{fourth} step, we let \(L_1 \neq L\) be an \(S-\)cycle contained in~\(L\) and prove that the outermost boundary~\(\partial_0(L_1)\) of~\(L_1\) is a single cycle that lies in the exterior of the outermost boundary cycle~\(\partial_0^+\) and shares edges with~\(\partial_0^+.\) We also enumerate the relevant properties of~\(\partial_0(L_1)\) needed for future use.

Our aim is to use the \(S-\)cycle \(L_1\) to obtain the required contradiction. In the \textbf{fifth} step in the proof, we merge the cycles~\(\partial_0(L_1)\) and~\(\partial_0^+\) to obtain a bigger  cycle~\(\partial_{temp}\) and again list the properties of~\(\partial_{temp}\) needed for the next step. In the \textbf{sixth} crucial step, we use the properties of \(\partial_0(L_1)\) and \(\partial_{temp}\) obtained above to prove that there is a vacant square~\(Z_j\)~of~\(L_1\) \emph{all of whose edges lie in the interior of~\(\partial_{temp}.\)}

In the \textbf{seventh} final step, we argue that \(\partial_{temp}\) is precisely the ``intermediate"
cycle obtained in the iterative algorithm of the first step, after merging only vacant squares in the subsequence~\(L_1\) with the cycle~\(\partial_0^+.\) Merging the remaining squares in \(L \setminus L_1\) with \(\partial_{temp}\) we get \(D_{fin}.\) All edges of the vacant square \(Z_j \in L_1 \subset L\) are contained in the interior of \(\partial_{temp}\) and therefore also in the interior of \(D_{fin}.\) This is a contradiction since \(D_{fin}\) contains at least one edge from each square in \(L.\)

We provide the details below.
%contains some vacant square of \(L_1 \subset L\) completely in the interior

\subsubsection*{1. Merge all vacant squares in \(\Lambda^+\) with \(\partial_0^+.\)}
Setting \(C_0 := \partial_0^+,\) we merge the vacant squares in \(\Lambda^+\) with \(C_0\) in an iterative manner to obtain a finite sequence of cycles~\(\{C_i\}_{i \geq 1}.\)

%This is analogous to the merging procedure of the proof of~\((i)\) of Theorem~\ref{thm_lr}. %We then see that the final cycle obtained is the desired outermost boundary~\(\partial_H.\)

In the first step of our iterative procedure, we choose the square \(Y_1\) and merge it with \(\partial_0^+\) as in Theorem~3 of Ganesan~(2017) to get the new cycle~\(C_1.\) For \(i = 1,\) the cycle \(C_i\) has the following properties:\\
\((a1.1)\) Every occupied square of \(C^+(0)\) lies in the interior of \(C_i.\)\\
\((a1.2)\) If edge \(e \in C_i,\) then either \(e \in \partial_0^+\) or \(e\) is an edge of the square \(Y_l\) for some \(1 \leq l \leq t.\) Every edge in \(\partial_0^+\) either belongs to the cycle \(C_i\) or is contained in the interior of~\(C_i.\) The edge \(e_{j}\) of \(\partial^+_0\) belongs to \(C_i\) if and only if the corresponding vacant square \(Y_{j}\) lies in the exterior of \(C_i.\)\\
\((a1.3)\) The vacant squares \(Y_l ,1 \leq l \leq i\) and the edges \(\{e_l, 1 \leq l \leq i\} \subseteq \partial_0^+\) are contained in the interior of the cycle \(C_i.\)\\
\emph{Proof of \((a1.1)-(a1.3)\) for \(C_1\)}: From Theorem~\(3\) of Ganesan~(2017), the interior of the new cycle~\(C_1\) contains the interior of the cycle~\(C_0\) and the interior of the square~\(Y_1.\) Also the cycle~\(C_1\) consists only of edges of~\(C_0\) and~\(Y_1.\) Thus~\((a1.1)\) and the first statement of~\((a1.2)\) are true.

The second statement in \((a1.2)\) is true by construction. To prove the third statement, suppose that edge~\(e_j\) of~\(\partial_0^+\) belongs to~\(C_1.\) We recall that~\(e_j\) is common to the occupied square~\(W_j \in C^+(0)\) and the vacant square~\(Y_j \in \Lambda^+.\) Therefore either~\(W_{j}\) or~\(Y_{j}\) lies in the interior of~\(C_1\) but not both. Since~\(W_j\) lies in the interior of~\(C_1\) (property~\((a1.1)\)), we have that~\(Y_j\) lies in the exterior of~\(C_1.\)

Conversely, suppose that \(e_j \notin C_1.\) From the first statement of \((a1.2),\) we have that \(e_j\) is contained in the interior of \(C_1.\) Thus both the squares containing \(e_j\) as an edge belong to the interior of \(C_1.\) This proves~\((a1.2).\) To see~\((a1.3)\) is true, we have from \((a1.2)\) that if \(e_1 \in C_1,\) then \(Y_1\) lies in the exterior of~\(C_1.\) But as mentioned above, the cycle \(C_1\) contains \(Y_1\) in its interior. Thus~\((a1.3)\) is true. \(\qed\)

Using property \((a1.2),\) we proceed to the next step of the iteration. Fix the least index \(j\) such that the edge \(e_{j}\) of \(\partial^+_0\) belongs to the cycle \(C_1.\) If there is no such \(j,\) we stop the procedure. If there exists such an index \(j,\) then from property \((a1.2)\) above, the corresponding vacant square \(Y_{j}\) lies in the exterior of \(C_1.\) Merge \(Y_{j}\) and \(C_1\) using Theorem~\(3\) of Ganesan (2017) to get the new cycle \(C_{2}.\)

%(CHECK REF...HERE.. for merged version)

The new cycle \(C_2\) also satisfies properties \((a1.1)-(a1.3).\)\\
\emph{Proof of \((a1.1)-(a1.3)\) for \(C_2\)}: It suffices to verify that \(Y_1\) and \(Y_2\) both belong to the interior of \(C_2.\) The rest of the proof is as above. Since \(C_1\) satisfies \((a1.1)-(a1.3)\) we have that \(Y_1\) is contained in the interior of \(C_1\) and therefore contained in the interior of \(C_2.\) If \(Y_2\) is also contained in the interior of \(C_1,\) then we are done.

If not, then \(Y_2\) is in the exterior of \(C_1\) and in our iteration, we choose the least index \(j\) such that the edge \(e_j\) belongs to \(C_1.\) Since \(Y_2\) is in the exterior of \(C_1,\) we have from property \((a1.2)\) of cycle \(C_1\) that the edge \(e_2 \in C_1.\) Thus \(j = 2\) here and we merge~\(Y_2\) with~\(C_1\) to get~\(C_2.\) Thus~\(Y_2\) lies in the interior of~\(C_2.\) \(\qed\)

%To see \((a1.3)\) is true, we have from \((a1.2)\) that if \(e_1 \in C_3,\) then \(Y_1\) lies in the exterior of \(C_3.\) But as mentioned above, the cycle \(C_2\) contains \(Y_1\) in its interior. Thus \((a1.3)\) is true. \(\qed\)

As before, we use property \((a1.2)\) of cycle \(C_2\) to proceed to the next step of iteration. This process proceeds for a finite number of steps and the cycle~\(C_i\) obtained at iteration step \(i \geq 1\) satisfies properties \((a1.1)-(a1.3).\) Let \(D_{fin}\) denote the final cycle obtained after the procedure stops. The cycle \(D_{fin}\) satisfies the following properties.\\
\((b1.1)\) The cycle~\(D_{fin}\) contains only edges of vacant squares in~\(\Lambda^+.\)\\
\((b1.2)\) If \(e\) is an edge of the outermost boundary cycle~\(\partial_0^+,\) then~\(e \notin D_{fin}\) and~\(e\) lies in the interior of~\(D_{fin}.\) \\
\((b1.3)\) Every square in~\(C^+(0)\cup \Lambda^+\) lies in the interior of~\(D_{fin}.\)\\
\((b1.4)\) The cycle~\(D_{fin}\) is unique in the sense that if a cycle~\(C\) satisfies the above mentioned properties~\((b1.1)-(b1.3),\) then \(C = D_{fin}.\)\\

We recall that we say an edge \(e\) lies in the interior of a cycle~\(C\) if both the squares containing~\(e\) lie in the interior of~\(C.\)
%WRT CRFLLY+ Et+C..

\emph{Proof of \((b1.1)-(b1.4)\)}: Properties \((b1.1)-(b1.2)\) and \((b1.5)\) follow by construction of the cycle \(D_{fin}.\) To see \((b1.3),\) we have from~\((b1.2)\) that the outermost boundary cycle~\(\partial_0^+\) and therefore all squares contained in the interior of~\(\partial_0^+,\) are also in the interior of~\(D_{fin}.\) By construction every (vacant) square in~\(\Lambda^+\) is contained in the interior of \(D_{fin}.\)

To see \((b1.4),\) we suppose there is a cycle~\(C\) distinct from~\(D_{fin} \) that satisfies \((b1.1)-(b1.3)\) and suppose that~\(C\) contains an edge~\(e\) in the exterior of~\(D_{fin}.\) This means that some vacant square \(Y_j \in \Lambda^+\) lies in the exterior of~\(D_{fin},\) a contradiction to the fact that~\(D_{fin}\) satisfies property~\((b1.3).\)~\(\qed\)

%//done...ok...seems...GTHALS LATER SINCE WE NEED TO SEE THAT PARTIAL 0 L1 NEQ DFIN ONLY USING PROTEOS B1.1 AND B1.2...PRKVMMM+eTC....

%//WRT IT MORE PRECICELY...

\subsubsection*{2. Obtaining the star connected sequence~\(L\) from \(D_{fin}\)}
In this subsection, we use the cycle \(D_{fin}\) obtained above to construct the star connected \(S-\)cycle \(H_{out}.\) Letting \(D_{fin} = (f_1,...,f_r),\) there exists a unique vacant square \(Z_1 \in \Lambda^+\) that contains the edge~\(f_1.\) Indeed if two vacant squares \(Z_1\) and \(Z'_1\) in \(\Lambda^+\) both contain the edge~\(f_1,\) then necessarily one of them is in the exterior of \(D_{fin}.\) This contradicts property~\((a3)\) of~\(D_{fin}.\)

Similarly there exists a unique vacant square \(Z_2 \in \Lambda^+\) that has edge \(f_2.\) If \(Z_2 = Z_1,\) we proceed to \(f_3,\) else \(Z_2\) is star adjacent to \(Z_1\) and we add \(Z_2\) to the existing sequence and obtain \((Z_1,Z_2).\) Continuing this way, we obtain a final sequence of squares \(L = (Z_1,...,Z_s)\) such that \(Z_i\) is star adjacent to~\(Z_{i+1}\) for \(1 \leq i \leq s-1.\) %and \(Z_s\) is adjacent to \(Z_1.\)

By construction, the sequence \(L\) obtained is unique and has the following property.  \\
\((a2.1)\) The outermost boundary \(\partial_0(L)\) of \(L\) is \(D_{fin}.\)\\
\((a2.2)\) Every square in the sequence~\(L=(Z_1,\ldots,Z_s)\) contains an edge in~\(D_{fin}.\)\\

\emph{Proof of~\((a2.1)-(a2.2)\)}: Property~\((a2.2)\) is true by construction. We prove~\((a2.1)\) below. Let~\(G^*_V\) denote the graph with vertex set as corners of vacant squares in \(\Lambda^+\) and edge set being the edges of squares in \(\Lambda^+.\) From Theorem~1 of Ganesan (2017), the outermost boundary \(\partial_0(L)\) is a connected union of cycles in \(G^*_V\) and every square in the component~\(L\) is contained in the interior of some cycle in~\(\partial_0(L).\) %As in the proof of Theorem~\ref{thm1} above, we label squares as \(0\) and \(1\) obtain the corresponding outermost boundary \(\partial_0(L_1).\) %by labelling all squares in \(L_1\) as \(1,\) all squares sharing a vertex with \(L_1\) as in Theorem~\ref{thm1} above.

By property~\((b1.3),\) all squares of~\(\Lambda^+\) are contained in the interior of~\(D_{fin}.\) Since the sequence \(L\) consists of squares in~\(\Lambda^+,\) we have that every edge in the outermost boundary~\(\partial_0(L)\) either belongs to the cycle \(D_{fin}\) or is contained in its interior. If there exists an edge~\(e\) of \(D_{fin}\) not in \(\partial_0(L),\) then the edge \(e\) necessarily lies in the exterior of all cycles in \(\partial_0(L).\) Any square containing \(e\) as an edge also lies in the exterior of all cycles in \(\partial_0(L).\)

By construction, there is a vacant square~\(Z_j \in L\) that contains \(e\) as an edge and this square \(Z_j\) lies in the exterior of all cycles of~\(\partial_0(L),\) a contradiction.~\(\qed\)

\subsubsection*{3. The sequence~\(L\) contains a \(S-\)cycle}
We construct the vacant square graph \(G_V\) whose vertex set consists of the \emph{centres} of the vacant squares in \(\Lambda^+.\) Suppose \(v_1\) and \(v_2\) are two vertices in~\(G_V\) and let \(Y_{i_1} \in \Lambda^+\) and \(Y_{i_2} \in \Lambda^+\) be the vacant squares with centres~\(v_1\) and~\(v_2,\) respectively. We draw an edge between~\(v_1\) and~\(v_2\) if and only if~\(Y_{i_1}\) and~\(Y_{i_2}\) are star adjacent (i.e., share a corner).

For the sequence \(L,\) let \(H_{L}\) be the subgraph of~\(G_V\) obtained in the same manner above but using only squares in the sequence~\(L = (Z_1,\ldots,Z_s);\)  in other words, the vertex set of \(H_L\) consists of the centres of the (vacant) squares in~\(L.\) Suppose~\(u_1\) and~\(u_2\) are two vertices in~\(H_L\) that are centres of the squares \(Z_{j_1}\) and \(Z_{j_2}.\) We draw an edge between~\(u_1\) and~\(u_2\) if and only if~\(Z_{j_1}\) and~\(Z_{j_2}\) are star adjacent.

The following is the main Lemma we prove in this subsection.
\begin{Lemma}\label{h_l_cyc_lem}
\text{The graph \(H_L \subset G_V\) contains a cycle}.
\end{Lemma}
If~\(H_L\) is itself a cycle, then the component~\(L\) is an \(S-\)cycle and we are done. In Steps~\(4-7,\) we assume this is not the case and arrive at a contradiction .

We prove the above Lemma by using induction on induced tree subgraphs of the graph~\(G_V.\) We have the following properties.\\\\
\((a3.1)\) The graph \(H_L\) is a connected induced subgraph of \(G_V.\)\\\\
Suppose \(X\) is a tree and an induced subgraph of~\(G_V.\) Every vertex in \(X\) is the centre of a vacant square in \(\Lambda^+.\) For vertex \(v \in X,\) let  \(J_v\) be the vacant square in \(\Lambda^+\) containing \(v\) as its centre. Since \(X\) is connected, the union of the squares \(\{J_v\}_{v \in X}\) is also connected and forms a star connected component~\(C_X.\) We denote the outermost boundary of~\(C_X\) as \(\partial_0(X).\) \\\\
\((a3.2)\) The outermost boundary~\(\partial_0(X)\) does not contain any (occupied) square of the plus connected component~\(C^+(0)\) in its interior.\\

The properties~\((a3.1)\) and \((a3.2)\) imply Lemma~\ref{h_l_cyc_lem}.\\\\
\emph{Proof of~\((a3.1)\)}: By construction, the set of squares~\(\{Z_i\}_{1 \leq i \leq s}\) form a star connected component and so the corresponding graph~\(H_L\) is a connected subgraph of~\(G_V.\) To see that~\(H_L\) is also an induced subgraph of~\(G_V,\)  we let~\(u_1\) and~\(u_2\) be two vertices in~\(H_L.\) For \(i = 1,2\) let \(Z_{j_i}\) be the vacant square in~\(L\) with~\(u_i\) as centre. If there is an edge between the vertices~\(u_1\) and~\(u_2\) in the graph~\(G_V,\) then~\(Z_{i_1}\) and~\(Z_{i_2}\) are star adjacent. Since both~\(Z_{i_1}\) and~\(Z_{i_2}\) belong to~\(L,\) there is also an edge between~\(u_1\) and~\(u_2\) in the graph~\(H_L.\) \(\qed\)

%// IMP...PARRKVUM...

In what follows we prove \((a3.2)\) using induction.\\
\emph{Proof of~\((a3.2)\)}:  Let \(\#X\) denote the number of vertices in~\(X.\) If~\(\#X = 1\) so that~\(X = \{u\}\) for some vertex \(u \in G_V,\) then the outermost boundary~\(\partial_0(X)\) is simply the union of all four edges of the square~\(J_u.\) Thus~\((a3.2)\) is true for \(\#X = 1.\)

For \(\#X \geq 2,\) we proceed in two steps. In the first step, we apply induction and obtain the outermost boundary \(\partial_0(X_1)\) for a subtree~\(X_1\) containing one less vertex than~\(X.\) In the next step, we construct the outermost boundary~\(\partial_0(X)\) for the tree~\(X\) using~\(\partial_0(X_1)\) and use induction assumption to obtain~\((a3.2).\)

\subsubsection*{Defining the subtree~\(X_1\) and obtaining~\(\partial_0(X_1)\)}
Suppose~\((a3.2)\) is true for all induced tree subgraphs of~\(G_V\) containing at most \(n\) vertices and consider an induced tree~\(X \subset G_V\) with \(n+1\) vertices. Let \(u\) be any leaf of the tree~\(X\) that is adjacent to a vertex~\(v \in X.\) Writing \(X = X_1 \cup \{u\}\) where \(X_1 = X\setminus \{u\},\) we have that~\(X_1\) is an induced tree subgraph of~\(G_V\) containing~\(n\) vertices.

%and therefore the square \(J_u\) shares a vertex with the cycle~\(D_{i_0}.\) Without loss of generality, we assume that \(i_0 = r.\)\\

Let \(G(X_1)\) be the graph formed by the squares \(\{J_{w}\}_{w \in X_1}.\) From Theorem~\ref{thm_out}, we have that the outermost boundary~\(\partial_0(X_1) = \cup_{i=1}^{r} D_i\) is a connected union of cycles in~\(G(X_1)\) with mutually disjoint interiors. We have the following properties regarding the square~\(J_u\) and the cycles~\(\{D_i\}.\)\\
\((b3.1)\) The vacant square~\(J_u \in \Lambda^+\) containing~\(u\) as the centre is star adjacent only to the square~\(J_v\) containing~\(v\) as the centre and no other square having centre in the tree~\(X_1.\)\\
\((b3.2)\) Every square \(J_w\) with centre~\(w \in X_1\) is contained in the interior of some cycle~\(D_i, 1 \leq i \leq r.\) If \( i \neq j\) then~\(D_i\) and~\(D_j\) have mutually disjoint interiors and share at most one vertex in common. No cycle in~\(\{D_i\}_{1 \leq i \leq r}\) contains an (occupied) square of~\(C^+(0)\) in its interior. \\
\((b3.3)\) The square~\(J_v\) lies in the interior of a unique cycle~\(D_{i_0}, 1 \leq i_0 \leq r.\) Without loss of generality, we assume that \(i_0 = r.\)\\
\((b3.4)\) No cycle \(D_i, 1 \leq i \leq r-1,\) contains a vertex of the square~\(J_u.\)\\\\
\emph{Proof of \((b3.1)-(b3.4)\)}: The property~\((b3.1)\) is true since the original tree~\(X\) is an induced subgraph of~\(G_V.\) The first and second statements of property~\((b3.2)\) are true by Theorem~\ref{thm_out}. The final statement of property~\((b3.2)\) is true by induction assumption. Property~\((b3.3)\) is true since the cycles~\(\{D_i\}\) have mutually disjoint interiors.

%To prove the second half, we argue as follows. Every edge in the cycle~\(D_{i_0}\) belongs to some square~\(J_w, w \in X_1.\) In particular, every vertex in~\(D_{i_0}\) belongs to some square~\(J_w, w \in X_1.\) Since the square~\(J_u\) shares an endvertex with~\(J_v,\) we obtain~\((q3).\)

To prove~\((b3.4),\) we argue by contradiction. Suppose some vertex~\(z \in D_i, i \neq r\) shares an endvertex with~\(J_u.\) The vertex~\(z\) belongs to an edge~\(e_z \in D_i\) and by Theorem~\ref{thm_out}, the edge~\(e_z\) is the edge of a square~\(J_{v(z)}\) lying in the interior of \(D_i.\) Also,  the centre~\(v(z)\) of the square~\(J_{v(z)}\) belongs to~\(X_1.\) We have that
\begin{equation}\label{v_z_not_v}
J_{v(z)} \neq J_v
\end{equation}
and so the square~\(J_u\) shares endvertices with two squares~\(J_{v(z)}\) and~\(J_v,\) both containing centres in~\(X_1.\) This contradicts property~\((b3.1)\) and so~\((b3.4)\) is true. The property~(\ref{v_z_not_v}) is a consequence of the following two statements. From~\((b3.2),\) the cycles~\(D_i\) and~\(D_r\) have mutually disjoint interiors. The square~\(J_{v(z)}\) is contained in the interior of the cycle~\(D_i, i \neq r\) and the square~\(J_v\) is contained in the interior of the cycle~\(D_{r}.\) \(\qed\)

\subsubsection*{Constructing~\(\partial_0(X)\) from \(\partial_0(X_1)\)}
We now construct the outermost boundary~\(\partial_0(X)\) using the cycles~\(\{D_i\}_{1 \leq i \leq r}.\)
We consider three possible cases and see that~\((a3.2)\) is satisfied in each case; i.e., no cycle of the outermost boundary~\(\partial_0(X)\) contains any occupied square of the plus connected component~\(C^+(0)\) in its interior.\\\\
\emph{\underline{Case \(I\)}:} The first case we consider is when the square~\(J_u\) lies in the interior of some cycle~\(D_i, 1 \leq i \leq r.\) We have the following properties.\\
\((I.1)\) The square~\(J_u\) lies in the strict interior of the cycle~\(D_r\) in the sense that~\(J_u\) lies in the interior of~\(D_r\) and no edge of~\(D_r\) belongs to~\(J_u.\)\\
\((I.2)\) The outermost boundary \(\partial_0(X) = \partial_0(X_1) = \cup_{i=1}^{r}D_i.\)\\
Thus by induction assumption, we have that~\(\partial_0(X)\) contains no occupied square of the component~\(C^+(0)\) in its interior. This proves~\((a3.2)\) for Case~\(I.\)\\\\
\emph{Proof of~\((I.1)\)}: Suppose first that the square~\(J_u\) lies in the interior of the cycle~\(D_i\) for some \(1 \leq i \leq r-1.\) We recall that the square \(J_u\) also shares an endvertex with the square~\(J_v,\) which in turn lies in the interior of the cycle~\(D_r\) (property~\((b3.3)\)). The cycles \(D_i\) and \(D_r\) have at most one vertex in common and have mutually disjoint interiors (property~\((b3.2)\)). Therefore~\(D_i\) and~\(D_r\) have exactly one vertex~\(z\) in common and the vertex~\(z\) is also common to~\(J_u\) and~\(J_v.\) But this means that~\(D_i\) contains a vertex of the square~\(J_u\) contradicting property~\((b3.4)\) above.  Therefore, we have that the square~\(J_u\) lies in the interior of the cycle~\(D_r.\)

The second part regarding strict interior containment is true because, every edge in the cycle~\(D_r \in \partial_0(X_1)\) belongs to a square with centre in~\(X_1 = X \setminus \{u\}\) and contained in the interior of~\(D_r\) (property~\((iv),\) Theorem~\ref{thm_out}). In other words, any square \(J_w,w \in X,\) lying in the interior of~\(D_r\) and containing an edge in~\(D_r,\) necessarily has its centre~\(w \in X_1 = X \setminus \{u\}.\)\(\qed\)\\\\ %Since \(J_u\) lies in the interior of \(D_r,\) we obtain~\((I.1)\).\(\qed\)\\\\
\emph{Proof of~\((I.2)\)}: We recall that \(G(X_1)\) is the graph formed by the union of the squares \(\{J_w\}_{w \in X_1}\) with centres in~\(X_1=X\setminus\{u\}.\) The union of the cycles \(\cup_{1 \leq j \leq r} D_j\) as a subgraph of \(G(X_1)\) satisfy properties~\((i)-(v)\) of Theorem~\ref{thm_out}. Therefore, in the graph~\(G(X),\) the cycles \(\cup_{1 \leq j \leq r} D_r\) satisfy properties~\((ii)-(v).\) It only remains to see that property~\((i)\) is true; i.e., every edge in \(\cup_{1 \leq j \leq r} D_j\) is also an outermost boundary edge in the graph~\(G(X).\)

%and a square~\(Q_e\) with centre  not in~\(X_1\) (property~\((b),\) Lemma~\(3,\) Ganesan~(2015)). The square~\(Q_e\) lies in the exterior of all the cycles~\(D_j,1 \leq j \leq r.\) Since the square~\(J_u\) lies in the interior of~\(D_r,\) we also have that~\(Q_e \neq J_u.\)

Fix an edge \(e \in D_i, 1 \leq i \leq r.\) We arrive at a contradiction supposing that \(e\) is not an outermost boundary edge in the graph~\(G(X).\) Using the fact that the cycle~\(D_i\) is an outermost boundary cycle in~\(G(X_1),\) we have that the edge~\(e\) is the edge of a square~\(J_{v(e)}\) with centre~\(v(e) \in X_1 = X \setminus \{u\}.\) The above statement follows from the fact that the cycle~\(D_i\) satisfies property~\((b),\) Lemma~\ref{outer}.

Let \(E_e\) be the outermost boundary cycle in the graph~\(G(X)\) containing the square~\(J_{v(e)}\) in its interior. We recall that we have only assumed that the edge~\(e\) is not an outermost boundary edge in the graph~\(G(X).\) Intuitively, this also implies that the outermost boundary cycle~\(E_e \neq D_i\) and we state related properties.\\\\
\((I.3)\) The square~\(J_{v(e)}\) lies in the interior of the cycle~\(D_i.\) Every edge of the cycle~\(D_i\) either belongs to~\(E_e\) or lies in the interior of~\(E_e.\) Also \(E_e \neq D_i\) and so at least one edge of \(E_e\) lies in the exterior of~\(D_i.\)\\
\((I.4)\) At least one edge of \(E_e \setminus D_i\) belongs to the square~\(J_u\) and the square~\(J_u\) lies in the interior of the cycle~\(E_e.\)\\ %and exterior of the cycle~\(D_i.\)\\
We recall from property~\((I.1)\) that the square~\(J_u\) lies in the interior of the cycle~\(D_r.\)\\
\((I.5)\) The cycles \(D_r\) and \(E_e\) have more than one vertex in common and at least one edge of~\(D_r\) lies in the exterior of~\(E_e.\)

We use property~\((I.5)\) to prove the property~\((I.2)\) regarding the outermost boundary~\(\partial_0(X).\)\\\\
\emph{Proof of \((I.3)-(I.5)\)}: The first statement of~\((I.3)\) is true as follows. The cycle~\(D_i\) is the outermost boundary cycle containing the square~\(J_{v(e)}\) in its interior, in the graph~\(G(X_1).\) Therefore~\(D_i\) satisfies property~\((b),\) Lemma~\ref{outer} and so \(J_{v(e)}\) lies in the interior of~\(D_i.\) For the second statement, we argue as follows. The cycle~\(E_e\) is the outermost boundary cycle containing the square~\(J_{v(e)}\) in its interior, in the graph \(G(X).\) So~\(E_e\) satisfies property~\((c),\) Lemma~\ref{outer} and so every edge in \(D_i\) either belongs to~\(E_e\) or lies in the interior of~\(E_e.\)

The final statement of~\((I.3)\) is true as follows. Since the edge~\(e\) is not an outermost boundary edge in the graph~\(G(X),\) some cycle~\(C \subset G(X)\) contains the edge~\(e\) and therefore the square~\(J_{v(e)}\) containing~\(e\) as an edge, in its interior. The cycle~\(D_i\) also contains the square~\(J_{v(e)}\) in its interior and so merging~\(D_i\) and~\(C\) if necessary, we assume that every edge of~\(D_i\) either belongs to or lies in the interior of~\(C.\) Since the edge~\(e\) belongs to~\(D_i\) and lies in the interior of~\(C,\) we also have that~\(C \neq D_i.\) Finally, the outermost boundary cycle~\(E_e\) satisfies property~\((c),\) Lemma~\ref{outer} and so every edge of~\(C\) either belongs to or lies in the interior of~\(E_e.\) In other words~\(E_e \neq D_i.\)

%// TO SEE CRFLLY +eTC...

To prove \((I.4),\) we argue as follows. Using properties~\((b3.4)\) and~\((I.1),\) we have that the cycle \(D_i\) consists only edges of squares with centre in~\(X_1 = X \setminus \{u\}.\) Suppose that all the edges in \(E_e \setminus D_i\) also belong to squares with centres in \(X_1 = X \setminus \{u\}.\) This means that all edges of the cycle~\(E_e\) belong to squares with centres in \(X_1 = X\setminus \{u\}.\) But~\(E_e \neq D_i\) and~\(E_e\) contains at least one edge in the exterior of the cycle~\(D_i\) (property~\((I.3)\)). This contradicts the fact that~\(D_i\) is the outermost boundary cycle containing the square~\(J_{v(e)}\) in the graph~\(G(X_1)\) and satisfies property~\((c),\) Lemma~\ref{outer}.

From the above paragraph, we have that there exists at least one edge \(f \in E_e \setminus D_i\) which belongs to the square~\(J_u.\) If \(J_u\) lies in the exterior of \(E_e,\) we would merge \(E_e\) and \(J_u\) using Theorem~\(3\) of Ganesan~(2017) to obtain a bigger cycle \(E_{eu}\) in \(G(X)\) with the following property. The cycle~\(E_{eu}\) would contain \(E_e\) in its interior and at least one edge of \(E_{eu}\) would lie in the exterior of \(E_e.\) This would contradict the fact that \(E_e\) is an outermost boundary cycle in the graph~\(G(X)\) and satisfies property~\((c),\) Lemma~\ref{outer}. Thus the square~\(J_u\) lies in the interior of~\(E_e\) and this proves~\((I.4).\)

We prove~\((I.5)\) as follows. From~\((I.1),\) we have that the cycle~\(D_r\) contains no edge of~\(J_u\) and the square \(J_u\) lies in the interior of~\(D_r.\) From~\((I.4),\) we obtain that the cycle~\(E_e\) contains at least one edge of the square~\(J_u\) and the square~\(J_u\) lies in the interior of~\(E_e.\) So there exists at least one edge of~\(D_r\) lying in the exterior of~\(E_e.\) The cycles~\(E_e\) and \(D_r\) cannot have mutually disjoint interiors since they both contain the square~\(J_u\) in their respective interiors and so~\(E_e\) and~\(D_r\) have more than one vertex in common.\(\qed\)

Using property~\((I.5)\), we merge paths of~\(D_r\) lying in the exterior of~\(E_e\) in an iterative manner as in Theorem~\(3\) of Ganesan~(2017) to obtain a final cycle~\(E''_e\) containing both the cycles~\(E_e\) and~\(D_r\) in its interior. The cycle \(E''_e \neq E_e\) and this contradicts the fact that~\(E_e\) is the outermost boundary cycle containing the square~\(J_{v(e)}.\) This proves that~\(E_e = D_i\) and so every edge in~\(D_i\) is an outermost boundary cycle in the graph~\(G(X).\) This completes the proof of~\((I.2).\)~\(\qed\)

%//TO WRITE BELOW ALSO AS PROPRTIES ETC... PRKVMM+ETC..\\\\
\emph{\underline{Case \(II\)}}: In this case, the square~\(J_u\) shares only a corner with the square~\(J_v\) and lies in the exterior of all cycles~\(\{D_i\}_{1 \leq i \leq r}.\) We have the following property.\\
\((II.1)\) The square~\(J_u\) shares a unique vertex with the cycle~\(D_r\) and does not share a vertex with any other cycle~\(D_i, 1 \leq i \leq r-1.\) The square~\(J_u\) does not share an edge with any cycle~\(D_i, 1 \leq i \leq r.\)\\
\((II.2)\) The outermost boundary~\(\partial_0(X) = \cup_{i=1}^{r}D_i \cup \{J_u\}.\)\\\\
In other words, the outermost boundary for the squares with centres in~\(X\) is simply the union of the outermost boundary~\(\partial_0(X_1)\) and the cycle formed by the four edges of~\(J_u.\) Since~\(J_u\) is vacant, we have from property~\((b3.2)\) that~\(\partial_0(X)\) contains no occupied square of the component~\(C^+(0)\) in its interior. This proves~\((a3.2)\) for Case~\((II).\)\\\\
\emph{Proof of~\((II.1)\)}: From property~\((b3.3),\) we have that the square \(J_v\) sharing a vertex with the square~\(J_u\) lies in the interior of the cycle~\(D_r.\) In this Case~\(II,\) we have assumed that \(J_v\) and \(J_u\) share a single vertex~\(z\) and the square~\(J_u\) lies in the exterior of all cycles~\(\{D_j\}_{1 \leq j \leq r}.\) In particular, the square~\(J_u\) lies in the exterior of~\(D_r\) and so shares the vertex~\(z\) with~\(D_r.\)

To see that the square \(J_u\) does not share any other vertex with the cycle~\(D_r,\) we argue as follows. Suppose that some vertex \(y \in D_r\) also belongs to~\(J_u\) and~\(y \neq z.\) The vertex~\(y\) belongs to an edge \(e_y \in D_r\) which is also the edge of a square \(J_{v(y)}\) with centre~\(v(y) \in X_1 = X\setminus\{u\}.\) The final statement is true since \(D_r\) is the outermost boundary cycle in the graph~\(G(X_1)\) and so satisfies property~\((b),\) Lemma~\ref{outer}.

If \(J_{v(y)} = J_v,\) we then obtain that \(J_u\) shares an edge with \(J_v,\) a contradiction. If \(J_{v(y)} \neq J_v,\) then \(J_u\) shares vertices with two distinct squares with centres in~\(X_1 = X\setminus \{u\},\) a contradiction to the fact that \(u\) is a leaf in~\(X.\) From the above argument, we also obtain that \(J_u\) does not share an edge with the cycle~\(D_r.\) Finally, from property~\((b3.4),\) we obtain that~\(J_u\) does not share a vertex with any other cycle~\(D_i, 1 \leq i \leq r-1.\) This proves~\((II.1).\)\(\qed\)\\\\
\emph{Proof of~\((II.2)\)}: It suffices to see that properties~\((i)-(v)\) in the statement of Theorem~\(1\) of Ganesan~(2017) holds. We first prove properties \((ii)-(v).\) Property~\((ii)\) holds since~\(\partial_0(X_1)\) is connected and the square~\(J_u\) shares a vertex with the cycle~\(D_r \in \partial_0(X_1)\) (see property~\((II.1)\) above).

Property~\((iii)\) holds since the square~\(J_u\) shares a single vertex with~\(D_r\) and does not share a vertex with any other cycle~\(D_i, 1 \leq i \leq r-1\) (property~\((II.1)\)). Property~\((iv)\) is true since every square with centre in~\(X_1 = X \setminus \{u\}\) is contained in some cycle \(D_i, 1 \leq i \leq r\) (property~\((b3.2)\)) and the square~\(J_u\) with centre \(u =  X\setminus X_1\) is contained in the cycle formed by the four edges of~\(J_u.\)

To see property~\((v)\) is true, we must see that every edge in~\(\cup_{1 \leq i \leq r} D_i \cup \{J_u\}\) is a boundary edge in the graph~\(G(X).\) Let~\(e \in D_i, 1 \leq i \leq r,\) be any edge. Applying~Theorem~\ref{thm_out}, property~\((iv),\) to the outermost boundary~\(\partial_0(X_1) = \cup_{1 \leq j \leq r} D_j,\) we have that~\(e\) is a boundary edge in the graph~\(G(X_1)\) and so~\(e\) is adjacent to a square~\(J_{v(e)}\) with centre~\(v(e) \in X_1\) and a square~\(Q_e\) with centre not in~\(X_1.\) We recall that~\(G(X_1)\) is the graph formed by the squares~\(\{J_{w}\}_{w \in X_1}\) with centres in~\(X_1.\) Moreover, the square~\(J_{v(e)}\) lies in the interior of~\(D_i\) and the square~\(Q_e\) lies in the exterior of all cycles~\(\{D_j\}_{1 \leq \leq r}.\)

%Since the square~\(Q_e\) does not have its centre in \(X_1 = X \setminus\{u\},\) the final statement of the previous paragraph would be true if \(Q_e \neq J_u.\) To see that~\(Q_e \neq J_u,\) we argue as follows.

We recall that to prove that property~\((v)\) holds for~\(\cup_{1 \leq j \leq r} D_j \cup \{J_u\},\) we need to show that the square~\(Q_e\) lies in the exterior of all the cycles in~\(\cup_{1 \leq i \leq r}\{D_i\} \cup \{J_u\}\) and also that~\(Q_e\) does not have centre in~\(X.\) Since~\(Q_e\) lies in the exterior of the cycles~\(\{D_j\}_{1 \leq j \leq r}\) (see final statement, previous paragraph), it is therefore enough to see that~\(Q_e \neq J_u\) and we argue as follows. The square~\(Q_e\) shares the edge~\(e\) with the square~\(J_{v(e)}, v(e) \in X_1.\) Since the square~\(J_u\) shares only a vertex with the square~\(J_v\) and does not share an edge with any square in~\(\{J_w\}_{w \in X_1},\) we have~\(Q_e \neq J_u.\) This proves that property~\((v)\) of Theorem~\ref{thm_out} is true for every edge~\(e \in D_i, 1 \leq i \leq r;\) i.e., every edge in~\(\cup_{1 \leq j \leq r}D_j\) is a boundary edge in the graph~\(G(X).\)

%It remains to see that~\(Q_e\) lies in the exterior of the cycles~\(\{D_j\}_{1 \leq j \leq r}.\) The cycles~\(\{D_j\}_{1 \leq j \leq r}\) have mutually disjoint interiors and share at most one vertex in common (property~\((b3.2)\)). Therefore, if the square~\(Q_e\) lies in the interior of some cycle~\(D_j,j \neq i\) then necessarily~\(D_i\) and~\(D_j\) must share the edge~\(e \in Q_e.\) This contradicts property~\((b3.2)\) and so \(Q_e\) lies in the exterior of the cycles~\(\cup_{1 \leq i \leq r} D_i.\)  This proves that property~\((v)\) is true for every edge~\(e \in D_i, 1 \leq i \leq r.\)

%Therefore the square~\(Q_e\) does not have its centre in~\(X\) and lies in the exterior of all cycles of~\(\cup_{1 \leq i \leq r}D_i \cup \{J_u\}.\) So the edge~\(e\) is also a boundary edge in the graph~\(G(X).\) This proves that every edge in the union of the cycles~\(\cup_{1 \leq i \leq r} D_i\) is a boundary edge.

%formed by the edges of the squares with centres in~\(X.\)

%// WRT ASSUMPTIOSN ETC...CRFLLY...

Consider now an edge~\(f \in J_u.\) The edge~\(f\) is also the edge of a square~\(Q_f\) and since~\(J_u\) does not share an edge with a square having centre in~\(X_1 = X \setminus \{u\},\) we have that \(Q_f\) does not have its centre in~\(X.\) Suppose now that~\(Q_f\) lies in the interior of some cycle~\(D_j, 1 \leq j \leq r.\) Since the square~\(J_u\) lies in the exterior of all cycles~\(\{D_j\}_{1 \leq j \leq r}\) by assumption, the cycle~\(D_j\) and the square~\(J_u\) must share an edge. But by property~\((II.1),\) the square~\(J_u\) does not share an edge with any cycle~\(D_i, 1 \leq i \leq r\) and so we obtain a contradiction. Thus the square~\(Q_f\) lies in the exterior of all cycles~\(\{D_j\}_{1 \leq j \leq r}\cup\{J_u\}.\) So every edge in~\(J_u\) is also a boundary edge in the graph~\(G(X).\) %This proves that property~\((v)\) holds for~\(\cup_{1 \leq i \leq r} D_i \cup \{J_u\};\) i.e., every edge in~\(\cup_{1 \leq j \leq r}D_j \cup \{J_u\}\) is a boundary edge in the graph~\(G(X).\)

It remains to see that property~\((i)\) of Theorem~\ref{thm_out} holds; i.e., every edge \(e \in \cup_{1 \leq i \leq r}D_i \cup \{J_u\}\) is an outermost boundary edge in the graph~\(G(X).\) Fix an edge \(e \in D_i, 1 \leq i \leq r.\) As in Case~\(I,\) we arrive at a contradiction supposing that~\(e\) is not an outermost boundary edge (in the graph~\(G(X)\)). We use the fact that the cycle~\(D_i\) is an outermost boundary cycle in~\(G(X_1)\) and so the edge~\(e\) is the edge of a square~\(J_{v(e)}\) with centre~\(v(e) \in X_1 = X \setminus \{u\}.\) The above statement follows from the fact that the cycle~\(D_i\) satisfies property~\((b),\) Lemma~\ref{outer}. Let \(E_e\) be the outermost boundary cycle in the graph~\(G(X)\) containing the square~\(J_{v(e)}\) in its interior. Intuitively, as in Case~\(I,\) we must have~\(E_e \neq D_i\) and we state related properties.\\\\
\((II.3)\) The square~\(J_{v(e)}\) lies in the interior of the cycle~\(D_i.\) Every edge of the cycle~\(D_i\) either belongs to~\(E_e\) or lies in the interior of~\(E_e.\) Also \(E_e \neq D_i\) and so at least one edge of~\(E_e\) lies in the exterior of~\(D_i.\)\\
\((II.4)\) At least one edge of \(E_e \setminus D_i\) belongs to the square~\(J_u\) and the square~\(J_u\) lies in the interior of the cycle~\(E_e.\)\\
\((II.5)\) The squares~\(J_{v(e)}\) and~\(J_u\) are distinct and the cycle~\(E_e\) contains both~\(J_{v(e)}\) and~\(J_u\) in its interior.

%The square~\(J_{v(e)}\) lies in the interior of the cycle~\(D_i.\) %and exterior of the cycle~\(D_i.\)\\

We use property~\((II.5)\) to prove the property~\((II.2)\) regarding the outermost boundary~\(\partial_0(X).\)\\\\
\emph{Proof of \((II.3)-(II.5)\)}: The proof of~\((II.3)-(II.4)\) is the same as the corresponding proofs of~\((I.3)-(I.4).\) The square~\(J_u \neq J_{v(e)}\) since the square~\(J_{v(e)}\) lies in the interior of the cycle~\(D_i\) and by assumption, the square~\(J_u\) lies in the exterior of all cycles~\(\{D_j\}_{1 \leq j \leq r}.\) This proves that~\(J_{v(e)} \neq J_u.\)

From property~\((II.4),\) we have that the square~\(J_u\) lies in the interior of the cycle~\(E_e.\) By definition, the square~\(J_{v(e)}\) also lies in the interior of~\(E_e.\) This proves the second part of~\((II.5).\)~\(\qed\)

%Since \(E_e\) is the outermost boundary cycle containing the square~\(J_{v(e)}\) in the graph \(G(X)\) and \(E_e \neq D_i,\) the first statement of property~\((I.3)\) follows from the fact that \(E_e\) satisfies property~\((c),\) Lemma~\(3,\) Ganesan~(2105). To see the second statement we argue as follows. Since the edge~\(e\) is not an outermost boundary edge in the graph~\(G(X),\) some cycle~\(C\) in~\(G(X)\) contains~\(e\) in its interior. In particular, the square~\(J_{v(e)}\) also lies in the interior of \(C.\)

%The cycles~\(D_i\) and~\(C\) both contain the square~\(J_{v(e)}\) in their interiors. Merging the cycles~\(C\) and \(D_i\) using Theorem~\(3\) of Ganesan~(2017) if necessary, we therefore assume that the cycle~\(D_i\) also lies in the interior of~\(C.\) Since the edge \(e\) lies in the interior of~\(C\) and \(e \in D_i,\) we also have that~\(C \neq D_i.\) Finally, the outermost boundary cycle~\(E_e\) satisfies property~\((c),\) Lemma~\(3,\) Ganesan~(2017) and so the cycle \(C\) in turn is contained in the interior of the outermost boundary cycle~\(E_e.\) Thus~\(E_e \neq D_i\) and this proves~\((I.3).\)

Using~\((II.3)-(II.5)\) we arrive at a contradiction as follows. Suppose now that the cycle~\(E_e\) contains \emph{only} the two squares~\(J_u\) and~\(J_{v(e)}\) in its interior. This means that \(J_u\) shares an edge with the square~\(J_{v(e)},\) a contradiction since we have assumed in this Case~\(II\) that the square~\(J_u\) shares only a vertex with~\(J_v\) and does not share an edge with any square in~\(\{J_w\}_{w \in X_1}.\)

\begin{figure}[tbp]
\centering
%\fbox{
\includegraphics[width=2in, trim= 100 450 300 100, clip=true]{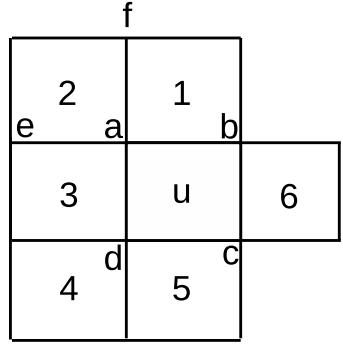}
%}
\caption{The squares~\(J_{q_1}\) and~\(J_u\) share the edge~\(bc.\)}
\label{fig_q1q2}
\end{figure}

Suppose now that~\(E_e = (g_1,\ldots,g_t)\) contains more than two squares in its interior. Let \(g_1 = f\) be an edge belonging to the square~\(J_u.\) We refer to Figure~\ref{fig_q1q2}, where the square \(J_u\) has corners~\(a,b,c\) and \(d\) and has label~\(u.\) Let \(g_1 = ad\) be the left edge of the square~\(J_u.\) By property~\((b),\) Lemma~\ref{outer}, every edge in~\(E_e\) is the edge of some square with centre in~\(X\) and a square with centre not in~\(X.\) Therefore, we have the following property.
\begin{equation}\label{three_sq}
\text{The square labelled~\(3\) in Figure~\ref{fig_q1q2} does not have its centre in~\(X.\)}
\end{equation}
Let \(i_1 \geq 2\) be the smallest index such that the edge~\(g_{i_1} \in E_e\) belongs to a square~\(J_{q_1}\) whose centre~\(q_1 \in X_1 = X\setminus \{u\}.\) Similarly let~\(i_2 \leq t\) be the largest index such that the edge~\(g_{i_2} \in E_e\) belongs to a square~\(J_{q_2}, q_2 \in X_1 = X \setminus \{u\}.\) The square~\(J_u\) shares a vertex with both~\(J_{q_1}\) and~\(J_{q_2}.\)

If \(q_1 \neq q_2,\) we obtain a contradiction, since the vertex~\(u\) is a leaf in~\(X\) and so the square~\(J_u\) shares only a vertex with the square~\(J_v, v \in X_1\) and does not share a vertex with any other square with centre in~\(X_1.\) If~\(q_1 = q_2,\) then we have the following property.
\begin{equation}\label{shr_edge}
\text{The squares~\(J_{q_1}\) and~\(J_u\) share the edge~\(bc\)}.
\end{equation}
The statement~(\ref{shr_edge}) is again a contradiction since we have assumed in this Case~\(II\) that the square~\(J_u\) shares only a vertex with~\(J_v\) and does not share an edge with any square with centre in~\(X.\) So~every edge \(e \in D_i, 1 \leq i \leq r\) is an outermost boundary edge in the graph~\(G(X).\)\\\\
\emph{Proof of~(\ref{shr_edge})}: Referring to Figure~\ref{fig_q1q2}, we first see that the subpath formed by the segments \(ab,ad\) and \(dc\) belongs to \(E_e;\) i.e., we must have \(g_2 = ab\) and \(g_t = cd.\)

Suppose \(g_2 \neq ab\) and \(g_t \neq cd.\) This means that either~\(g_2 = ae\) or~\(g_2 = ad.\) In either case, the edge \(g_2\) belongs to a square~\(J_{q_1}\) that shares the vertex~\(a\) with~\(J_u.\) Arguing similarly, the square~\(J_{q_2}\) shares the vertex \(d\) with \(J_u.\)
If \(q_1 = q_2,\) then \(J_{q_1}\) necessarily contains the edge \(ad\) as an edge and so \(J_{q_1}\) is the square labelled~\(3\) in Figure~\ref{fig_q1q2}.
This contradicts property~(\ref{three_sq}).

From the discussion in the above paragraph, we have that either \(g_2 = ab\) or \(g_t = cd.\) Suppose now that \(g_2 = ab\) but \(g_t \neq cd.\) Arguing as before, the square \(J_{q_1}\) then contains the vertex \(b\) and the square~\(J_{q_2}\) contains the vertex~\(d.\) But since \(q_1 = q_2\) and \(q_1 \neq u,\) this is a contradiction. An analogous contradiction is obtained if we assume that \(g_2 \neq ab\) and \(g_t = cd.\) From the above two paragraphs, we therefore obtain that~\(g_2 = ab\) and~\(g_t = cd.\) Thus the square~\(J_{q_1} = J_{q_2}\) contains both~\(b\) and~\(c\) as vertices and so the square~\(J_{q_1}\) is the square labelled~\(6\) in Figure~\ref{fig_q1q2}.\(\qed\)

%//TO WRT MORE ABOVE +eTC.. PRKVMMM+eTC...

%// TO ALSO WRT THAT NO OTHER POSSBLE BOUNDARY EDGE ETC...

It remains to see that every edge in the square~\(J_u\) is also an outermost boundary edge in the graph~\(G(X).\) Suppose~\(e \in J_u\) is not an outermost boundary edge in~\(G(X)\) and let~\(E_u\) be the outermost boundary cycle in~\(G(X),\) containing the square~\(J_u\) in its interior. We have the following properties.\\
\((II.6)\) The cycle~\(E_{u} \neq J_u\) and at least one edge of the cycle~\(E_u\) belongs to the square~\(J_u.\)\\
\((II.7)\) The cycle~\(E_u\) contains~\(J_u\) and at least one other square~\(J_w, w \in X_1\) in its interior.\\
\emph{Proof of~\((II.6)-(II.7)\)}: By definition, there is a cycle~\(C\) in the graph~\(G(X)\) containing the edge~\(e\) in its interior. The cycle~\(C\) therefore contains the square~\(J_u\) in its interior. Also the other square~\(Q_u\) containing~\(e\) as an edge also lies in the interior of~\(C.\) Since~\(E_u\) satisfies property~\((c),\) Lemma~\ref{outer}, we have that~\(C\) itself is contained in the interior of~\(E_u\) and so~\(E_u \neq J_u.\)

We prove the remaining part of~\((II.6)\) as follows. Suppose that the cycle~\(E_u\) contains no edge of the square~\(J_u.\) We then have that~\(E_u\) contains only edges of squares with centres in~\(X_1.\) Fix an edge~\(f \in E_u.\) Since~\(E_u\) satisfies property~\((b),\) Lemma~\ref{outer}, the edge~\(f\) belongs to a square~\(J_{v(f)}\) with centre in~\(X.\) Moreover, the square~\(J_{v(f)}\) lies in the interior of~\(E_u.\) Since \(f \notin J_u,\) we have that~\(J_{v(f)} \neq J_u.\) Therefore~\(J_{v(f)}\) has its centre~\(v(f) \in X_1 = X \setminus \{u\}.\) From property~\((iv),\) Theorem~\ref{thm_out}, we have that the square~\(J_{v(f)}\) is contained in the interior of one of the cycles~\(D_{j_0}, 1 \leq j_0 \leq r\) of the outermost boundary~\(\partial_0(X_1).\)
From property~\((d)\) following Lemma~\ref{outer}, we have that the outermost boundary cycle in the graph~\(G(X_1)\) containing~\(J_{v(f)}\) in its interior, is also~\(D_{j_0}.\)

The cycle~\(D_{j_0}\) satisfies property~\((c),\) Lemma~\ref{outer} and so the cycle~\(E_u\) containing~\(J_u\) in its interior, must itself lie in the interior of~\(D_{j_0}.\) This means that the square~\(J_u\) is contained in the interior of~\(D_{j_0},\) a contradiction to our assumption that~\(J_u\) lies in the exterior of all cycles~\(\{D_j\}_{1 \leq j \leq r}.\) Thus~\(E_u\) contains at least one edge of the square~\(J_u\) and this proves~\((II.6).\)

We prove property~\((II.7)\) as follows. Using property~\((II.6),\) let \(g \in E_u \setminus J_u\) be any edge. Since~\(E_u\) satisfies property~\((b),\) Lemma~\ref{outer}, the edge~\(g\) is the edge of a square~\(J_{v(g)}\) with centre~\(v(g) \in X_1 = X \setminus~\{u\}.\) The square~\(J_{v(g)}\) lies in the interior of~\(E_u.\) Since~\(g \notin J_u,\) we have that~\(J_{v(g)} \neq J_u.\)~\(\qed\)

Using~\((II.6)-(II.7)\) and an analogous analysis as following~\((II.3)-(II.5),\) we obtain a contradiction and so every edge in~\(J_u\) is an outermost boundary edge in the graph~\(G(X).\) Finally, every other edge of~\(G(X)\) not present in~\(\cup_{1 \leq j \leq r} D_j \cup \{J_u\}\) necessarily lies in the interior of one of the cycles in~\(\{D_j\}.\) This proves property~\((i)\) and therefore~\((II.2).\)

%Suppose now that~\(E_u\) contains only the two squares~\(J_{w}\) and~\(J_u\) in its interior. This means that the square~\(J_u\) shares an edge with

%//WRT MORE ABV +etC..

%now proceed in an analogous manner as before, to obtain the desired contradiction. //WRT MORE +ETC...\(\qed\)\\\\
\emph{\underline{Case \(III\)}}: In this case, the square~\(J_u\) shares an edge~\(e_{uv}\) with the square~\(J_v\)  and lies in the exterior of all cycles~\(D_i, 1 \leq i \leq r.\) We have the following property.\\
\((III.1)\) The edge~\(e_{uv}\) belongs to the cycle~\(D_r \in \partial_0(X_1)\) and does not belong to any other cycle~\(D_i, 1 \leq i \leq r-1.\) Apart from the endvertices of~\(e_{uv},\) the square~\(J_u\) does not share any other vertex with~\(D_r.\)

%The square~\(J_u\) lies in the exterior of all cycles~\(D_i, 1 \leq i \leq r.\)

Using property~\((III.1),\) we construct~\(\partial_0(X)\) from \(\partial_0(X_1)\) as follows. We set \(E_i = D_i\) for~\(1 \leq i \leq r-1\) and let \begin{equation}\label{er_def}
E_{r} = (D_r\setminus e_{uv}) \cup (J_u \setminus e_{uv}).
\end{equation}
Here \(J_u \setminus e_{uv}\) is the path formed by the three edges of the square~\(J_u\) apart from the edge~\(e_{uv}.\) We have the following additional properties.\\
\((III.2)\) The graph \(E_r\) is a cycle and contains the square~\(J_u\) in its interior. The interior of the cycle~\(E_r\) is the union of the interior of the cycles~\(D_r\) and~\(J_u.\)\\
\((III.3)\) The outermost boundary~\(\partial_0(X) = \cup_{1 \leq i \leq r} E_i.\)\\
In this case, a square~\(Q\) belongs to the interior of a cycle of~\(\partial_0(X)\) if and only if either~\(Q = J_u\) or~\(Q\) lies in the interior of a cycle of~\(\partial_0(X_1).\) The square~\(J_u\) is vacant and so by induction assumption, the outermost boundary~\(\partial_0(X)\) contains no occupied square of the component~\(C^+(0)\) in its interior. This proves~\((a3.2)\) for Case~\((III).\)\\\\
\emph{Proof of \((III.1)-(III.2)\)}: To prove~\((III.1),\) we argue as follows. From property~\((b3.3)\) above, we have that the square~\(J_v\) lies in the interior of~\(D_r\) and by assumption, the square~\(J_u\) lies in the exterior of all the cycles~\(\{D_i\}_{1 \leq i \leq r}\) and shares an edge with the square~\(J_v.\) Therefore, the cycle~\(D_r\) contains the edge~\(e_{uv} \in D_r.\) Moreover, the cycle~\(D_r\) shares at most one vertex with~\(D_j, 1 \leq j \leq r-1\) (property~\((b3.2)\)) and so we have that~\(e_{uv} \notin D_j, 1 \leq j \leq r-1.\) Also, arguing as in the proof of property~\((b3.4),\) we obtain that the square~\(J_u\) does not share any other vertex with~\(D_r.\) %since otherwise, we would obtain that~\(J_u\) is also adjacent to some

To see the first part of~\((III.2)\) that the graph~\(E_r\) is a cycle, we argue as follows. The square~\(J_u\) does not share an endvertex with~\(J_w, w \neq v, w \in X_1.\) Therefore apart from the endvertices, the path~\(D_r \setminus e_{uv}\) does not share any other vertex with the square~\(J_u\) and so~\(E_r\) is a cycle. The final statement is true since every vertex of the path~\(J_u \setminus e_{uv}\) lies in the exterior of the cycle~\(D_r.\)~\(\qed\)  %   //WRT MORE ?? PRKMVMM +etC...\(\qed\)

\emph{Proof of~(III.3)}: As before, we first prove properties \((ii)-(v)\) of Theorem~\ref{thm_out}. Property~\((ii)\) is true as follows. The union~\(\cup_{i=1}^{r}D_i = \partial_0(X_1)\) is connected and we remove an edge \(e_{uv}\) belonging to the cycle~\(D_r.\) Therefore the graph~\(\cup_{i=1}^{r-1}D_i \cup (D_r \setminus e_{uv})\) is still connected. We then add the path \(J_u\setminus e_{uv}\) with same endvertices as the edge~\(e_{uv}.\) Therefore, \(\cup_{i=1}^{r-1}D_i \cup (D_r \setminus e_{uv}) \cup (J_u\setminus e_{uv})\) is also connected.

We prove property~\((iii)\) as follows. Fix \(1 \leq i \neq j \leq r-1.\) The cycles~\(E_i = D_i\) and~\(E_j = D_j\) have disjoint interiors and have at most one vertex in common (see property~\((b3.2)\)). Suppose~\(1 \leq i \leq r-1\) and~\(j = r.\) The cycle~\(D_r\) shares at most one endvertex with~\(D_i\) and so~\(D_r \setminus e_{uv}\) shares at most one endvertex with~\(D_i.\) From property~\((b3.4),\) the square~\(J_u\) shares no endvertex in common with~\(D_i.\) Therefore \(E_r = (D_r\setminus e_{uv}) \cup (J_u \setminus e_{uv})\) shares at most one endvertex with~\(D_i.\)

To see property~\((iv)\) is true we argue as follows. Let \(J_w\) be any square with centre \(w \in X.\) If \(w = u,\) then \(J_u\) is contained in the interior of the cycle~\(E_r\) by property~\((III.2).\) If \(w \neq u,\) then~\(w \in X_1 = X \setminus \{u\}\) and so the square~\(J_w\) is contained in the interior of one of the cycles~\(D_i, 1 \leq i \leq r.\) This is seen by applying Theorem~\ref{thm_out}, property~\((iv)\) to the outermost boundary~\(\partial_0(X_1) = \cup_{1 \leq j \leq r} D_j.\) If \(1 \leq i \leq r-1,\) then~\(J_w\) is also contained in the interior of the cycle~\(E_i = D_i.\) If \(i = r,\) then~\(J_w\) is contained in the interior of the cycle~\(D_r,\) which in turn is contained in the interior of the cycle~\(E_r\)~(property~\((III.2)\)).

%because of the following two statements. For \(1 \leq i \leq r-1,\) we have that \(E_i = D_i\) have identical interiors. The cycle~\(E_r\) contains the cycle \(D_r\) and the square \(J_u\) in its interior. %and the square~\(J_u\) is contained in the interior of the cycle~\(D_r.\)

To see property~\((v)\) is true, we need to see that every edge in~\(\cup_{1 \leq j \leq r} E_j\) is a boundary edge adjacent to one square with centre in~\(X\) and one square with centre not in~\(X.\) Suppose first that~\(e \in E_i \cap D_i\) for some~\(1 \leq i \leq r;\) i.e., either~\(e \in D_i = E_i\) for some~\(1 \leq i \leq r-1\) or~\(e \in E_r \cap D_r = D_r \setminus e_{uv}.\) We use the fact that~\(D_j, 1 \leq j \leq r\) is an outermost boundary cycle in the graph~\(G(X_1).\) Here~\(G(X_1)\) is the graph formed by the edges of squares with centres in~\(X_1.\)

Applying Theorem~\ref{thm_out}, property~\((v)\) to the edge~\(e\) of the outermost boundary~\(\partial_0(X_1),\) we have that~\(e\) belongs to a square~\(J_{v(e)}\) with centre~\(v(e) \in X_1\) and a square~\(Q_e\) with centre not in~\(X_1.\) We have the following properties.
\begin{equation}\label{j_int}
\text{The square~\(J_{v(e)}\) lies in the interior of the cycles~\(D_i\) and~\(E_i.\)}
\end{equation}
\begin{eqnarray}
&&\text{The square~\(Q_e\) lies in the exterior of all the cycles~\(\{D_j\}_{1 \leq j \leq r}\)} \nonumber\\
&&\;\;\;\;\;\text{and in the exterior of all the cycles~\(\{E_j\}_{1 \leq j \leq r}.\)}\label{q_ext}
\end{eqnarray}
From the above two properties, we have that~\(e\) is a boundary edge in the graph~\(G(X)\) and so satisfies property~\((v)\) of Theorem~\ref{thm_out}.\\\\
\emph{Proof of~(\ref{j_int})}: From property~\((v),\) Theorem~\ref{thm_out} applied to the outermost boundary~\(\partial_0(X_1) = \cup_{1 \leq i \leq r} D_i,\) we have that the square~\(J_{v(e)}\) lies in the interior of the cycle~\(D_i.\) To see that~\(J_{v(e)}\) also lies in the interior of the cycle~\(E_i,\) we argue as follows. If \(1 \leq i \leq r-1,\) then~\(J_{v(e)}\) lies in the interior of the cycle~\(E_i = D_i.\) If~\(i = r,\) then~\(J_{v(e)}\) lies in the interior of the cycle~\(D_r,\) which in turn lies in the interior of the cycle~\(E_r\) (property~\((III.2)\)).\(\qed\)\\
\emph{Proof of~(\ref{q_ext})}: From property~\((v),\) Theorem~\ref{thm_out} applied to the outermost boundary~\(\partial_0(X_1) = \cup_{1 \leq i \leq r}D_i,\) we have that the square~\(Q_e\) lies in the exterior of all the cycles~\(\{D_i\}_{1 \leq i \leq r}.\) Since~\(E_j = D_j\) for~\(1 \leq j \leq r-1,\) the square~\(Q_e\) also lies in the exterior of the cycles~\(\{E_j\}_{1 \leq j \leq r-1}.\) The square~\(Q_e\) lies in the exterior of the cycle~\(D_r.\)

%we show that~\(Q_e \neq J_u\) and use property~\((III.1)\) that the interior of the cycle~\(E_r\) is the union of the interiors of the cycle~\(D_r\) and the square~\(J_u.\)

If~\(Q_e \neq J_u,\) then~\(Q_e\) also lies in the exterior of the cycle~\(E_r,\) since the interior of the cycle~\(E_r\) is the union of the interiors of the cycle~\(D_r\) and the square~\(J_u\) (property~\((III.1)\)). To see that the square~\(Q_e \neq J_u,\) we argue as follows. The square~\(J_u\) shares only the edge~\(e_{uv}\) with the cycle~\(D_r\) and does not any other vertex with~\(D_r\) (see property~\((b3.4)\)). Since~\(Q_e\) shares the edge~\(e \neq e_{uv}\) with~\(D_r,\) we have that~\(Q_e \neq J_u.\)~\(\qed\)

%and from property~\((III.2)\) and the fact that~\(Q_e \neq J_u,\) we also have that the square~\(Q_e\) also lies in the exterior of the cycle~\(E_r.\) %\(\qed\)

%we have that~\(Q_e \neq J_u.\)
%\begin{equation}\label{qe_n_ju}
%\text{The square~\(Q_e \neq J_u.\)}
%\end{equation}
%Thus~\(Q_e\) lies in the exterior of all cycles~\(\{E_j\}_{1 \leq j \leq r}\) and so~\(e\) is a boundary edge in the graph~\(G(X).\)
%\emph{Proof of~(\ref{qe_n_ju})}: The square~\(J_u\) shares only the edge~\(e_{uv}\) with the square~\(J_v\) and does not share a vertex with any other square~\(J_w, w\in X_1.\) But \(e_{uv} \in D_r\)  and \(e_{uv} \notin D_i, 1 \leq i \leq r-1\) (see property~\((III.1)\)). In particular, this means that the squares~\(J_{v(e)}\neq J_v\) and \(Q_e \neq J_u.\)\(\qed\)

%We prove that the square~\(Q_e \neq J_u\) as follows.

%We now consider the case when the edge~\(e \in E_r.\) If~\(e \in E_r \cap D_r = D_r \setminus e_{uv},\) then~\(e\) is adjacent to a square~\(J_{v(e)}\) with centre in~\(X_1\) and a square~\(Q_e\) with centre not in~\(X_1.\) The square~\(Q_e\) lies in the exterior of all the cycles~\(\{D_j\}_{1 \leq j \leq r}.\) As in the previous paragraph, the squares~\(J_e \neq J_v\) and~\(Q_e \neq J_u,\) since \(e \neq e_{uv}.\) Thus~\(Q_e\) does not have centre in~\(X\) and lies in the exterior of all the cycles~\(\{D_j\}_{1 \leq j \leq r}.\) This implies that the edge~\(e\) is a boundary edge in~\(G(X).\)

We now consider the case \(e \in E_r\setminus D_r = J_{u} \setminus e_{uv}.\) The edge~\(e\) is then the edge of the square~\(J_u\) and a square~\(Q_e\) lying in the exterior of the cycle~\(E_r.\) If the square~\(Q_e\) lies in the interior of some cycle~\(E_j = D_j, 1 \leq j \leq r-1,\) then the square~\(J_u\) would share an edge with the cycle~\(D_j,\) a contradiction to property~\((b3.4).\) Thus~\(Q_e\) lies in the exterior of the cycles~\(\{E_j\}_{1 \leq j \leq r}\) and so property~\((v),\) Theorem~\ref{thm_out} is also satisfied in this case.

 %As in the previous paragraph, the square \(Q_e \neq J_v\) since \(e \neq e_{uv}.\) But apart from~\(J_v,\) the square~\(J_u\) does not share an endvertex with any other square with centre in~\(X.\) Thus \(Q_e\) does not have centre in~\(X\) and lies in the exterior of all cycles~\(\{E_j\}_{1 \leq j \leq r}.\) We again obtain that~\(e\) is a boundary edge in~\(G(X).\) This proves that property~\((v),\) Theorem~\ref{thm_out} holds for all edges in~\(\cup_{1 \leq j \leq r} E_j.\)

%//CHNGE FROM DOWN AND TO SEE ABV CRFLLY +etC...

It remains to see that property~\((i)\) holds, i.e., every edge in~\(\cup_{1 \leq j \leq r} E_j\) is an outermost boundary edge in the graph~\(G(X).\) We assume other wise and arrive at a contradiction. Suppose that some edge~\(e \in \cup_{1 \leq j \leq r} E_j\) is not an outermost boundary edge in the graph~\(G(X).\) We first consider the case that~\(e \notin J_u.\) Using property~\((v)\) proved above, the edge~\(e\) is adjacent to a square~\(J_{v(e)}\) with centre in~\(v(e) \in X_1 = X \setminus \{u\}\) and a square~\(Q_e\) with centre not in \(X.\) Let~\(E_e\) be the outermost boundary cycle for the square~\(J_{v(e)}\) in the graph~\(G(X).\) %We refer to Lemma~\(3\) of Ganesan~(2017) for properties of the outermost boundary cycle.

%//WRTE FROM BELOW +etC...

As before, we must have~\(E_e \neq D_i\) and we state the related properties.\\
\((III.4)\) The square~\(J_{v(e)}\) lies in the interior of the cycle~\(D_i.\) Every edge of the cycle~\(D_i\) either belongs to~\(E_e\) or lies in the interior of~\(E_e.\) Also \(E_e \neq D_i\) and so at least one edge of~\(E_e\) lies in the exterior of~\(D_i.\)\\
\((III.5)\) At least one edge of \(E_e \setminus D_i\) belongs to the square~\(J_u\) and the square~\(J_u\) lies in the interior of the cycle~\(E_e.\)\\
\((III.6)\) The square~\(J_{v(e)} \neq J_u\) and the cycle~\(E_e\) contains both~\(J_{v(e)}\) and~\(J_v\) in its interior.

\emph{Proof of \((III.4)-(III.6)\)}: The proof of~\((III.4)-(III.5)\) is the same as before. To prove~\((III.6),\) we argue as follows. The square \(J_{v(e)} \notin \{J_{u}\}\) since the edge~\(e \) does not belong to~\(J_u.\) %or~\(J_v.\)

By definition, the square~\(J_{v(e)}\) lies in the interior of the cycle~\(E_e.\) Suppose that the square~\(J_v\) lies in the exterior of~\(E_e.\) Since the square~\(J_u\) lies in the interior of~\(E_e\) (property~\((III.5)\)), we then have that the edge~\(e_{uv}\) common to~\(J_u\) and~\(J_v\) belongs to~\(E_e.\) Merging~\(E_e\) and~\(J_v\) we then get a bigger cycle~\(C_e\) in the graph~\(G(X)\) containing~\(E_e\) in its interior, a contradiction to the fact that~\(E_e\) is the outermost boundary cycle in the graph~\(G(X)\) and satisfies property~\((c),\) Lemma~\ref{outer}. Thus the square~\(J_v\) also lies in the interior of the cycle~\(E_e.\)~\(\qed\)

%//done...CHK ABOVE +eTC..

%From property~\((III.4),\) some edge~\(f\) of the cycle~\(E_e\) belongs to the square~\(J_u\) and lies in the exterior of the cycle~\(D_i.\) The edge~\(f\) necessarily belongs to~\(J_u\) since~\(D_i\) is the outermost boundary cycle in the graph~\(G(X_1)\) formed by the squares with centres in~\(X_1 = X \setminus \{u\}.\)

%Suppose that the cycle~\(E_e\) contains only the two squares~\(J_u\) and~\(J_{v(e)}\) in its interior. This means that~\(J_u\) shares an edge with the square~\(J_{v(e)}.\)  But in this case~\((III),\) the square~\(J_u\) shares only an edge with~\(J_v\) and does not share a vertex with any other square in~\(\{J_w\}_{w \neq v}.\) Since~\(J_{v(e)} \neq J_v,\) (property~\((III.6)\)), this is a contradiction.

\begin{figure}[tbp]
\centering
%\fbox{
\includegraphics[width=2.5in, trim= 140 470 200 110, clip=true]{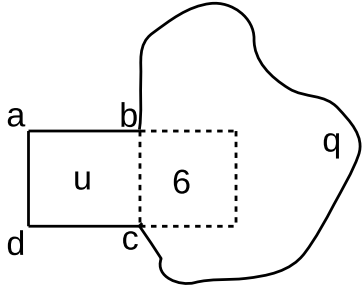}
%}
\caption{The squares~\(J_{q_1} = J_v\) (labelled \(6\)) and~\(J_u\) (labelled \(u\)) contained in the interior of the outermost boundary cycle~\(E_e\) represented by the wavy curve~\(abqcda.\)}
\label{fig_q1q2_ext}
\end{figure}

Using properties~\((III.4)-(III.6),\) we arrive at a contradiction. Let~\(E_e = (g_1,\ldots,g_t)\) and let~\(g_1 = f\) be the edge belonging to the square~\(J_u\) represented by the edge~\(ad\) in Figure~\ref{shr_edge}. Arguing as in Case~\((II),\) we obtain that the square labelled~\(6\) is the square~\(J_v\) sharing an edge with the square~\(J_u\) and the remaining edges~\(ad,ab\) and~\(cd\) belong to the cycle~\(E_e.\) We illustrate this in Figure~\ref{fig_q1q2_ext}, where we have shown only the squares labelled~\(u\) and~\(6\) representing the squares~\(J_u\) and~\(J_v,\) respectively. The outermost boundary cycle~\(E_e\) is the wavy curve~\(abqcda\) and the edge~\(bc\) is the edge~\(e_{uv}\) common to the squares~\(J_u\) and~\(J_v.\)

%By property~\((b),\) Lemma~\(3\) of Ganesan~(2017), every edge in~\(E_e\) is the edge of some square with centre in~\(X.\) Let \(i_1 \geq 2\) be the smallest index such that the edge~\(g_{i_1}\) belongs to a square~\(J_{q_1}, q_1 \in X_1 = X\setminus \{u\}.\) Similarly let~\(i_2 \leq t\) be the largest index such that the edge~\(g_{i_2}\) belongs to a square~\(J_{q_2}, q_2 \in X_1 = X \setminus \{u\}.\) The square~\(J_u\) shares an endvertex with both~\(J_{q_1}\) and~\(J_{q_2}.\)

%If \(q_1 \neq q_2,\) we obtain a contradiction, since the square~\(J_u\) shares only an edge with the square~\(J_v, v \in X_1\) and does not share an endvertex with any other square with centre in~\(X_1.\) If~\(q_1 = q_2,\) we proceed as follows. First, arguing as in the proof of~(\ref{shr_edge}), we have that the squares~\(J_{q_1}\) and~\(J_u\) share an edge and the edges~\(ab, ad\) and~\(dc\) of the square~\(J_u\) (see Figure~\ref{fig_q1q2}) form consecutive edges in the cycle~\(E_e.\) In particular, we have that~\(q_1 = v,\) the edge \(bc\) is the edge \(e_{uv} \in D_r\) and the square \(J_v\) lies in the interior of the cycle~\(E_e.\)

 %To argue CRFLLY THAT BX AND CY BELONG ETC... to see this carefully...The edges  \(bx\) and \(cy\) also belong to \(D_e\) and we argue as follows. If \(bx\) does not belong to \(D_e,\) then \(bx\) lies in the interior of \(D_e\) and the square lying above \(J_v\) also lies in the interior of \(D_e.\)

Consider now the graph \(D_{mod} = (E_e\setminus (ab,ad,dc)) \cup \{bc\}\) obtained by removing the subpath formed by the edges \(ab,ad\) and \(dc\) and adding the edge \(bc.\) The graph \(D_{mod}\) is also a cycle and consists only of edges belonging to squares with centres in \(X_1.\) By construction, the graph~\(D_{mod}\) is a cycle in the graph \(G(X_1)\) and we also have the following property.
\begin{eqnarray}
&&\text{The edge~\(bc = e_{uv}\) and the cycle~\(D_{mod} = D_r,\) the outermost} \nonumber\\
&&\;\;\;\text{boundary cycle in~\(\partial_0(X_1)\) containing the square~\(J_v\) in its interior.} \nonumber\\
\label{d_mod}
\end{eqnarray}
\emph{Proof of~(\ref{d_mod})}: From the discussion in the previous paragraph, we have that the edge~\(bc = e_{uv}.\) To prove the second part we argue as follows. By construction the cycle~\(D_{mod}\) is a cycle in the graph~\(G(X_1)\) and contains the square~\(J_v\) in its interior. Therefore~\(D_{mod}\) lies in the interior of the cycle~\(D_r,\) since the square~\(J_v\) lies in the interior of the cycle~\(D_r\) and~\(D_r\) is the outermost boundary cycle containing~\(J_v\) in its interior (property~\((e)\) following Lemma~\ref{outer}).

If the cycle~\(D_{mod} \neq D_r,\) then the merging of~\(D_r\) with the square~\(J_u\) would result in a cycle~\(C_e \neq E_e\) which contains~\(E_e\) in its interior. We recall that the square~\(J_{v(e)}\) also lies in the interior of~\(D_{mod}\) and therefore also lies in the interior of~\(C_e.\) But this is a contradiction since the cycle~\(E_e\) is the outermost boundary cycle containing the square~\(J_{v(e)}\) in its interior and so satisfies property~\((c),\) Lemma~\ref{outer}.~\(\qed\)

From~(\ref{d_mod}), we obtain that~\(E_{e} = E_r,\) since the cycle~\[E_r = (D_r \setminus \{e_{uv}\}) \cup (J_u \setminus \{e_{uv}\})\] is also obtained by merging~\(D_r\) with the square~\(J_u.\) From property~\((d)\) following Lemma~\ref{outer}, we obtain that~\(e\) is an outermost boundary edge in the graph~\(G(X),\) a contradiction.~\(\qed\)

It remains to see that every edge in~\(E_r \setminus D_r = J_u \setminus \{e_{uv}\}\) is also an outermost boundary edge in the graph~\(G(X).\) Suppose~\(e \in J_u\) is not an outermost boundary edge in~\(G(X)\) and let~\(E_u\) be the outermost boundary cycle in~\(G(X),\) containing the square~\(J_u\) in its interior. Arguing as in Case~\((II),\) we have the following properties.\\
\((III.6)\) The cycle~\(E_{u} \neq J_u\) and at least one edge of the cycle~\(E_u\) belongs to the square~\(J_u.\)\\
\((III.7)\) The cycle~\(E_u\) contains~\(J_u\) and at least one other square~\(J_w, w \in X_1\) in its interior.\\
%//PRKVMM +eTC...

Using~\((III.6)-(III.7)\) and an analogous analysis as following~\((III.3)-(III.5),\) we obtain a contradiction and so every edge in~\(J_u \setminus\{e_{uv}\}\) is also an outermost boundary edge in the graph~\(G(X).\) Finally, every other edge of~\(G(X)\) not present in~\(\cup_{1 \leq j \leq r} E_j\) necessarily lies in the interior of one of the cycles in~\(\{E_j\}.\) This proves property~\((i)\) and therefore~\((III.2).\)~\(\qed\)

\subsubsection*{4. Outermost boundary \(\partial_0(L_1)\) of a sub-\(S-\)cycle \(L_1 \subset L\)}
We recall the iterative process followed in obtaining the sequence \(L\) of vacant squares from the outermost boundary \(\partial_0(L) = D_{fin} = (f_1,\ldots,f_r)\) and suppose that \(i_2\) is the ``first time" we obtain a \(S-\)cycle; i.e., \[i_2 = \min\{i \geq 2 : (Z_i,Z_{i-1},\ldots,Z_j) \text{ form an }S-\text{cycle for some } j < i\}.\] Let \(i_1\) be the largest index such that \((Z_{i_2},Z_{i_2-1},\ldots,Z_{i_1})\) is a \(S-\)cycle.

%We now use \(D_{fin}\) to obtain the star connected \(S-\)cycle \(H_{out}.\) Letting \(D_{fin} = (f_1,...,f_r),\) there exists a unique square \(Z_1 \in \Lambda^+_1\) that has edge \(f_1.\) Because, if two squares \(Z_1\) and \(Z'_1\) in \(\Lambda_1^+\) share the edge \(f_1,\) then necessarily one of them is in the exterior of \(D_{fin}.\) This contradicts property \((a3)\) of \(D_{fin}.\)

For convenience let \(L_1 = (Z_{i_1},Z_{i_1+1},\ldots,Z_{i_2})\) be the \(S-\)cycle and suppose~\(L_1\) does not contain all squares of the sequence~\(L.\) The sequence \(L_1\) is a star connected component and let \(\partial_0(L_1)\) denote its corresponding outermost boundary. We use the following properties of~\(\partial_0(L_1)\) in the remaining steps.\\
\((a4.1)\) There are at least three distinct squares in the sequence~\(L_1.\) The outermost boundary~\(\partial_0(L_1)\) of the sequence~\(L_1\) is a cycle and is not equal to the outermost boundary \(D_{fin} = \partial_0(L)\) of the sequence \(L\) obtained in Step~\(1.\)\\
\((a4.2)\) Every vacant square~\(Z_j, i_1 \leq j \leq i_2\) of the sequence~\(L_1\) is contained in the interior of the cycle~\(\partial_0(L_1)\) and every edge~\(e \in \partial_0(L_1)\) belongs to a unique vacant square~\(Z_{j(e)}, i_0 \leq j(e) \leq i_1\) of the sequence~\(L_1.\) Moreover, the square~\(Z_{j(e)}\) lies in the interior of~\(\partial_0(L_1).\)\\
\((a4.3)\) Suppose two edges \(f_1, f_2\in \partial_0(L_1)\) belong to some vacant square \(Z_j \in L_1\) but do not share an endvertex; i.e., they are ``opposite" edges of~\(Z_j.\) If~\(g_1\) and~\(g_2\) denote the remaining edges of \(Z_j,\) then either \((f_1,g_1,f_2)\) or \((f_1,g_2,f_2)\) is a subpath of~\(\partial_0(L_1).\)\\

We recall that \(C^+(0)\) is the plus connected occupied component containing the origin and the outermost boundary \(\partial^+_0\) of~\(C^+(0)\) is a single cycle denoted by~\(\partial^+_0 = (e_1,\ldots,e_t)\) (see (\ref{do_plus})). The relation between~\(C^+(0),\partial_0^+\) and the cycle~\(\partial_0(L_1)\) is as follows.\\
\((b4.1)\) All (occupied) squares in \(C^+(0)\) lie in the exterior of the cycle~\(\partial_0(L_1).\)\\
\((b4.2)\) The outermost boundary cycle~\(\partial_0^+\) of the component~\(C^+(0)\) and the cycle~\(\partial_0(L_1)\) have mutually disjoint interiors.\\
\((b4.3)\) There are three distinct edges \(\{e_{j_k}\}_{1 \leq k \leq 3} \subset \partial_0^+\) and three distinct squares~\(\{Z_{j_k}\}_{1 \leq k \leq 3} \subset L_1\) with the following property. For \(k = 1,2,3,\) the edge~\(e_{j_k}\) belongs to~\(\partial_0(L_1)\) and is an edge of the square~\(Z_{j_k}.\)\\

The final property says that there are at least three distinct edges belonging to both~\(\partial_0^+\) and \(\partial_0(L_1).\) Moreover, these edges belong to distinct vacant squares of~\(L_1.\)

%//CHECK HOW TO DO WITH THREE AND GET CONTRA AT END..PRKVM...
%\((b3)\) Either all squares in \(C^+(0)\) lie in the interior of \(\partial_0(L_1)\) or all the squares in \(C^+(0)\) lie in the exterior of \(\partial_0(L_1).\)\\

%//WRITE B2 AS A SEPARATE RESULT ?? PARKKVM...
%//CHECK PROPERTIES BELOW AGAIN...SINCE RENUMBEREDD...ETC..

\emph{Proof of \((a4.1)-(a4.2)\)}: For the first statement of \((a4.1),\) we let~\(H_{L_1}\) be the subgraph of \(G_V\) obtained for the sequence \(L_1\) in the same manner as the subgraph~\(H_L.\) Since~\(H_{L_1}\) is a cycle in~\(G_V,\) there are at least three vertices in~\(H_{L_1}\) and so there are at least three distinct squares in~\(L_1.\)

For the second statement of \((a4.1),\)  we label the vacant squares in the sequence~\(L_1\) as with label~\(1\) and every other square with label~\(0.\) We then apply Theorem~\ref{thm_out} with ``occupied" replaced by label~\(1\) and ``vacant" replaced by label~\(0.\) From Theorem~\ref{thm6}, we obtain that the outermost boundary~\(\partial_0(L_1)\) of the \(S-\)cycle~\(L_1\) is a single cycle. To see that~\(D_{fin} \neq \partial_0(L_1),\) we recall from property~\((a2.1)\) above that the cycle~\(D_{fin}\) is the outermost boundary of the sequence~\(L = (Z_1,\ldots,Z_s).\) Moreover, by construction~\(D_{fin}\) has at least one edge from every \(Z_i, 1 \leq i \leq s\) (see property~\((b1.5)\) of~\(D_{fin}\) in Step~\(1\)). Since the sequence~\(L_1\) does not contain all the squares of the sequence~\(L,\) we obtain the second statement of \((a4.1).\)

%//DRAW FIG HERE...PRKVMM...

The first statement \((a4.2)\) follows from property~\((iv),\) Theorem~\ref{thm_out}. To see the uniqueness of~\(Z_{j(e)},\) suppose that there are two squares~\(Z_{j(e)}\) and~\(Z_q\) belonging to~\(L_1\) and containing~\(e\) as an edge. We then have that exactly one of~\(Z_{j(e)}\) or~\(Z_q\) lies in hte interior of~\(\partial_0(L_1)\) but not both. This contradicts the first statement of~\((a4.2).\)\(\qed\)

\emph{Proof of~\((a4.3)\)}: Suppose that two opposite edges of a vacant square~\(Z_k \in L_1\) belong to the outermost boundary cycle~\(\partial_0(L_1).\) We arrive at a contradiction assuming that the remaining two edges of~\(Z_k\) do not belong to~\(\partial_0(L_1).\)

%For convenience, we let \(h(t),h(r), h(b)\) and \(h(l)\) be the top, right, bottom and left edges of \(Z_k.\) Let \(z_{tr}\) be the vertex common to~\(h(t)\) and~\(h(r)\) and \(z_{rb}\) be the vertex common to \(h(r)\) and \(h(b).\) Similarly define the other two vertices \(z_{bl}\) and \(z_{lt}.\)
\begin{figure}[tbp]
\centering
%\fbox{
\includegraphics[width=1.6in, trim= 150 340 240 170, clip=true]{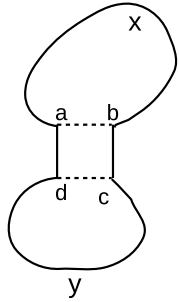}
%}
\caption{Opposite edges~\(ad\) and~\(bc\) of the vacant square~\(Z_k =abcd\) in the sequence~\(L_1\) belonging to the outermost boundary cycle~\(\partial_0(L_1) =axbcyda.\) The cycles~\(C_1 = axba\) and \(C_2 = cydc\) have mutually disjoint interiors and share an edge each, with the square~\(Z_k = abcd.\)}
\label{opp_edge_fig}
\end{figure}

%//DRAW FIG HERE...

%~\(Z_{j_2}\) lies in the interior of~\(C_2,\) then~\(P_{12}\) necessarily contains~\(Z_k\) and so \(j_1 < k < j_2.\)

The square~\(Z_k\) lies in the interior of the cycle~\(\partial_0(L_1)\) (property~\((a4.2)\)) and so the edges of~\(Z_k\) not in~\(\partial_0(L_1)\) lie in the interior of~\(\partial_0(L_1).\) This is illustrated in Figure~\ref{opp_edge_fig}, where the outermost boundary cycle~\(\partial_0(L_1)\) is denoted by the cycle~\(axbcyda.\) The vacant square~\(Z_k\) denoted by the square~\(abcd\) has opposite edges~\(ad\) and~\(bc\) in the cycle~\(\partial_0(L_1).\) The remaining edges~\(ab\) and~\(cd\) lie in the interior of~\(\partial_0(L_1).\)

The subpath~\(P_1 = axb \subset \partial_0(L_1)\) together with the edge~\(ab \notin \partial_0(L_1)\) forms a cycle~\(C_1 = axba.\) Similarly, the subpath~\(P_2 = dyc \subset \partial_0(L_1)\) and the edge~\(cd \notin \partial_0(L_1)\) form a cycle~\(C_2 = cydc.\) We have the following properties.\\
\((c4.1)\) The cycles~\(C_1\) and~\(C_2\) have non empty mutually disjoint interiors and do not share any vertex in common. Moreover, for~\(i = 1,2,\) the square~\(Z_k\) lies in the exterior of~\(C_i\) and shares exactly one edge with \(C_i.\) In Figure~\ref{opp_edge_fig}, the cycle~\(C_1\) shares the edge~\(ab\) with the square~\(Z_k = abcd\) and the cycle~\(C_2\) shares the edge~\(cd\) with~\(Z_k.\)\\
\((c4.2)\) If a square~\(Q\) lies in the interior of the cycle~\(\partial_0(L_1),\) then either~\(Q = Z_k\) or~\(Q\) lies in the interior of either~\(C_1\) or~\(C_2.\) In particular, for \(i = 1,2,\) the only square in~\(L_1\) lying in the exterior of \(C_i\) and sharing a vertex with~\(C_i\) is~\(Z_k.\) \\
\((c4.3)\) Each cycle~\(C_i,  i = 1,2\) contains at least one vacant square belonging to the \(S-\)cycle~\(L_1.\) Without loss of generality, suppose that the square~\(Z_{i_1} \in L_1\) lies in the interior of~\(C_1.\)\\
\((c4.4)\) All the squares~\(\{Z_j\}_{i_1 \leq j \leq k-1} \subset L_1 \) lie in the interior of~\(C_1\) and all the squares~\(\{Z_{j}\}_{k+1 \leq j \leq i_2}\) lie in the interior of~\(C_2.\)\\\\
\emph{Proof of~\((c4.1)-(c4.4)\)}:  Properties~\((c4.1)-(c4.2)\) are true by construction. We prove property~\((c4.3)\) as follows. Since the \(S-\)cycle~\(L_1\) contains at least three distinct squares~(property~\((a4.1)\)), at least one of the cycles~\(C_1\) or~\(C_2\) contains a (vacant) square of~\(L_1\) in its interior. Without loss of generality, we assume that the square~\(Z_{i_1}\) lies in the interior of~\(C_1.\)

The cycle~\(C_2\) has non empty interior and every edge of~\(C_2 \setminus \{cd\}\) belongs to the outermost boundary cycle~\(\partial_0(L_1).\) Fix an edge~\(e \neq cd\) of the cycle~\(C_2.\) From Theorem~\ref{thm_out}, property~\((v),\) the edge~\(e\) is the edge of a square~\(Z_j \in L_1.\) The square~\(Z_j\) lies in the interior of the cycle~\(\partial_0(L_1).\) From property~\((c4.2),\) the square~\(Z_j\) therefore either equals~\(Z_k\) or lies in the interior of the cycle~\(C_1\) or lies in the interior of the cycle~\(C_2.\)

The square~\(Z_j \neq Z_k\) since \(e \neq cd\) and the square~\(Z_k\) shares only the edge~\(cd\) with~\(C_2\) (property~\((c4.1)\)). Also the edge~\(e \in C_1\) and the cycles~\(C_1\) and \(C_2\) share no vertex in common~(property~\((c4.1)\)).  Thus~\(Z_j\) does not lie in the interior of~\(C_1\) and so lies in the interior of~\(C_2.\) This proves~\((c4.3).\)

%in the exterior of~\(C_2.\) Thus from property~\((c4.2),\) we have that~\(Z_j\) lies in the

%//...WRT MORE?? PRKVMM+ETC..

%Both squares~\(C_1\) and~\(C_2\) contain at least one square of the cycle~\(C_1.\) Without loss of generality, suppose that~\(Z_{i_1}\) lies in the interior of~\(C_1.\)

To prove property~\((c4.4),\) we use the following additional property.\\
\((c4.5)\) Let~\(P_{12} = (Z_{j_1},Z_{j_1+1},\ldots,Z_{j_2}), \) be a subpath of the \(S-\)cycle~\(L_1.\) If~\(Z_{j_1}\) lies in the interior of~\(C_1\) and~\(Z_{j_2}\) lies in the interior of~\(C_2,\) then~\(Z_k \in P_{12}\) and so \(\min(j_1,j_2) < k < \max(j_1,j_2).\) \\
\emph{Proof of~\((c4.5)\)}: We assume that \(j_1 < j_2\) and an analogous argument holds otherwise. Since the cycles~\(C_1\) and~\(C_2\) have mutually disjoint interiors~(property~\((c4.1)\)), the path~\(P_{12}\) crosses the cycle~\(C_1\) at some point, in the sense that some square \(Z_{y}\in P_{12}\) lies in the exterior of~\(C_1.\) From property~\((c4.1),\) we have that~\(Z_k = Z_y\) and this proves~\((c4.5).\) \(\qed\)

We prove~\((c4.4)\) by applying property~\((c4.5)\) repeatedly. First, applying~\((c4.5),\) we obtain that all squares~\(Z_{j}, i_1 \leq j \leq k-1\) lie in the interior of~\(C_1.\) Using property~\((c4.3),\) suppose that some square~\(Z_{r}, k < r \leq i_2\) belonging to~\(L_1\) lies in the interior of the cycle~\(C_2.\) If a square~\(Z_{j}, r < j \leq i_2\) of~\(L_1\) lies in the interior of~\(C_1,\) we apply property~\((c4.5)\) with \(j_1 = r\) and \(j_2 = j\) to get that \(r < k < j,\) a contradiction. Thus all squares~\(Z_j, r \leq j \leq i_2\) lie in the interior of~\(C_2.\) Finally, if some square~\(Z_j, k+1 \leq j \leq r-1\) belongs to the interior of \(C_1,\) then we again apply property~\((c4.5)\) to get that~\(j < k < r,\) a contradiction. This proves~\((c4.4).\)~\(\qed\)

Using property~\((c4.4)\) and the fact that cycles~\(C_1\) and \(C_2\) do not have any vertex in common (property~\((c4.1)\)), we must have that~\(Z_{i_1}\) and \(Z_{i_2}\) do not share any vertex in common. But this is a contradiction since~\(L_1\) is a star connected \(S-\)cycle and so the first and the last squares of~\(L_1,\) \(Z_{i_1}\) and \(Z_{i_2},\) are star connected; i.e., they share a vertex. We recall that we have obtained the above contradiction assuming that both the edges~\(ab\) and \(cd\) in Figure~\ref{opp_edge_fig} do not belong to the cycle~\(\partial_0(L_1).\) Thus either the edge~\(ab\) or~\(cd\) belongs to~\(\partial_0(L_1)\) and this proves~\((a4.3).\) \(\qed\)

%For the second statement, suppose that \(h(t), h(r)\) and \(h(b)\) belong to the outermost boundary cycle~\(\partial_0(L_1).\) We note that every vertex in the cycle \(\partial_0(L_1)\) has degree \(2\) and therefore the vertex \(z_{tr}\) common to \(h(t)\) and \(h(r)\) cannot be adjacent to any other edge of \(\partial_0(L_1).\) Similarly, the vertex \(z_{rb}\) common to \(h(r)\) and \(h(b)\) cannot contain any other edge of \(\partial_0(L_1).\) Thus the edges~\((h(t),h(r),h(b))\) form a subpath in~\(\partial_0(L_1).\)

%//SEE CRFLLY AGN ABV +eTC...

\emph{Proof of \((b4.1)-(b4.3)\)}: To see~\((b4.1),\) we first suppose there is an occupied square~\(S_{k_1} \in C^+(0)\) lying in the interior of the cycle~\(\partial_0(L_1)\) and another occupied square \(S_{k_2} \in C^+(0)\) lying in the exterior of~\(\partial_0(L_1).\) There is a plus connected \(S-\)path \(P = (S_{k_1} = J_1,J_2,\ldots,J_{t-1},J_t = S_{k_2}), J_i \in C^+(0)\) of occupied squares connecting~\(S_{k_1}\) and~\(S_{k_2}.\) This path~\(P\) necessarily crosses the cycle~\(\partial_0(L_1)\) in the sense that there is an integer~\(m \geq 1\) such that the following three statements hold. The square~\(J_m \in P\) lies in the interior of the cycle \(\partial_0(L_1).\) The square \(J_{m+1} \in P\) lies in the exterior of~\(\partial_0(L_1).\) The edge~\(e\) common to~\(J_m\) and~\(J_{m+1}\) belongs to~\(\partial_0(L_1).\) But by property \((a4.2),\) this means that the square~\(J_{m}\) is vacant, a contradiction.

From the above paragraph, we therefore have that either all squares of~\(C^+(0)\) lie in the exterior of the cycle~\(\partial_0(L_1)\) or all squares of~\(C^+(0)\) lie in the interior of~\(\partial_0(L_1).\) Suppose that all occupied squares of~\(C^+(0)\) lie in the interior of~\(\partial_0(L_1).\) We arrive at a contradiction by showing that~\(\partial_0(L_1) \neq D_{fin}\) satisfies the properties~\((b1.1)-(b1.3)\) of the cycle~\(D_{fin}\) (see Step~\(1\)).

Property~\((b1.1)\) holds for~\(\partial_0(L_1)\) since~\(\partial_0(L_1)\) contains only edges of the (vacant) squares~\(\{Z_{j}\}_{i_0 \leq j \leq i_1}\subset \Lambda^+\) (see property~\((a4.2)\) above). To see that the cycle~\(\partial_0(L_1)\) also satisfies property~\((b1.2),\) we argue as follows. By assumption, all the occupied squares of the component~\(C^+(0)\) lie in the interior of the cycle~\(\partial_0(L_1)\) and by Theorem~\ref{thm_out}, property~\((iv),\) every edge of the outermost boundary~\(\partial^+_0\) belongs to some occupied square in~\(C^+(0).\) Thus every edge of~\(\partial^+_0\) either belongs to or lies in the interior of~\(\partial_0(L_1).\)

Suppose~\(e \in \partial^+_0\) also belongs to the cycle~\(\partial_0(L_1)\) and let \(V_e\) be the square containing the edge~\(e\) lying in the interior of~\(\partial_0(L_1).\) Applying Theorem~\ref{thm_out}, property~\((iv)\) to~\(\partial_0^+,\) we obtain that \(V_e\) is an occupied square of~\(C^+(0).\) Applying Theorem~\ref{thm_out}, property~\((iv)\) to~\(\partial_0(L_1),\) we obtain that~\(V_e\) is a vacant square belonging to the sequence~\(L_1.\) This is a contradiction and so \(e\) lies in the interior of~\(\partial_0(L_1)\) and so~\(\partial_0(L_1)\) also satisfies property~\((b1.2)\) of cycle~\(D_{fin}.\)

To see that the cycle~\(\partial_0(L_1)\) satisfies property~\((b1.3),\) it is enough to see that all squares in~\(\Lambda^+\) lie in the interior of~\(\partial_0(L_1).\) We recall that~\(\Lambda^+\) is the set of all vacant squares sharing an edge with some occupied square of the component~\(C^+(0).\) We assume otherwise and arrive at a contradiction. Suppose some square~\(Y_j \in \Lambda^+\) lies in the exterior of~\(\partial_0(L_1).\) By definition, the square~\(Y_j\) contains an edge~\(e_j \in \partial_0^+\) and the edge~\(e_j\) also belongs to an occupied square~\(A_j \in C^+(0)\) (Theorem~\ref{thm_out}, property~\((v)\)).

If the square~\(A_j\) lies in the interior of~\(\partial_0(L_1),\) then the edge~\(e_j\) belongs to~\(\partial_0(L_1)\) and this contradicts property~\((a4.2).\) Thus~\(A_j\) lies in the exterior of~\(\partial_0(L_1)\) and this contradicts the assumption that all occupied squares of~\(C^+(0)\) lie in the interior of~\(\partial_0(L_1).\) So all squares in~\(\Lambda^+\) lie in the interior of~\(\partial_0(L_1)\) and therefore~\(\partial_0(L_1)\) also satisfies property~\((b1.3)\) of the cycle~\(D_{fin}.\) Since \(D_{fin} \neq \partial_0(L_1)\) (see property~\((a4.1)\) above) we obtain a contradiction to the uniqueness property~\((b1.4)\) of the cycle~\(D_{fin}\) and so all occupied squares of the component~\(C^+(0)\) lie in the exterior of~\(\partial_0(L_1).\) This proves property~\((b4.1).\)

%Suppose further that the square~\(Y_j\) shares an edge~\(e\) with~\(\partial_0(L_1).\)

%does not share any edge with~\(\partial_0(L_1)\) since if \(e \in \partial_0(L_1),\) then the square containing~\(e\) and lying in the exterior of~\(\partial_0(L_1)\) does not belong to~\(L_1.\) The above statement is true by applying Theorem~\ref{thm_out}, property~\((iv)\)~\(\partial_0(L_1).\)

%From the above paragraph, we obtain that all edges of the square~\(Y_j\) lie in the exterior of the cycle~\(\partial_0(L_1).\) But by definition, the square~\(Y_j\) contains an edge~\(e_j \in \partial_0^+\) and we have proved above that~\(\partial_0(L_1)\) satisfies property~\((b1.2)\) of the cycle~\(D_{fin};\) i.e., all edges of~\(\partial_0^+\) lie in the interior of~\(\partial_0(L_1).\) This leads to a contradiction and so all squares in~\(\Lambda^+\) lie in the interior of~\(\partial_0(L_1)\) and therefore~\(\partial_0(L_1)\) also satisfies property~\((b1.3)\) of the cycle~\(D_{fin}.\) Since \(D_{fin} \neq \partial_0(L_1)\) (see property~\((a4.1)\) above) we obtain a contradiction to the uniqueness property~\((b1.4)\) of the cycle~\(D_{fin}\) and so all occupied squares of the component~\(C^+(0)\) lie in the exterior of~\(\partial_0(L_1).\) This proves property~\((b4.1).\)

We prove property~\((b4.2)\) as follows. We first have from property \((b4.1)\) above that every occupied square of \(C^+(0)\) lies in the exterior of the cycle~\(\partial_0(L_1).\) Therefore every edge in the outermost boundary \(\partial_0^+\) either belongs to or lies in the exterior of~\(\partial_0(L_1).\) If \(\partial_0^+\) and \(\partial_0(L_1)\) do not have mutually disjoint interiors, then~\(\partial_0(L_1)\) lies in the interior of~\(\partial_0^+;\) i.e., every edge of~\(\partial_0(L_1)\) either belongs to or lies in the interior of~\(\partial^+_0.\)

%//WRT BELOW +eTC CARFLLY...

This leads to a contradiction as follows. From property \((b1.2)\) of the cycle~\(D_{fin}\) (see Step~\(1\)), we recall that every edge in the outermost boundary~\(\partial_0^+\) lies in the interior of~\(D_{fin}.\) In particular, no edge of~\(\partial_0^+\) belongs to~\(D_{fin}.\) From the above paragraph, we therefore have that every edge in~\(\partial_0(L_1)\) also lies in the interior of~\(D_{fin}.\) From Theorem~\ref{thm_out}, property~\((ii)\) applied to the outermost boundary cycle~\(\partial_0(L_1),\) we have that all squares~\(\{Z_i\}_{i_1 \leq i \leq i_2} \subset L_1\) lie in the interior of~\(\partial_0(L_1).\) In particular, this means that every edge of the square~\(Z_{i_1}\) is contained in the \emph{interior} of \(D_{fin}.\) This contradicts the property~\((b1.5)\) of~\(D_{fin}\) that every edge in the sequence \(L = (Z_1,\ldots,Z_s)\) contains at least one edge in~\(D_{fin}.\) Thus the cycles \(\partial_0^+\) and~\(\partial_0(L_1)\) have mutually disjoint interiors and this proves~\((b4.2).\)

%//SEE FROM BELW++ETC...

%//WRT CRFLLY BELOW...SHW THAT THREE EDGES OF THREE DISTINCT SQUARES....

It remains to prove property \((b4.3).\) We recall that every vacant square \(Z_k, i_1 \leq k \leq i_2\) of the \(S-\)cycle~\(L_1\) belongs to~\(\Lambda^+\) and so by definition, shares an edge~\(f_{k}\) with the outermost boundary cycle~\(\partial_0^+\) of the plus connected component~\(C^+(0).\) From property~\((a4.1),\) we have that there are at least three distinct squares in~\(L_1.\) Thus~\(f_{i_1},f_{i_1+1}\) and \(f_{i_1+2}\) exist. We first prove that each~\(f_{j}, i_1 \leq j \leq i_1+2\) also belongs to the cycle~\(\partial_0(L_1)\) and then show that they are distinct.

First, applying Theorem~\ref{thm_out}, property~\((iv),\) to the outermost boundary cycle~\(\partial_0^+,\) we have that the edge~\(f_{i_1}\) is also the edge of an occupied square~\(A_1 \in C^+(0)\) lying in the interior of~\(\partial_0^+.\) So~\(f_{i_1}\) cannot lie in the interior of the cycle~\(\partial_0(L_1)\) since if it does, then the occupied square~\(A_1 \in C^+(0)\) also would also lie in the interior of~\(\partial_0(L_1),\) contradicting~property~\((b4.1).\) Thus~\(f_{i_1} \in \partial_0(L_1)\) and similarly~\(f_{i_1+1}\) and \(f_{i_1+2}\) also belong to~\(\partial_0(L_1).\) %for~\(1 \leq k \leq 3.\)

To see that the edges~\(\{f_{j}\}_{i_1 \leq j \leq i_1+2}\) are all distinct, we argue as follows. Suppose~\(f_{i_1}\) belongs to two distinct squares~\(Z_{i_1}\) and~\(Z_k,\) both belonging to~\(L_1.\) Since~\(f_{i_1} \in \partial_0(L_1),\) exactly one of the squares in~\(\{Z_{i_1},Z_k\}\) lies in the interior of~\(\partial_0(L_1)\) and the other square lies in the exterior of~\(\partial_0(L_1).\) This is a contradiction to property~\((a4.2)\) that all squares of the sequence~\(L_1\) lie in the interior of~\(\partial_0(L_1).\) Thus the edges~\(\{f_{j}\}_{i_1 \leq j \leq i_1+2}\) are all distinct and this proves~\((b4.3).\)\(\qed\)

%The edges~\(\{e_{j_k}\}_{1 \leq k \leq 3}\) are also distinct,

%We now see that each edge~\(e_{j_k},1 \leq k \leq 3,\) belongs to the outermost boundary~\(\partial_0(L_1).\)  Since every vacant square in the sequence~\(L_1\) lies in the interior of~\(\partial_0(L_1),\) we have that the edge~\(e_{j_k},1 \leq k \leq 3\) either belongs to~\(\partial_0(L_1)\) or lies in the interior of~\(\partial_0(L_1).\) If \(e_{j_k}\) lies in the interior of~\(\partial_0(L_1),\) then

%To see that the edges are distinct, we suppose that edge~\(e_{j_1}\) is common to both~\(Z_1\) and~\(Z_j\) for some index~\(i_1 \leq j \leq i_2.\) Since the outermost boudnary The edge \(e_{j_k}\) is common to the square~\(Z_k\) lying in the exterior of~\(\partial_0^+\) and an occupied square~\(W_k \in C^+(0)\) contained in the interior of \(\partial_0^+.\) In particular, this means that the edges \(\{e_{j_k}\}_{i_0 \leq k \leq i_1}\) are all distinct. Since the \(S-\)cycle \(L_1\) contains at least three squares (property~\((a4.1)\)), we obtain \((b4.3).\) ?? WHY ??? CHK PRKVMM+eTC...\(\qed\)

\subsubsection*{5. Merging cycles \(\partial_0(L_1)\) and \(\partial_0^+\) to get cycle \(\partial_{temp}\)}
In this step of our argument, we merge the cycles~\(\partial_0(L_1)\) and~\(\partial^+_0.\) From property~\((b4.3),\) there is at least one edge common to the~\(\partial_0(L_1)\) and~\(\partial_0^+.\) Since the cycles~\(\partial_0(L_1)\) and \(\partial^+_0\) have mutually disjoint interiors (property~\((b4.2)\)), we have from Theorem~\ref{thm8} that the merged cycle~\(\partial_{temp}\) is the union of two unique bridges~\(A_{out}\) and~\(B_{out}\) such that~\(A_{out} \subset \partial^+_0\) is a bridge for~\(\partial_0(L_1)\) and~\(B_{out} \subset \partial_0(L_1)\) is a bridge for~\(\partial^+_0.\)

%//WE WANT TO GET THAT THERE IS A SQUARE CONTAINING ALL ITS **EDGES** IN THE INTERIOR OF DTEMP?? PRKVMM +etC...

%SEE BELOW +ETC.. CRFLY...

%// TRY TO RELATE BIN WITH PAR ZERO PLUS CAP PAR ZERO L1

We use the cycle \(\partial_{temp}\) to arrive at a contradiction in the next step. In this step, we state and prove the properties of the cycle~\(\partial_{temp}\) needed for future use.\\
\((a5.1)\) The interior of the cycles~\(\partial_0^+\) and~\(\partial_0(L_1)\) is contained in the interior of the cycle~\(\partial_{temp} = A_{out} \cup B_{out}.\) Every vacant square of the sequence~\(L_1\) and every occupied square of the component~\(C^+(0)\) lies in the interior of~\(\partial_{temp}.\) Every edge in~\(\partial_{temp}\) belongs to either~\(A_{out} \in \partial_0^+\) or~\(B_{out} \subset \partial_0(L_1)\) but not both.\\
\((a5.2)\) Every edge in the path \(B_{in} = \partial_0(L_1) \setminus B_{out}\) lies in the interior of the cycle~\(\partial_{temp}.\) At most two edges of~\(B_{in}\) contain an endvertex in~\(\partial_{temp}\) and every other edge of~\(B_{in}\) has both its endvertices in the interior of~\(\partial_{temp}.\)\\
The above  two properties are used to prove the following crucial property used in the next step.\\
\((a5.3)\) There are three edges~\(h_1,h_2,h_3 \in B_{in}\) satisfying the following properties. For~\( i =1,2,3,\) the edge~\(h_i\) belongs to~\(\partial_0^+ \cap \partial_0(L_1)\) and is the edge of a vacant square~\(Z_{j_i} \in L_1.\) The indices~\(j_1 \neq j _2 \neq j_3\) and the edge~\(h_1\) has \emph{both} its endvertices in the interior of the cycle~\(\partial_{temp}.\)\\
%There are exactly two edges in~\(B_{in}\) that share an end vertex with~\(B_{out}\) (and therefore with~\(\partial_{temp}\)). Every other edge of \(B_{in}\) has both its endvertices in the interior of \(\partial_{temp}.\)\\

%In particular, the edge~\(e\) cannot belong to~\(\partial_0(L_1) \cap \partial_0^+.\)

%If \(e \in \partial_0(L_1),\) then~\(e\) is the edge of a vacant square in~\(L_1.\) If an edge~\( h \in \partial_0^+ \cup \partial_0(L_1)\)  but does not belong to the cycle~\(\partial_{temp},\) then~\(h\) lies in the interior of~\(\partial_{temp}.\)\\

%TRY TO PUT ABOVE INSIDE...

%WRT BETTER ETC +...

%From property~\((b4.3),\) we recall that there are at least three distinct edges belonging to both~\(\partial_0^+\) and~\(\partial_0(L_1).\) We have the following properties.\\
%\((a5.3)\) Every edge in \(\partial_0^+ \cap \partial_0(L_1)\) belongs to~\(B_{in}\) and therefore lies in the interior of the cycle~\(\partial_{temp}.\) There is an edge \(g_0 \in \partial_0^+ \cap \partial_0(L_1) \subset B_{in}\) that has \emph{both} its endvertices in the interior of~\(\partial_{temp}.\)\\

%contains at least three edges with the property that each edge belongs to a \emph{distinct} vacant square of the sequence \(L_1.\)

%SEE IF WE NEED MORE PROPERTIES...COULD WE PUT BRIDGES ETC..INSIDE THE PROOF??PRKVMM...
%ARRANGE LATER...

%\((a5.2)\) If \(e \in \partial_{temp}\) is an edge, then either \(e \in A_{out} \subset \partial_0^+\) or~\(e \in B_{out} \subset \partial_0(L_1)\) but not both.\\

\emph{Proof of \((a5.1)-(a5.3)\)}:  Property~\((a5.1)\) and the first statement of property~\((a5.2)\) follow from the properties of the merged cycle in Theorem~\(3,\) Ganesan~(2017). To see the second statement, let \(h \in B_{in}\) be any vertex that shares an endvertex with~\(\partial_{temp} = A_{out} \cup B_{out}.\) All four paths~\(A_{in}, A_{out}, B_{in}\) and~\(B_{out}\) have common endvertices and have no other vertex in common. Therefore the edge~\(h\) necessarily contains an endvertex of the path~\(B_{in}.\) Since~\(B_{in}\) has two endvertices, there are at most two possible choices for~\(h.\) Every other edge of~\(B_{in}\) has both its endvertices in the interior of~\(\partial_{temp}.\) This proves~\((a5.2).\)

%But from property \((a5.6)\) of the cycle~\(\partial_0(L_1)\) (see Step~\(5\)), we have that there are at least three distinct edges in~\(\partial_0^+ \cap \partial_0(L_1).\) Thus at least one edge in~\(\partial_0^+ \cap \partial_0(L_1)\) has both its endvertices in the interior of~\(\partial_{temp}.\)

To prove~\((a5.3),\) we first prove the following property.\\
\((a5.4)\) If an edge~\(e \in \partial_0^+ \cap \partial_0(L_1),\) then~\(e \in B_{in}.\) \\
\emph{Proof of~\((a5.4)\)} consider an edge \(e \in \partial_0^+ \cap \partial_0(L_1).\) Applying Theorem~\ref{thm_out}, property~\((iv)\) to \(e \in \partial_0^+,\) we have that the edge~\(e\) is the edge of an occupied square~\(K_e\) of the plus connected component~\(C^+(0)\) and is also the edge of a vacant square~\(Y_e.\) Applying Theorem~\ref{thm_out}, property~\((iv)\) to~\(e \in \partial_0(L_1),\) we obtain that the vacant square~\(Y_e\) belongs to the sequence~\(L_1.\) Therefore if \(e \in \partial_{temp},\) then necessarily one of~\(K_e\) or~\(Y_e\) lies in the exterior of the cycle~\(\partial_{temp},\) contradicting property~\((a5.1).\)  Thus~\( e \in B_{in}.\)\(\qed\)

The first statement of~\((a5.3)\) follows from property~\((b4.3)\) of Step~\(4\) which states that there are at least three distinct edges~\(h_1,h_2,h_3 \in \partial_0^+ \cap \partial_0(L_1),\) each adjacent to a distinct vacant square of the sequence~\(L_1.\) From property~\((a5.4),\) the edges~\(h_i, 1 \leq i \leq 3\) belong to~\(B_{in}.\) From the final statement of property~\((a5.2),\) we have that at least one of edges in~\(\{h_1,h_2,h_3\}\) contains both its endvertices in the interior of~\(\partial_{temp}.\)~\(\qed\)

\subsubsection*{6. A vacant square of~\(L_1\) in the interior of~\(\partial_{temp}\)}
We prove the following property in this subsection. We recall that~\(\partial_{temp}\) is the cycle obtained in Step~\(5\) by merging the cycles~\(\partial_0^+\) and~\(\partial_0(L_1).\)\\
\((a6.1)\) There is a vacant square~\(Z_{j_0}\) belonging to the sequence~\(L_1\) \emph{all} of whose edges lie in the interior of the cycle~\(\partial_{temp}.\)\\

%//NEED TO USE PROPRTY...IF OPP EDGES IN DHO L1 THEN THREE EDGES BELONGS ETC...

\begin{figure}[tbp]
\centering
%\fbox{
\includegraphics[width=1.6in, trim= 150 340 240 170, clip=true]{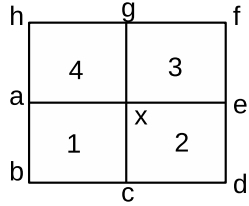}
%}
\caption{The edge~\(ax\) belongs to the vacant square~\(Z_{j_1} = abcd \in L_1\) and has both its endvertices in the interior of the merged cycle~\(\partial_{temp}.\)}
\label{z_int_fig}
\end{figure}

%CHNG CPTN ETC.. PRKVMM +ETC..

\emph{Proof of \((a6.1)\)}:  From property~\((a5.3),\) there exists an edge \(e \in \partial_0^+ \cap \partial_0(L_1)\) that has both its endvertices in the interior of the cycle~\(\partial_{temp}.\) Suppose that the edge~\(e\) belongs to the vacant square~\(Z_{j_1} \in L_1.\) We illustrate this in Figure~\ref{z_int_fig}, where the square~\(Z_{j_1}\) is the square labelled~\(1\) and the edge~\(e = ax.\)

%The edges~\(h_1(t) = ax = e,h_1(r) = xc, h_1(b) = bc\) and~\(h_1(l) = ab.\)

The edges~\(ab\) and~\(xc\) share the endvertices~\(a\) and~\(x,\) respectively, with the edge~\(e = ax.\) Since vertices~\(a\) and~\(x\) lie in the interior of the cycle~\(\partial_{temp},\) the edges~\(ab\) and~\(xc\) also lie in the interior of~\(\partial_{temp}.\) If the edge~\(bc\) also is in the interior of~\(\partial_{temp},\) then the square~\(Z_{j_1}\) has all its edges in the interior of~\(\partial_{temp}\) and the property~\((a6.1)\) is true in this case.

%We continue the exploration process for one more step and look at the ``next" edge after~\(ax\) in the cycle~\(\partial_0(L_1).\)
%In other words, exactly one of \((h_1(t),h_1(l),h_1(b))\) or \((h_1(t),h_1(r),h_1(b))\) forms consecutive edges in the cycle \(\partial_0(L_1).\)\\

We now assume that the edge~\(bc\) does not lie in the interior of~\(\partial_{temp}\) and proceed with the argument. We recall that the edge~\(e = ax\) belongs to~\(\partial_0^+ \cap \partial_0(L_1) \subset \partial_0(L_1).\) We use the following additional properties.\\
\((b6.1)\) The edge~\(bc\) belongs to the cycles~\(\partial_{temp}\) and~\(\partial_0(L_1).\) In general, if an edge~\(e\) belongs to a vacant square of the sequence~\(L_1\) and does not lie in the interior of~\(\partial_{temp},\) then~\(e\) belongs to both the cycles~\(\partial_{temp}\) and~\(\partial_0(L_1).\)\\
\((b6.2)\) Either~\((ax,ab,bc)\) or~\((ax,xc,bc)\) forms a subpath in the cycle~\(\partial_0(L_1)\) but not both. If~\((ax,ab,bc)\) is a subpath of~\(\partial_0(L_1),\) then the edge~\(cx\) lies in the interior of~\(\partial_0(L_1).\)\\
Suppose that~\((ax,ab,bc)\) forms a subpath of~\(\partial_0(L_1).\) We recall that the cycle~\(\partial_{temp}\) is the union of two bridges~\(A_{out} \subset \partial_0^+\) and~\(B_{out} \subset \partial_0(L_1).\) \\
\((b6.3)\) The vertex~\(b\) is an endvertex of the bridge~\(B_{out} \subset \partial_{temp}\) and therefore~\(b\) is also an endvertex for the path~\(B_{in} = \partial_0(L_1) \setminus B_{out}.\)\\
\((b6.4)\) The first two edges in the path~\(B_{in} = \partial_0(L_1) \setminus B_{out}\) are~\(ab\) and~\(ax.\)\\

\begin{figure}[ht!]
    \centering
    %\fbox{
    \begin{subfigure}[t]{0.4\textwidth}
    \centering
        \includegraphics[width=2.5in, trim= 130 330 150 240, clip=true]{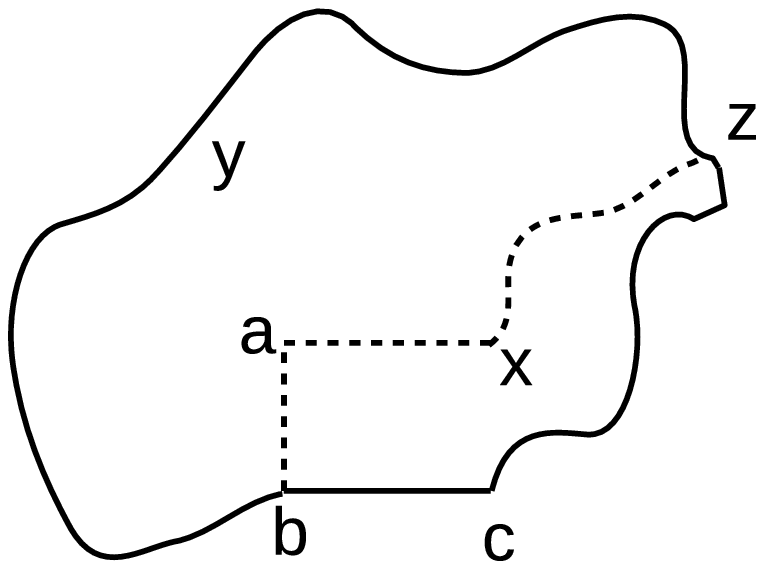}
        \caption{}
     \end{subfigure} %}%
     ~
     %\fbox{ %to see the bounding boxes
     \begin{subfigure}[t]{0.4\textwidth}
     \centering
        \includegraphics[width=2.5in, trim= 140 270 130 280, clip=true]{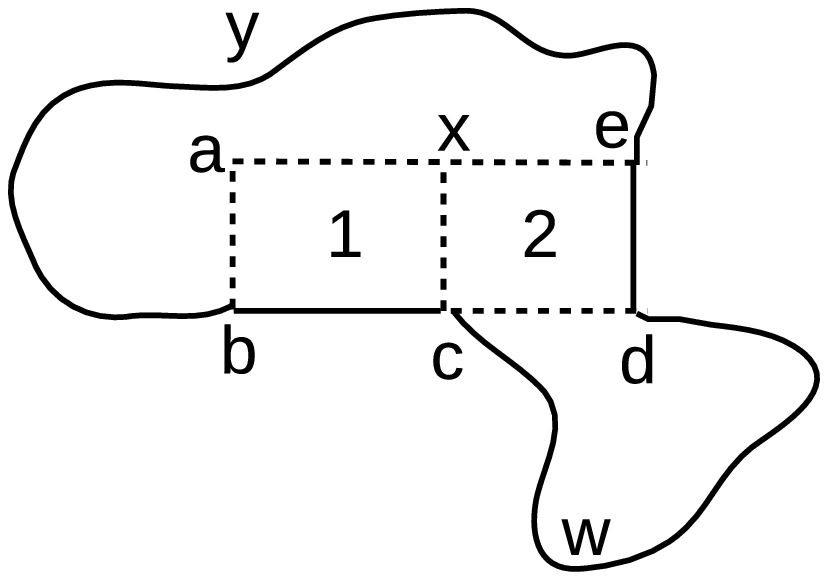}
        \caption{}
    \end{subfigure}
    %}

\caption{\((a)\) The cycle~\(\partial_{temp}\) represented by the cycle~\(byzcb.\) The bridge~\(B_{out} = bcz\) and the path~\(B_{in}\) represented by the dotted curve~\(baxz\) lies entirely in the interior of~\(\partial_{temp}.\) \((b)\) If the edge~\(ed\) belongs to the cycle~\(\partial_{temp} = bcwdeyb,\) then the path~\(B_{in} = baxe\) contains only edges of two distinct vacant squares of the sequence~\(L_1,\) represented by labels~\(1\) and~\(2\) (see proof of Property~\((c6.2)\) below). The bridge~\(B_{out}\) is represented by the path~\(bcwde.\)}
\label{z_int_fig2}
\end{figure}

The above properties are illustrated in Figure~\ref{z_int_fig2}\((a),\) where the merged cycle~\(\partial_{temp} = bczyb.\) The bridge~\(B_{out} = bcz\) and the path~\(B_{in} = baxz\) contains the two edges~\(ab\) and~\(ax.\) We use property~\((b6.4)\) along with the fact that the path~\(B_{in}\) contains edges from three distinct vacant squares of the sequence~\(L_1\) (see property~\((a5.3)\) of Step~\(5\)) to arrive at a contradiction.\\\\
\emph{Proof of \((b6.1)-(b6.2)\)}: To see~\((b6.1)\) is true, we argue as follows. The edge~\(bc\) belongs to the vacant square~\(Z_{j_1}\) of the sequence~\(L_1\) and the square~\(Z_{j_1}\) lies in the interior of the cycle~\(\partial_{temp}\) (property~\((a5.1)\) of Step~\(5\)). Since the edge~\(bc\) does not lie in the interior of~\(\partial_{temp},\) we have that~\(bc\) belongs to~\(\partial_{temp}.\)

%We recall from Step~\(5\) that the cycle~\(\partial_{temp}\) is the union of two bridges~\(A_{out} \subset \partial_0^+\) and~\(B_{out} \subset \partial_0(L_1).\) From property~\((a5.1)\) of Step~\(5\) regarding the cycle~\(\partial_{temp},\) we therefore have that the edge~\(bc\) belongs to either~\(\partial_0^+\) or~\(\partial_0(L_1)\) but not both.

Suppose now that the edge~\(bc\) does not belong to the outermost boundary cycle~\(\partial_0(L_1).\) The edge~\(bc\) belongs to the vacant square~\(Z_{j_1} \in L_1\) and the square~\(Z_{j_1}\) lies in the interior of the outermost boundary cycle~\(\partial_0(L_1)\) (see property~\((a4.1)\) of Step~\(4\)). Thus the edge~\(bc\) lies in the interior of~\(\partial_0(L_1)\) and so lies in the interior of the bigger cycle~\(\partial_{temp}\) (property~\((a5.1)\)). This is a contradiction since~\(bc \in \partial_{temp}\) and so the edge~\(bc\) belongs to~\(\partial_0(L_1).\) This proves the first statement of~\((b6.1)\) and the proof of the second statement is analogous.

We prove~\((b6.2)\) as follows. Using property~\((b6.1)\) we have that the opposite edges~\(ax\) and~\(bc\) of the vacant square~\(Z_{j_1} \in L_1\) belong to the outermost boundary cycle~\(\partial_0(L_1).\) From property~\((a4.3)\) in Step~\(4\) regarding opposite edges of the cycle \(\partial_0(L_1),\) we then have that either~\((ax,ab,bc)\) or~\((ax,xc,bc)\) forms a subpath in the cycle~\(\partial_0(L_1).\)  If both the subpaths belong to \(\partial_0(L_1),\) then the sequence~\(L_1\) contains the single square~\(Z_{j_1}\) in its interior. This contradicts the property \((a4.1)\) of~\(\partial_0(L_1)\) which states that the sequence~\(L_1\) has at least three distinct vacant squares. This proves the first statement of~\((b6.2).\)

Suppose now that the path~\((ax,ab,bc)\) is a subpath of~\(\partial_0(L_1).\) The edge~\(cx\) does not belong to~\(\partial_0(L_1)\) and since the square~\(Z_{j_1}\) containing the edge~\(cx\) lies in the interior of the cycle~\(\partial_0(L_1),\) (see property~\((a4.1)\) of~\(\partial_0(L_1)\) in Step~\(4\)) we have that the edge~\(cx\) also lies in the interior of~\(\partial_0(L_1).\)

%In what follows we assume that \((h_1(t),h_1(l),h_1(b))\) form consecutive edges in \(\partial_0(L_1).\) We then continue the exploration process and look at the next edge \(h_{next}\) after \(h_1(t)\) in the cycle \(\partial_0(L_1)\) containing the endvertex \(x_2.\) An analogous analysis holds with endvertex~\(x_1\) for the other case where \((h_1(t),h_1(r),h_1(b))\) form consecutive edges in \(\partial_0(L_1).\)

To prove~\((b6.3),\) we argue as follows. Using the fact that~\(bc \in \partial_{temp}\) (property~\((b6.1)\)), there are two edges in the cycle~\(\partial_{temp}\) containing the vertex~\(b,\) one of which is the edge~\(bc.\) Let \(g \in \partial_{temp}\) be the other edge containing~\(b\) as an endvertex. The cycle~\(\partial_{temp}\) is obtained by merging the cycles~\(\partial_0^+\) and~\(\partial_0(L_1)\) and so the edge~\(g\) therefore belongs to either~\(\partial_0^+\) or~\(\partial_0(L_1)\) (property~\((a5.1)\) of~\(\partial_{temp}\)). The edge~\(bc\) belongs to~\(\partial_{0}(L_1)\) (property~\((b6.1)\)) and by assumption, the edge~\(ab\) also belongs to~\(\partial_0(L_1).\) So the edge~\(g\) cannot belong to the cycle~\(\partial_0(L_1)\) and therefore belongs to~\(\partial_0^+.\) This also means that the vertex~\(b\) is an endvertex of the bridge~\(B_{out} \subset \partial_0(L_1).\) This proves~\((b6.3).\)

We prove~\((b6.4)\) as follows. By assumption, the edges~\(ab\) and~\(ax\) are consecutive edges of the cycle~\(\partial_0(L_1)\) and also, the vertices~\(a\) and~\(x\) lie in the interior of the merged cycle~\(\partial_{temp}.\) This means that both~\(ax\) and~\(ab\) lie in the interior of~\(\partial_{temp}\) and therefore necessarily belong to the path~\(B_{in} = \partial_0(L_1) \setminus B_{out}.\)

From property~\((b6.3),\) the vertex~\(b\) is an endvertex of the bridge~\(B_{out}\) and is therefore also the endvertex of the path \(B_{in}\) i.e., the paths~\(B_{out}\) and~\(B_{in}\) meet at~\(b.\) Therefore the edges~\(ab\) and~\(ax\) are the first two edges of~\(B_{in}.\)  This proves~\((b6.4).\)~\(\qed\)

%// CHNG PRTY IV TO PROTY V EVERY WHERE....

%//SEE ABV + BELOW...

From property~\((b6.4),\) we have that the edges~\(ab\) and~\(ax\) form consecutive edges of the path~\(B_{in}.\) Let \(h_{next}\) be the ``next" edge in~\(B_{in};\) i.e., the edge in~\(B_{in} \setminus \{ax\}\) containing~\(x\) as an endvertex. The right edge~\(cx\) of the square~\(Z_{j_1}\) lies in the interior of~\(\partial_0(L_1)\) (property~\((b6.2)\)) and so cannot belong to~\(B_{in} \subset \partial_0(L_1).\) Thus~\(h_{next} \in \{gx, xe\}\) and we discuss each case separately.\\\\
\underline{\emph{Case~\((i)\)}}: The edge~\(h_{next} = xe.\) We then have the following properties.\\
\((c6.1)\) The square labelled~\(2\) in Figure~\ref{z_int_fig} is vacant and belongs to the sequence~\(L_1.\)\\
We call the vacant square labelled~\(2\) as~\(Z_{j_2} \in L_1.\)\\
\((c6.2)\) All four edges of~\(Z_{j_2}\) lie in the interior of the cycle~\(\partial_{temp}.\)\\

%We recall that our aim is to find a vacant square in~\(L_1\) that has all its edges in the interior of the merged cycle~\(\partial_{temp}.\) We have the following properties.\\

%\((c6.2)\) The edges~\(xe\) and~\(cx\) of the square~\(Z_{j_2}\) lie in the interior of the cycle~\(\partial_{temp}.\)\\
%If the edges~\(cd\) and \(ed\) of~\(Z_{j_2}\) also lie in the interior of~\(\partial_{temp},\) then~\(Z_{j_2}\) is the desired vacant square.\\
%\((c6.3)\) If the edge~\(cd\) does not lie in the interior of~\(\partial_{temp},\) then~\(cd \in \partial_{temp}.\)

\emph{Proof of~\((c6.1)-(c6.2)\)}: To see~\((c6.1),\) we argue as follows. Since~\((ax,ab,bc)\) is a subpath of~\(\partial_0(L_1),\) the edge~\(cx\) lies in the interior of~\(\partial_0(L_1)\) (see property~\((b6.2)\)). This means that the square labelled~\(2\) in Figure~\ref{z_int_fig} also lies in the interior of~\(\partial_0(L_1).\) In this Case~\((i),\) the edge~\(xe \in B_{in} \subset \partial_0(L_1)\) also belongs to~\(\partial_0(L_1)\) and so the square labelled~\(3\) lies in the exterior of~\(\partial_0(L_1).\) Applying Theorem~\ref{thm_out}, property~\((v)\) to the edge~\(xe \in \partial_0(L_1),\) we obtain the the square labelled~\(2\) belongs to the sequence~\(L_1\) and is therefore vacant.

We prove~\((c6.2)\) as follows. By assumption, the edge~\(xe\) belongs to the path~\(B_{in}\) and from property~\((a5.2),\) every edge in~\(B_{in}\) lies in the interior of the merged cycle~\(\partial_{temp}.\) In particular, so does~\(xe.\) From property~\((b6.2),\) the edge~\(cx\) lies in the interior of the cycle~\(\partial_0(L_1)\) and the merged cycle~\(\partial_{temp}\) contains~\(\partial_0(L_1)\) in its interior (property~\((a5.1)\) of Step~\(5\)) and so the edge~\(cx\) also lies in the interior of~\(\partial_{temp}.\)

Suppose that the edge~\(cd\) does not lie in the interior of~\(\partial_{temp}.\) We arrive at a contradiction as follows. The edge~\(cd\) belongs to the vacant square~\(Z_{j_2} \in L_1\) and so from property~\((b6.1),\) we have that~\(cd\) belongs to the cycles~\(\partial_{temp}\) and~\(\partial_0(L_1).\) By assumption, the edge~\(xe \in B_{in} \subset \partial_0(L_1).\) In particular, this means that the opposite edges~\(xe\) and~\(cd\) both belong to~\(\partial_0(L_1).\) From property~\((a4.3)\) of~\(\partial_0(L_1)\) in Step~\(4,\) this means that either~\((xe,cx,cd)\) or~\((xe, de, cd)\) forms a subpath of~\(\partial_0(L_1).\) But from property~\((b6.2),\) we have that the edge~\(cx\) lies in the interior of  the cycle~\(\partial_0(L_1)\) and so~\(cd \notin \partial_0(L_1).\) This means \(de \in \partial_0(L_1)\) and the edges~\((xe,de,cd)\) form a subpath of~\(\partial_0(L_1).\)

We recall that we have assumed that~\((ax,ab,bc)\) already forms a subpath of~\(\partial_0(L_1).\) From the final statement in the above paragraph, we then have that~\(\partial_0(L_1) = (ax,xe,ed,dc,cb,ba);\) i.e., the cycle~\(\partial_0(L_1) = axedcba\) and there are exactly two squares in the interior of~\(\partial_0(L_1).\)  But this is a contradiction since the outermost boundary cycle~\(\partial_0(L_1)\) contains all squares of the sequence~\(L_1\) in its interior and the sequence~\(L_1\) contains at least three distinct vacant squares (property~\((a4.1),\) Step~\(4\)). Thus the edge~\(cd\) also lies in the interior of the cycle~\(\partial_{temp}.\)

%this means that~\(\partial_0(L_1)\) contains only two vacant squares in its interior, a contradiction to the property~\((a4.1),\) Step~\(4,\) that~\(L_1\) contains at least three distinct vacant squares. Thus~\(cd\) lies in the interior of~\(\partial_{temp}.\)

Finally, suppose that the edge~\(ed\) does not lie in the interior of the cycle~\(\partial_{temp}.\) Using property~\((b6.1)\) again, we have that~\(ed\) belongs to the cycles~\(\partial_{temp}\) and~\(\partial_0(L_1).\) Since~\(\partial_0(L_1)\) is the union of the path~\(B_{in} \cup B_{out},\) the edge~\(ed\) belongs to either~\(B_{in}\) or~\(B_{out}.\) But from property~\((a5.2),\) every edge of~\(B_{in}\) lies in the interior of the cycle~\(\partial_{temp}.\) This necessarily means that the edge~\(ed\) belongs to the bridge~\(B_{out} \subset \partial_0(L_1).\) But the edge~\(xe\) belongs to the path~\(B_{in}\) and so~\(e\) is an endvertex of the path~\(B_{in}\) (and~\(B_{out}\)).

%// SE ABV +eTC..

From property~\((b6.3),\) we already have that vertex~\(b\) is also an endvertex of~\(B_{in}.\) From~\((b6.4),\) we also have that the first two edges of~\(B_{in}\) are~\(ab\) and~\(ax.\) Thus we must have~\(B_{in} = (ab,ax,xe).\) We refer to Figure~\ref{z_int_fig2}\((b)\) for illustration. The path~\(B_{in} = baxe\) and the bridge~\(B_{out} = bcwde.\) The merged cycle~\(\partial_{temp}\) is the union of~\(B_{out}\) and the wavy path~\(bye.\) The square labelled~\(2\) is the vacant square~\(Z_{j_2}\) and as proved above, has its edges~\(cx, xe\) and~\(cd\) in the interior of~\(\partial_{temp}.\)

From the above paragraph, we obtain that the path~\(B_{in}\) contains only edges of two distinct vacant squares of the sequence~\(L_1.\) In Figure~\ref{z_int_fig2}, the vacant squares are represented by~\(axcb\) and~\(cxed.\) This is a contradiction to property~\((a5.3)\) of Step~\(5\) which states that the path~\(B_{in}\) contains edges of at least three distinct squares of the sequence~\(L_1.\) So the edge~\(ed\) also lies in the interior of~\(\partial_{temp}.\)~\(\qed\)

We now consider the case \((ii).\)\\
\underline{\emph{Case~\((ii)\)}}: The edge~\(h_{next} = xg.\) We then have the following properties.\\
\((d6.1)\) The square labelled~\(3\) in Figure~\ref{z_int_fig} is vacant and belongs to the sequence~\(L_1.\)\\
We call the vacant square labelled~\(3\) as~\(Z_{j_3} \in L_1.\)\\
\((d6.2)\) All four edges of~\(Z_{j_3}\) lie in the interior of the cycle~\(\partial_{temp}.\)\\

\emph{Proof of~\((d6.1)-(d6.2)\)}: To see~\((d6.1),\) we argue as follows. We recall that the edge~\(ax \in \partial_0(L_1) \cap \partial_0^+\) and the square labelled~\(4\) is occupied and the square~\(Z_{j_1}\) labelled~\(1\) is vacant and belongs to the sequence~\(L_1.\) Since~\(gx \in B_{in} \subset \partial_0(L_1),\) the square labelled~\(3\) in Figure~\ref{z_int_fig} is also vacant and belongs to the sequence~\(L_1.\) This is seen by applying Theorem~\ref{thm_out}, property~\((v)\) to the edge~\(gx.\) This proves~\((d6.1).\)

To prove~\((d6.2),\) we argue as follows. The edge~\(gx \in B_{in}\) lies in the interior of~\(\partial_{temp}\) (property~\((a5.2)\)).

To prove that the edge~\(xe\) lies in the interior of the cycle~\(\partial_{temp},\) we argue as follows. The edge~\(gx\) also belongs to the outermost boundary cycle~\(\partial_0^+\) because if not, then the square~\(Z_{j_3} \in L_1\) would lie in the interior of~\(\partial_0^+\) contradicting property~\((b4.2).\) From the discussion in the previous paragraphs, we therefore have that both the edges~\(ax\) and~\(gx\) belong to~\(\partial_0^+ \cap \partial_0(L_1).\) Thus the edge~\(ex \notin \partial_0^+ \cup \partial_0(L_1)\) and therefore does not belong to~\(\partial_{temp}.\) Since the square~\(Z_{j_3}\) containing~\(ex\) as an edge lies in the interior of~\(\partial_{temp},\) the edge~\(xe\) also lies in the interior of~\(\partial_{temp}.\)

If the edge~\(gf \in \partial_{temp},\) then it belongs to either~\(\partial_0^+\) or~\(\partial_0(L_1)\) but not both~(property~\((a5.1)\)). If it belongs to~\(\partial_0^+,\) then the square above~\(gf\) is occupied and since~\(Z_{j_3} \in L_1\) (property~\((d6.1)\)), this would mean that~\(gf\) also belongs to~\(\partial_0(L_1).\) Because if~\(gf\) lies in the interior of~\(\partial_0(L_1),\) then the (occupied) square above~\(gf\) would also lie in the interior of~\(\partial_0(L_1),\) a contradiction to property~\((b4.1).\) Thus~\(gf \in \partial_0(L_1)\) and since by assumption the edge~\(gf\) also belongs to~\(\partial_0^+,\) this contradicts the first statement of this paragraph.

From the above paragraph, we have that~\(gf \notin \partial_0^+.\) If~\(gf \in \partial_0(L_1),\) then~\(gf \in \partial_0(L_1) \setminus \partial_0^+\) and so necessarily belongs to~\(B_{out}.\) This is because using properties~\((a5.1)\) and~\((a5.2)\) we obtain that every edge in~\(\partial_0(L_1) \cap \partial_0^+\) belongs to~\(B_{in} \subset \partial_0(L_1).\) Since the edge~\(gx\) belongs to~\(B_{in},\) the vertex~\(g\) is an endvertex of~\(B_{out}\) and~\(B_{in}\) and so~\(B_{in}\) in this case is the path~\((ba,ax,xg).\) This is a contradiction since~\(B_{in}\) contains edges from at least three distinct squares of the sequence~\(L_1\) (property~\((a5.3)\)) and the path~\((ba,ax,xg)\) contains only edges from the squares~\(Z_{j_1} \in L_1\) (labelled~\(1\) in Figure~\ref{z_int_fig}) and the square~\(Z_{j_3} \in L_1\) (labelled~\(3\) in Figure~\ref{z_int_fig}).

From the discussion in the above paragraph, the edge~\(gf\) belongs to the interior of~\(\partial_{temp}.\) We then finally consider the edge~\(fe.\) If~\(fe \in \partial_{temp},\) then arguing as in the case of the edge~\(gf\) above, we have that~\(fe \in \partial_0(L_1) \setminus \partial_0^+\) and so~\(fe \in B_{out} \subset \partial_0(L_1).\) We recall that we have assumed that~\(gx \in B_{in} \subset \partial_0(L_1)\) in this case and so in particular the opposite sides~\(fe\) and~\(gx\) of the square~\(Z_{j_3} \in L_1\) belong to~\(\partial_0(L_1).\)

Using the opposite sides property~(property~\(a4.3\)), the edge~\(gf\) or~\(xe\) also must belong to~\(\partial_0(L_1).\) But~\(xe\) cannot belong to~\(\partial_0(L_1)\) since both the edges~\(ax\) and~\(gx\) belong to~\(B_{in} \subset \partial_0(L_1).\) This~\(gf \in \partial_0(L_1)\) and arguing as in two paragraphs above, we get that~\(gf \in B_{out}\) and that~\(B_{in} = (ba,ax,xg).\) Again using property~\((a5.3)\) that~\(B_{in}\) contains edges of at least three distinct squares in the sequence~\(L_1,\) we obtain a contradiction. Thus~\(fe \notin \partial_{temp}\) and since the square~\(Z_{j_3} \in L_1\) containing~\(fe\) as an edge lies in the interior of~\(\partial_{temp},\) the edge~\(fe\) also lies in the interior of~\(\partial_{temp}.\)~\(\qed\)

%// TO SEE CRFLLY +etC...

\subsubsection*{7. Arriving at a contradiction}
We use the above properties of the cycles \(\partial_0(L_1)\) and \(\partial_{temp}\) to arrive at a contradiction. We recall the iterative merging procedure in Step \(1\) of the proof. Following the same procedure, we merge the vacant squares \(\{Z_{i_0}\}_{i_0 \leq i \leq i_1}\) in the sequence \(L_1\) one by one in an iterative fashion with the outermost boundary cycle~\(\partial_0^+.\) We start with \(Z_{i_0}\) merge with \(\partial_0^+\) to get a new cycle. We then merge~\(Z_{i_0+1}\) with the new cycle and so on. Let \(R_1,R_2,\ldots,R_T\) denote the intermediate cycles obtained with \(R_1 = \partial_0^+\) and \(R_T\) denoting the final cycle.

The cycles \(\{R_i\}_{1 \leq i \leq T}\) satisfy properties \((x1)-(x3)\) (see Step 1) with \(\{Y_l\}_{1 \leq l \leq i}\) replaced with \(\{Z_l\}_{i_0 \leq l \leq i}.\) The final cycle \(R_T\) therefore satisfies the following properties.\\
\((a7.1)\) Every occupied square of \(C^+(0)\) is contained in the interior of \(R_T.\) Every vacant square in the sequence \(L_1\) is contained in the interior of \(R_T.\)\\
\((a7.2)\) If edge \(e \in R_T,\) then either \(e \in \partial_0^+\) or \(e\) belongs to some vacant square in the sequence \(L_1.\)\\
\((a7.3)\) The cycle~\(R_T\) is unique in the sense that if some cycle \(C\) satisfies properties \((a7.1)-(a7.2)\) above, then \(C = R_T.\)\\

\emph{Proof of \((a7.1)-(a7.3)\)}: Properties \((a7.1)-(a7.2)\) follow from properties \((a1.1)-(a1.3)\) of Step~\(1.\) It suffices to prove \((a7.3).\) If \(C \neq R_T,\) suppose there exists an edge \(e\) of \(C\) in the exterior of \(R_T.\) From property \((a7.2),\) we have that \(e \in \partial_0^+\) or \(e\) belongs to a vacant square~\(Z_j \in L_1.\) This means that either some occupied square of \(C^+(0)\) or some vacant square of~\(L_1\) lies in the exterior of \(R_T,\) a contradiction to property \((a7.1).\) An analogous argument holds if~\(R_T\) contains an edge in the exterior of~\(C.\)\(\qed\)

From property~\((a5.1)\) of the cycle \(\partial_{temp}\) in Step~\((5)\) above, we have that~\(\partial_{temp}\) satisfies properties~\((a7.1)-(a7.2)\) and so~\(\partial_{temp} = R_T.\) We recall that~\(\Lambda^+\) denotes the set of all vacant squares sharing an edge with some occupied square of the component~\(C^+(0)\) and that all the squares in the sequence~\(L_1 \subset L \subset \Lambda^+.\) If we now continue the iterative procedure starting with~\(R_T = \partial_{temp}\) and merge the vacant squares in~\(\Lambda^+ \setminus L_1\) one by one, then arguing as above, we have that the final cycle obtained~\(R_{fin}\) satisfies properties~\((b1.1)-(b1.3)\) of the cycle~\(D_{fin}\) mentioned in Step~\(1.\) Thus using the uniqueness property~\((b1.4),\) we have that~\(R_{fin} = D_{fin}.\) %CHK CRFLLY HERE +eTC....

The interior of~\(\partial_{temp}\) is contained in the interior of~\(R_{fin} = D_{fin}\) and so by property~\((a6.1),\) the vacant square~\(Z_{j_0}\) of the sequence~\(L_1 \subset L\) has all its edges in the interior of~\(D_{fin}.\) But by construction, every vacant square in the original sequence \(L = (Z_1,\ldots,Z_s)\) contains an edge in \(D_{fin}\) (see property~\((a2.2)\) in Step~\(2\) above). This is a contradiction and so the assumption we made in Step~\(4\) that~\(L_1 \neq L\) is false. Thus the sequence \(L\) is itself a \(S-\)cycle. This completes the proof of Theorem~\ref{thm5}.~\(\qed\)

\renewcommand{\theequation}{\thesection.\arabic{equation}}

\subsection*{Acknowledgement}
I thank Professors Rahul Roy and Federico Camia for crucial comments and for my fellowships. I also thank NISER for my fellowship.

\bibliographystyle{plain}

\begin{thebibliography}{10}
%\bibitem{Azuma} K. Azuma. (1967).
%\newblock {Weighted Sums of Certain Dependent Random Variables}.
%\newblock {\em Tohoku Math. Journ.}, \textbf{19}, 357--367.


%\bibitem{boll} N. H. Bingham, C. M. Goldie and J. L. Teugels. (1989).
%\newblock {\em Regular Variation}.
%\newblock {Cambridge University Press}.

%\bibitem{cox} J. T. Cox and R. Durrett. (1981).
%\newblock {Some Limit Theorems for Percolation Processes with Necessary and Sufficient Conditions}.
%\newblock {\em Annals Prob.}, \textbf{9}, 583--603.

\bibitem{boll} B. Bollobas. (2001).
\newblock {\em Modern Graph Theory}.
\newblock {Springer}.

\bibitem{boll} B. Bollobas and O. Riordan. (2006).
\newblock {\em Percolation}.
\newblock {Academic Press}.

\bibitem{gane} G. Ganesan. (2017).
\newblock {Duality in pecolation via outermost boundaries I: Bond percolation}.
\newblock {\emph{Arxiv Link}:  https://arxiv.org/abs/1704.00461}.


\bibitem{penrose} M. Penrose. (2003).
\newblock {\em Random Geometric Graphs}.
\newblock {Oxford}.


%\bibitem{galam} J. Galambos. (1978).
%\newblock {\em Asymptotic Theory of Extreme Order Statistics}.
%\newblock {Wiley}.


%\bibitem{ganesh} A. J. Ganesh, L. Massoulie, and D. F. Towsley. (2005).
%\newblock {The effect of network topology on the spread of epidemics}.
%\newblock {\em Proc. IEEE Infocomm}, pp. 1455--1466.

\end{thebibliography}

\end{document}